\documentclass[10pt]{article}

\usepackage[utf8]{inputenc}
\usepackage{cmap}
\usepackage[T1]{fontenc}
\input{glyphtounicode}
\usepackage{lmodern}
\usepackage{pdflscape}
\usepackage[hmargin=1in,vmargin=0.75in]{geometry}

\usepackage{amsmath,amssymb,mathtools}
\usepackage{graphicx}
\usepackage{booktabs}
\usepackage{longtable}
\usepackage{caption}
\usepackage[compact]{titlesec}
\usepackage{array}
\usepackage{ragged2e,placeins}
\usepackage{xcolor}
\usepackage[expansion=false]{microtype}
\usepackage[authoryear,round]{natbib}
\usepackage{tikz}
\usetikzlibrary{arrows.meta,positioning,shapes.geometric,fit,backgrounds,calc}
\usepackage[colorlinks=true,linkcolor=blue!50!black,citecolor=blue!50!black,urlcolor=blue!50!black]{hyperref}
\usepackage{enumitem}
\AtBeginDocument{\providecommand{\doi}{}\renewcommand{\doi}[1]{\href{https://doi.org/\detokenize{#1}}{\mbox{doi:\detokenize{#1}}}}}

\newlength{\tablewidth}
\newcommand{\SurveyEligibilityDate}{October 12, 2026}

\newcommand{\SurveyVerificationDate}{October 12, 2026}

\newcommand{\SurveyTitle}{Quantum-Classical Hybrid Methods for Optimization}
\newcommand{\SurveySubtitle}{Lessons from Classical Decomposition}
\newcommand{\SurveyAuthorName}{David E. Bernal Neira}
\newcommand{\SurveyDepartment}{Davidson School of Chemical Engineering}
\newcommand{\SurveyInstitution}{Purdue University}
\newcommand{\SurveyCity}{West Lafayette}
\newcommand{\SurveyState}{IN}
\newcommand{\SurveyCountry}{USA}
\newcommand{\SurveyEmail}{dbernaln@purdue.edu}
\newcommand{\SurveyTutorialNote}{This paper develops and extends a tutorial delivered at the CPAIOR 2024 Master Class ``Quantum Computing for CP, AI, and OR, and vice-versa'' \citep{bernal2024cpaior} and is prepared for the associated special issue of \emph{Constraints}.}
\newcommand{\software}[1]{\texttt{#1}}

\newcommand{\SurveyAbstract}{%
Quantum and physics-inspired solvers provide specialized sampling interfaces for Ising and quadratic unconstrained binary optimization models, but constrained optimization requires broader coordination and certification.
This survey uses classical decomposition to analyze how such solvers can be integrated into hybrid quantum-classical algorithms.
We distinguish formal decomposition, reformulation and reduction, and computational orchestration.
For each interface, we ask what the transformation preserves, which assumptions make the outer method converge, and which certificates the implemented inner routine supplies.
The central questions are whether that information meets the classical algorithm's requirements and whether hybridization improves solution quality or total cost enough to justify its coordination overhead.
Benders decomposition, column generation, graph reduction, and copositive methods illustrate this distinction by separating heuristic proposals from classical certification.
We also distinguish component effects within a fixed wrapper, reduction--solver interactions, and comparisons between complete quantum-enabled and classical systems.
The analyzed studies provide conditional component benefits and favorable co-design evidence, but do not establish transferable end-to-end superiority over strong classical comparators.
The resulting framework identifies opportunities for constraint programming and operations research in structure-preserving design, safeguarded hybrid algorithms, and controlled benchmarking.%
}
\newcommand{\SurveyKeywords}{quantum optimization; hybrid quantum-classical algorithms; decomposition; constraint programming; certification; benchmarking.}

\title{\bfseries \SurveyTitle:\\[2pt]
\large \SurveySubtitle}
\author{\SurveyAuthorName\thanks{\SurveyDepartment, \SurveyInstitution, \SurveyCity, \SurveyState, \SurveyCountry.
Email: \texttt{\SurveyEmail}.
\SurveyTutorialNote}}

\date{}
\hypersetup{
  pdftitle={Quantum-Classical Hybrid Methods for Optimization: Lessons from Classical Decomposition},
  pdfauthor={David E. Bernal Neira},
  pdfsubject={Survey of hybrid quantum-classical optimization and the attribution of reported improvements},
  pdfkeywords={\SurveyKeywords},
  pdflang={en-US}
}

\begin{document}
\maketitle
\begin{abstract}
\noindent\SurveyAbstract

\noindent\textbf{Keywords:} \SurveyKeywords
\end{abstract}
%
%
\section{Introduction}\label{sec:introduction}

Classical optimization makes discrete problems tractable by exploiting structure, coordinating tasks that are easier for suitable solvers, and combining their results.
Quantum optimization adds specialized solvers that encode an objective in a physical system or a sequence of controlled quantum operations, then return sampled candidate assignments.
Section~\ref{sec:quantum-background} introduces these quantum terms for readers without a quantum information science background.
Using them for general constrained models requires additional encoding, classical coordination, or both.

Decomposition connects specialized solvers to broader algorithms through interfaces such as pricing, recourse evaluation, and separation.
The output required at an interface may be a useful candidate, a valid inequality, or a certified bound.
We ask which outputs a quantum routine can supply, which classical checks remain necessary, and whether the resulting hybrid improves quality or total computational cost, including coordination overhead, over suitable classical alternatives.

Existing reviews describe Ising machines \citep{mohseni2022ising}, the quantum approximate optimization algorithm (QAOA) and its variants \citep{blekos2024review}, the wider opportunities and obstacles in quantum optimization \citep{abbas2024challenges}, and, in a recent preprint, learned control of adaptive quantum workflows \citep{lau2026beyond}.
Our contribution is to analyze these topics jointly through \emph{interface requirements, guarantees, and performance attribution}.

Tables~\ref{tab:hybrid-roles} and~\ref{tab:guarantee-requirements} apply this analysis to the surveyed algorithm families, distinguishing returned candidates from the certificates required for guarantees.
Section~\ref{sec:performance-attribution} then assesses the reported studies through component comparisons, reduction--solver interactions, and complete-system comparisons.
The resulting conclusions concern the purposive corpus described in Appendix~\ref{app:survey-method}, with eligibility through \SurveyEligibilityDate{}.

Section~\ref{sec:classical-foundations} introduces the classical foundations; Section~\ref{sec:solver-interfaces} maps mixed-integer models to Ising and QUBO form and introduces the quantum algorithms and hardware that accept them.
Section~\ref{sec:hybrid-methods} develops the hybrid architectures, Section~\ref{sec:performance-attribution} evaluates their evidence, and Section~\ref{sec:outlook} identifies design questions for constraint programming and operations research.


%
\section{Discrete nonlinear optimization and the classical success story}\label{sec:classical-foundations}

\subsection{What can be expressed, and why it is hard}\label{sec:discrete-nonlinear-models}

Mixed-integer nonlinear programming (MINLP) is a very expressive modeling language.
A general MINLP has the form

\begin{equation}\label{eq:minlp}
\begin{aligned}
\underset{x,y}{\operatorname{minimize}}\quad &f(x,y)\\
\operatorname{subject\ to}\quad &g(x,y)\le0,\\
&x\in\mathbb{R}^{n_c},\quad y\in\mathbb{Z}^{n_i},
\end{aligned}
\end{equation}

with \(f\) and \(g\) possibly nonlinear and nonconvex.
The \(n_c\) continuous variables represent quantities, the \(n_i\) integer variables represent choices, and nonlinear terms may represent physical relationships.
Applications include process design, scheduling, and engineering optimization \citep{biegler2004retrospective}.
The convex-MINLP solver review \citep{kronqvist2019convex} and the fifty-year perspective on MINLP and disjunctive programming \citep{kronqvist2026fifty} trace the algorithms and formulations behind this modeling capability.
Even linear integer programming is NP-hard; practical methods exploit instance structure, whereas quantum Ising/quadratic unconstrained binary optimization (QUBO) interfaces require additional encodings to represent the constrained form above.

This survey focuses on a single objective.
Quantum multiobjective methods can instead explore trade-offs by approximately solving or sampling scalarized problems with different objective weights \citep{kotil2025multiobjective}.
Many interface and certification considerations below apply to those subproblems, although solving individual scalarizations does not by itself certify a complete Pareto front.

\subsection{Decomposition as the engine of progress}\label{sec:classical-decomposition}

%
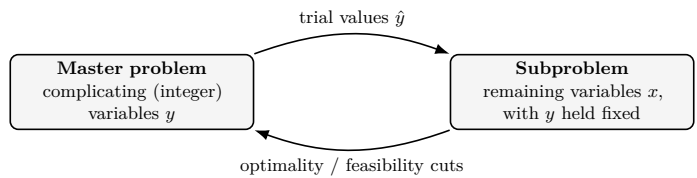
\begin{figure}[tbp]
\centering
\scalebox{0.76}{\begin{tikzpicture}[
  box/.style={rectangle, rounded corners, draw=black, thick, align=center, minimum height=1.3cm, text width=4.0cm, fill=gray!8},
  every node/.style={font=\small}
]
  \node[box] (master) {\textbf{Master problem}\\ complicating (integer)\\ variables $y$};
  \node[box, right=3.4cm of master] (sub) {\textbf{Subproblem}\\ remaining variables $x$,\\ with $y$ held fixed};
  \draw[-{Latex[length=2.5mm]}, thick] (master.north east) to[bend left=20]
     node[above, midway]{trial values $\hat{y}$} (sub.north west);
  \draw[-{Latex[length=2.5mm]}, thick] (sub.south west) to[bend left=20]
     node[below, midway]{optimality / feasibility cuts} (master.south east);
\end{tikzpicture}}
\caption{Classical Benders coordination: the master proposes complicating variables \(y\); classical recourse in \(x\) returns valid cuts.
Hybrid interfaces are compared in Section~\ref{sec:hybrid-methods}.}\label{fig:benders-coordination}
\end{figure}


Classical solver studies report large aggregate gains from algorithms and hardware \citep{bixby2012brief,koch2022progress}; they do not isolate decomposition's causal share.
Decomposition is one important way to exploit problem structure.
The idea: split a problem into smaller or more tractable parts, then coordinate them.
In Benders decomposition, fixing complicating variables leaves a linear subproblem whose dual generates cuts for a master problem \citep{benders1962partitioning}.
Generalized Benders extends this construction to nonlinear subproblems through convex duality, making it relevant to MINLP \citep{geoffrion1972generalized}.
Dantzig--Wolfe decomposition reformulates a structured linear program over the extreme points and, when necessary, extreme rays of its blocks, pricing columns as needed \citep{dantzig1960decomposition}; this underlies column generation and branch-and-price \citep{barnhart1998branch,lubbecke2005column}.
For convex MINLP under standard subproblem and regularity assumptions, outer approximation alternates between a mixed-integer linear master and continuous nonlinear subproblems, accumulating valid supporting linearizations until the bounds meet \citep{duran1986outer}.
For nonconvex MINLP, local linearizations need not give valid global bounds; deterministic global methods combine valid relaxations with domain partitioning and bound tightening to certify a solution \citep{burer2012nonconvex}.
In Figure~\ref{fig:benders-coordination}, the outgoing assignment fixes the complicating variables \(y\); the return cuts exclude infeasible master choices or bound their recourse cost.
These cuts improve the master's relaxation, whereas a candidate assignment alone supplies no such bound.
Other decompositions use different master variables and oracle interfaces, as summarized in Section~\ref{sec:hybrid-methods}.

Logic-based Benders provides a bridge to constraint programming: cuts follow from an inference dual, allowing subproblems whose proofs need not come from linear programming duality \citep{hooker2003logic}.
Surveys and textbooks develop the theory and implementation of these methods \citep{rahmaniani2017benders,conejo2006decomposition,constante2025relaxation,conforti2014integer}.
Their interfaces require more than candidate assignments: cut validity, bounds, and termination each have proof obligations.

We adapt a classical decomposition taxonomy \citep{constante2025relaxation} into three composable structural descriptions:
\begin{itemize}
\item \emph{Formal decomposition} splits variables, constraints, or the search domain.
Its duality-based forms include Benders \citep{benders1962partitioning}, Dantzig--Wolfe \citep{dantzig1960decomposition}, and augmented Lagrangian methods such as the alternating direction method of multipliers (ADMM) \citep{boyd2011distributed}; its recursive forms include branching with bound-based pruning, as in branch and bound \citep{land1960automatic}.
\item \emph{Reformulation and reduction} transform the problem through penalty mappings, conic lifts, coarsening, or elimination.
An \emph{exact reformulation} preserves the relevant solutions and objective; an \emph{approximation} may change them.
A \emph{relaxation} supplies a lower bound in minimization (upper in maximization) by enlarging the feasible region or underestimating (overestimating) the objective.
\item \emph{Orchestration} coordinates routines through samples, parameters, or candidate solutions.
For example, an outer optimizer can tune the parameters of a quantum routine using the objective values estimated from its outputs; this parameter-tuning loop can itself be nested inside a decomposition.
\end{itemize}

Section~\ref{sec:hybrid-methods} uses these categories to distinguish formal master--subproblem methods, graph and conic reformulations, and sampling or parameter-tuning orchestration.
They identify algorithmic structure, not mutually exclusive guarantee classes.
An exact conic reformulation can support cutting-plane methods \citep{brown2024copositive} or Frank--Wolfe updates \citep{frank1956algorithm,yurtsever2022qfw}, but heuristic inner solves require separate convergence justification.
Conversely, a formal decomposition may lack its classical guarantee when convexity, oracle accuracy, or certified pruning assumptions fail.
We therefore assess three questions independently: what the transformation preserves, which assumptions make the outer method converge, and which certificates the implemented inner routine supplies.

The practical lesson is to assign each solver a task that exploits structure and whose output can be verified.
Smaller subproblems can help, but duplicated work, many iterations, or lost global structure can erase that benefit.
Dimension reduction alone establishes neither computational savings nor a potential quantum contribution.
\par


%
\section{Quantum and physics-inspired solvers for Ising and QUBO problems}\label{sec:solver-interfaces}

This section first states the classical models that Ising and quantum solvers accept and how mixed-integer problems can be translated into them (Sections~\ref{sec:ising-qubo} and~\ref{sec:encoding-costs}).
It then introduces the quantum concepts used in the rest of the article (Section~\ref{sec:quantum-background}), the main quantum optimization algorithms (Section~\ref{sec:quantum-algorithms}), and the hardware and physics-inspired platforms that run them (Section~\ref{sec:analog-solvers}).

\subsection{The Ising and QUBO substrate}\label{sec:ising-qubo}

Here the Ising objective is the energy of a finite, classical, pairwise spin model on a graph \(G=(V,\mathcal{E})\), with vertices \(V\) and edges \(\mathcal{E}\):

\begin{equation}\label{eq:ising-energy}
E(s)=\sum_{i\in V}h_i s_i+\sum_{\{i,j\}\in\mathcal{E}}J_{ij}s_i s_j,
\qquad s_i\in\{-1,1\}.
\end{equation}

The real coefficients \(h_i\) are external fields, or single-spin biases, and \(J_{ij}\) are pairwise couplings; their signs are included in the coefficients.
Fixing spins in a set \(V_{\mathrm{fix}}\subseteq V\) leaves each remaining spin \(i\) with the \emph{effective local field} \(h_i+\sum_{j\in V_{\mathrm{fix}}:\{i,j\}\in\mathcal{E}}J_{ij}s_j\), its new linear coefficient.
Under the displayed plus-sign convention, \(J_{ij}<0\) favors aligned spins (a ferromagnetic coupling), whereas \(J_{ij}>0\) favors opposite spins (an antiferromagnetic coupling).
The graph and coefficients are problem data, with no restriction to a regular interaction graph or couplings of one sign.
This specifies the optimization model within the broader physics tradition of Ising ferromagnets and disordered spin glasses \citep{ising1925ferromagnetismus,sherrington1975solvable}.
Finding a ground-state configuration means minimizing \(E(s)\).
The affine substitution \(s_i=2x_i-1\), with inverse \(x_i=(s_i+1)/2\), maps spin assignments bijectively to bits and converts the energy into an equivalent quadratic unconstrained binary optimization (QUBO) problem

\begin{equation}\label{eq:qubo}
\min_{x\in\{0,1\}^n}\quad x^\top Qx,
\end{equation}

where \(n=|V|\) and \(Q\in\mathbb{R}^{n\times n}\) is symmetric.
Off-diagonal entries are half the coefficients of the equivalent upper-triangular convention because \(x^\top Qx\) counts each pair twice.
Linear terms occupy the diagonal, since \(x_i^2=x_i\); additive constants do not change the minimizer and are dropped.
Conversely, substituting \(x_i=(s_i+1)/2\) into \(x^\top Qx\) gives~\eqref{eq:ising-energy}, up to a constant, with \(J_{ij}=Q_{ij}/2\) for \(i<j\) and \(h_i=\frac12\sum_j Q_{ij}\).
The maximum-cut problem partitions vertices into two sets to maximize the weight of crossing edges.
With nonnegative edge weights \(w_{ij}\), it maximizes \(\sum_{\{i,j\}\in\mathcal{E}}w_{ij}(1-s_i s_j)/2\), since an edge contributes exactly when its spins differ \citep{goemans1995improved}.
Negation gives, up to a constant, a zero-field Ising minimization with \(J_{ij}=w_{ij}/2\); the affine substitution above produces a QUBO.
A physical solver returns samples from a device- and protocol-dependent distribution, typically without an optimality certificate.
Ising-oriented devices include quantum annealers, optical and analog machines, and digital accelerators, with applications in optimization and engineering \citep{mohseni2022ising,bernal2022perspectives}; Section~\ref{sec:analog-solvers} describes these platforms.

\subsection{From mixed-integer programs to QUBO}\label{sec:encoding-costs}

%
\begin{figure}[tbp]
\centering
\scalebox{0.76}{\begin{tikzpicture}[
  b/.style={rectangle, rounded corners, draw=black, thick, align=center, minimum height=1.9cm, text width=3.3cm},
  every node/.style={font=\small}
]
  \node[b, fill=blue!6] (c) {\textbf{Constrained problem}\\[2pt] $\min\; f(x)$\\ s.t. $Ax \le b$\\ $x \in \{0,1\}^{n}$};
  \node[b, fill=orange!10, right=2.2cm of c] (p) {\textbf{Penalized form}\\[2pt] $\min\; f(x)$\\ $+\,\rho\,\lVert Ax + Sw - b\rVert^{2}$\\ slack bits $w$, $u = Sw$};
  \node[b, fill=red!8, right=2.2cm of p] (q) {\textbf{QUBO}\\[2pt] $\min\; (x,w)^{\top} Q\, (x,w)$\\ additional bits,\\ coefficient range};
  \draw[-{Latex[length=2.5mm]}, thick] (c) -- node[above, align=center, font=\scriptsize]{constraints\\into objective} (p);
  \draw[-{Latex[length=2.5mm]}, thick] (p) -- node[above, align=center, font=\scriptsize]{absorb slacks,\\tune $\rho$} (q);
\end{tikzpicture}}
\caption{Constraint-to-QUBO mapping for integer or scaled-rational data and a quadratic objective (after any required quadratization).
Exact slack encoding \(u=Sw\), with binary slack bits \(w\) and encoding matrix \(S\), requires bounded ranges commensurate with the encoding; arbitrary real coefficients need not admit it.
Here \(Q\) acts on the concatenated bits \((x,w)\).
Slack bits are auxiliaries; squaring the penalty introduces pairwise couplings, and weights satisfying (\ref{eq:exact-penalty}) preserve optima.
The larger model and coefficient range can worsen conditioning and increase the selected solver's resource requirements.}\label{fig:constraint-to-qubo}
\end{figure}
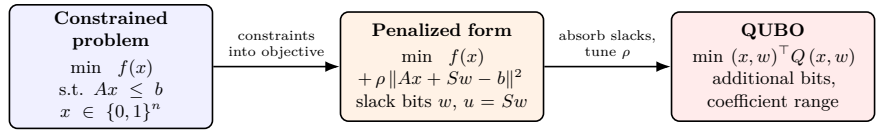


The MINLP~\eqref{eq:minlp} is the natural starting point for these mappings: the task is to represent general integer, continuous, and nonlinear structure, not only to penalize constraints.
Standard catalogs and modeling tutorials describe the resulting formulations \citep{lucas2014ising,glover2019qubo}.
Reaching~\eqref{eq:qubo} involves four steps, each with its own cost.
Below, \(x\) denotes the encoded binary vector.
First, each integer variable is written in binary: for \(z\in\{L,\ldots,U\}\), a binary encoding is \(z=L+\sum_{k=0}^{m-1}2^k b_k\), with \(m=\lceil\log_2(U-L+1)\rceil\) and \(b_k\in\{0,1\}\).
If the bit range exceeds \(U-L\), the constraint \(\sum_k2^k b_k\le U-L\) or a bounded encoding must exclude excess values.
Second, continuous variables must either be discretized, which introduces an additional approximation, or remain with a classical solver; the latter is the decomposition route of Section~\ref{sec:hybrid-methods}.
Third, nonlinear terms must become quadratic.
Once the variables are binary, every real-valued function of them has an exact multilinear polynomial representation, for example \(e^{b}=1+(e-1)\,b\) for \(b\in\{0,1\}\), although the representation can need exponentially many terms; auxiliary variables that represent products, enforced by penalty terms, then reduce it to quadratic form \citep{boros2002pseudo}.
Approximating a function to keep the model small is a separate modeling choice, distinct from the approximation introduced by discretizing continuous variables.
Fourth, constraints can become objective penalties, but there is no universal encoding.
A linear equality can use \(\rho(a^\top x-b)^2\); an inequality can first introduce a nonnegative slack, whose exact positional binary encoding requires bounded, commensurate values, as in integer or scaled-rational data.
Reading Figure~\ref{fig:constraint-to-qubo} from left to right reveals two resource costs: slack bits enlarge the assignment, and squaring a constraint introduces pairwise couplings between its variables.
A compact constrained model can therefore become a denser QUBO with a wider coefficient range.
Direct penalties can avoid slack bits when the constraint has suitable structure: for binary variables, \(\rho x_i x_j\) penalizes exactly the violation of \(x_i+x_j\le1\).
One-hot encodings assign one bit to each allowed integer value and require exactly one selected bit; they use more bits than positional encodings but can simplify selection constraints.
Higher-order violation penalties may instead require quadratization auxiliaries \citep{lucas2014ising,glover2019qubo,quintero2022characterization}.
Eliminating variables, preserving feasibility in the quantum algorithm (Section~\ref{sec:quantum-algorithms}), or retaining constraints in a classical master (Section~\ref{sec:hybrid-methods}) offers alternatives to a single penalized QUBO.
These choices trade variable count, coupling density, coefficient range, and proof requirements; large penalty weights stretch the range a physical solver must resolve, and samples can remain infeasible.
Penalization changes the feasible domain, even when it preserves global optima \citep{constante2025relaxation}.

Let \(\mathcal{X}\subseteq\{0,1\}^n\) be the finite encoded domain, let \(F\subseteq\mathcal{X}\) be the nonempty feasible set, let \(P(x)\ge0\) be a penalty function vanishing exactly on \(F\), and let \(f^\star=\min_{x\in F}f(x)\).
Then, for \(\mathcal{X}\setminus F\ne\varnothing\), the penalized problem \(\min_{x\in\mathcal{X}}f(x)+\rho\,P(x)\) has the same global minimizers and optimal value as the constrained problem whenever

\begin{equation}\label{eq:exact-penalty}
\rho>\max\left\{0,\max_{x\in\mathcal{X}\setminus F}
\frac{f^\star-f(x)}{P(x)}\right\}.
\end{equation}

The threshold is finite but involves the unknown \(f^\star\); an arbitrary penalty weight need not preserve the optimum.
If \(\mathcal{X}=F\), any \(\rho\ge0\) suffices.
One sufficient, possibly conservative bound is \(\rho>(\overline f-\underline f)/\delta\), where \(\underline f\le f(x)\) on \(\mathcal{X}\), \(\overline f\) is the value of a known feasible solution, and every infeasible point has penalty at least \(\delta>0\).
Problem-specific formulations can give sharper weights and smaller encodings; catalogs discuss many examples \citep{lucas2014ising,glover2019qubo}, and the maximum \(k\)-colorable subgraph analysis explicitly characterizes penalty choices and auxiliary-variable trade-offs \citep{quintero2022characterization}.
The broader theory of Lagrangian reformulation distinguishes equality of optimal values from equality of minimizer sets and gives sufficient conditions for finite exact penalties \citep{quintero2024lagrangian}.
In particular, a positive minimum penalty over infeasible points is automatic on the finite domain used in~\eqref{eq:exact-penalty}, whereas continuous domains can require limiting penalty sequences rather than a finite exact weight.

Maximum independent set (MIS) illustrates the mapping.
Given a graph \(G=(V,\mathcal{E})\), MIS asks for as many vertices as possible with no edge between two selected vertices; as a mixed-integer program (MIP),
\[
\max_{x\in\{0,1\}^{V}}\ \sum_{v\in V}x_v
\qquad\text{s.t.}\qquad x_u+x_v\le1\quad\text{for all }\{u,v\}\in\mathcal{E}.
\]
Penalizing each edge with both endpoints selected gives the QUBO
\[
\min_{x\in\{0,1\}^{V}}\ -\sum_{v\in V}x_v+\rho\sum_{\{u,v\}\in\mathcal{E}}x_ux_v.
\]
Any \(\rho>1\) satisfies~\eqref{eq:exact-penalty}, because deleting one endpoint of each such edge yields an independent set that is smaller by at most the number of penalized edges.
Arrays of individually controlled neutral atoms can instead realize MIS natively on unit-disk graphs; the interaction mechanism and its restrictions are explained in Section~\ref{sec:analog-solvers}.

Modeling software automates parts of this translation and differs in what it delivers.
\software{QUBO.jl} supports mathematical-programming reformulations and solver access \citep{xavier2026qubojl}; \software{Qiskit Optimization} converts quadratic programs into QUBO form when they have no continuous variables, no quadratic constraints, and integer-valued coefficients in all linear constraints \citep{qiskit_optimization_docs}; and D-Wave's open-source \software{dimod} library represents binary and constrained quadratic models, converting a constrained model into a penalized binary model only when its constraints are linear and its integer variables have zero lower bounds \citep{dwave_dimod_docs}.
\software{PyQUBO} supports symbolic QUBO construction and constraint-aware decoding \citep{zaman2022pyqubo}, and \software{qubovert} manages polynomial and constrained Boolean or spin formulations \citep{qubovert_repo}.
\software{Fixstars Amplify SDK} combines modeling and conversion with access to several quantum and classical solver backends \citep{amplify_sdk_docs}.
Such libraries build and transform models locally, whereas hosted solver services, discussed in Section~\ref{sec:decomposition-interfaces}, also execute them.
The interface automates parts of the translation; its choice alone establishes neither exactness of a penalty model nor a quantum contribution.

\subsection{Quantum computation for optimization readers}\label{sec:quantum-background}

A classical bit takes value 0 or 1; a qubit is a quantum system with two distinguished states, written \(|0\rangle\) and \(|1\rangle\).
For \(n\) qubits, the \emph{computational basis} consists of \(2^n\) mutually orthogonal vectors \(|x\rangle\), one for each binary assignment \(x\); the ket notation \(|\cdot\rangle\) denotes a state vector \citep{nielsen2010quantum}.
An ideal pure state is a normalized superposition \(|\psi\rangle=\sum_x\alpha_x|x\rangle\), where the complex amplitudes satisfy \(\sum_x|\alpha_x|^2=1\).
Measuring in this basis returns one assignment \(x\) with probability \(|\alpha_x|^2\), not all assignments at once.
Repeated preparation and measurement produce samples.
Unlike classical probabilities, amplitudes have phases: when combined, they can reinforce or cancel, producing \emph{interference} that changes the sampling probabilities.
A \emph{product state} factors into individual-qubit states; a pure state that does not is \emph{entangled}.

A \emph{Hamiltonian} \(H\) is a Hermitian matrix (equal to its conjugate transpose) representing energy, with real eigenvalues.
For a classical cost \(f(x)\), the diagonal problem Hamiltonian satisfies \(H_C|x\rangle=f(x)|x\rangle\); its lowest-energy, or ground, states encode minimizing assignments.

Quantum evolution is the change of the amplitudes over time.
For constant \(H\), ideal closed-system evolution applies \(U(t)=e^{-\mathrm{i}tH}\), where \(\mathrm{i}^2=-1\), setting the reduced Planck constant \(\hbar=1\); this linear map is unitary, meaning reversible and norm-preserving.
Optimization algorithms choose these operations to make good assignments more likely upon measurement; doing so efficiently is the algorithmic challenge.
Analog devices physically implement a Hamiltonian and vary its coefficients continuously.
Figure~\ref{fig:adiabatic-sweep} illustrates how energy levels change along such a path; Section~\ref{sec:quantum-algorithms} explains the annealing method.
Gate-model computers compose \emph{quantum gates}, unitary operations on one or a few qubits, into a \emph{circuit}, followed here by measurement.
A designed Hamiltonian evolution is implemented or approximated by such a gate sequence.
The upper part of Figure~\ref{fig:qaoa-loop} illustrates this circuit notation: each horizontal wire carries one qubit through operations ordered from left to right, and the final measurements produce the bits of one assignment.
The gates marked \(\mathsf H\) prepare \((|0\rangle+|1\rangle)/\sqrt{2}\) from each \(|0\rangle\); these Hadamard gates are distinct from the Hamiltonian \(H\).
Circuit width counts qubits; gate depth counts sequential layers of operations that cannot run in parallel.
Gate platforms include superconducting circuits, trapped ions, photonics \citep{bharti2022nisq}, and neutral atoms \citep{henriet2020neutral}.
Noise and control error perturb the ideal operations; fault-tolerant computation uses many physical qubits to encode logical qubits and correct errors, at substantial resource cost.
Progress from the noisy intermediate-scale quantum (NISQ) regime, whose gates lack error-correction protection \citep{preskill2018nisq} toward early fault tolerance now includes demonstrations of error correction and logical operations, while large-scale fault-tolerant systems remain a development goal \citep{bluvstein2026fault}.

\subsection{Quantum optimization algorithms}\label{sec:quantum-algorithms}

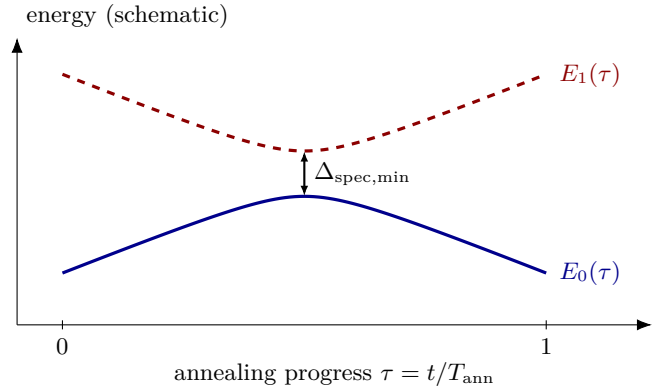
\begin{figure}[htbp]
\centering
\begin{tikzpicture}[every node/.style={font=\small}]
  \draw[-{Latex[length=2mm]}] (0,0) -- (8.4,0);
  \node at (4.2,-0.65) {annealing progress $\tau=t/T_{\mathrm{ann}}$};
  \draw[-{Latex[length=2mm]}] (0,0) -- (0,3.8)
    node[above right] {energy (schematic)};
  \draw[very thick, blue!55!black, domain=0.6:7, samples=101]
    plot (\x,{2-sqrt(((\x-3.8)/2.5)^2+0.09)});
  \draw[very thick, red!55!black, dashed, domain=0.6:7, samples=101]
    plot (\x,{2+sqrt(((\x-3.8)/2.5)^2+0.09)});
  \draw[{Latex[length=1.4mm]}-{Latex[length=1.4mm]}, thick]
    (3.8,1.7) -- (3.8,2.3) node[midway, right] {$\Delta_{\mathrm{spec},\min}$};
  \node[blue!55!black, anchor=west] at (7.05,0.69) {$E_0(\tau)$};
  \node[red!55!black, anchor=west] at (7.05,3.31) {$E_1(\tau)$};
  \draw (0.6,0.06) -- (0.6,-0.06) node[below] {$0$};
  \draw (7,0.06) -- (7,-0.06) node[below] {$1$};
\end{tikzpicture}
\caption{Illustrative instantaneous energy levels along an annealing path.
Only the ground level \(E_0\) (solid) and first excited level \(E_1\) (dashed) are shown; their closest approach marks the minimum spectral gap \(\Delta_{\mathrm{spec},\min}\).
These curves depict a schematic spectrum, not measured device performance.
TikZ code prepared with OpenAI Codex assistance.}\label{fig:adiabatic-sweep}
\end{figure}

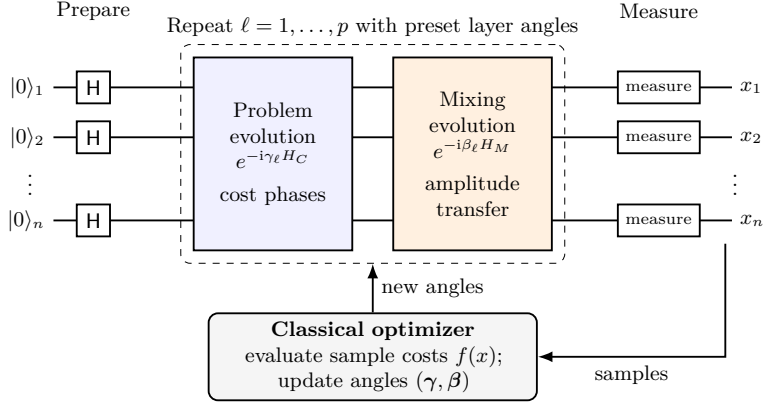
\begin{figure}[htbp]
\centering
\scalebox{0.88}{\begin{tikzpicture}[
  block/.style={rectangle, draw=black, thick, align=center,
                minimum height=2.9cm, text width=2.15cm},
  gate/.style={rectangle, draw=black, thick, fill=white,
               minimum height=0.48cm, minimum width=0.5cm},
  every node/.style={font=\small},
  flow/.style={-{Latex[length=2.2mm]}, thick}
]
  \foreach \y/\i in {1/1,0.25/2,-1/n} {
    \draw[thick] (-0.3,\y) node[left] {$|0\rangle_{\i}$} -- (9.9,\y);
    \node[gate] at (0.3,\y) {$\mathsf H$};
    \node[gate, minimum width=0.9cm] at (8.8,\y) {\scriptsize measure};
    \node[anchor=west] at (9.9,\y) {$x_{\i}$};
  }
  \node at (-0.65,-0.36) {$\vdots$};
  \node at (9.95,-0.36) {$\vdots$};
  \node[align=center] at (0.3,2.15) {Prepare};
  \node[align=center] at (8.8,2.15) {Measure};
  \node[block, fill=blue!6] (cost) at (3,0)
    {Problem evolution\\ $e^{-\mathrm{i}\gamma_\ell H_C}$\\[3pt] cost phases};
  \node[block, fill=orange!12] (mix) at (6,0)
    {Mixing evolution\\ $e^{-\mathrm{i}\beta_\ell H_M}$\\[3pt] amplitude transfer};
  \node[draw, dashed, rounded corners, inner sep=0.18cm,
        fit=(cost)(mix), label=above:{Repeat $\ell=1,\ldots,p$ with preset layer angles}] (layers) {};
  \node[draw, rounded corners, thick, align=center, text width=4.7cm,
        fill=gray!8] (classical) at (4.5,-3.05)
    {\textbf{Classical optimizer}\\ evaluate sample costs $f(x)$;\\
     update angles $(\boldsymbol\gamma,\boldsymbol\beta)$};
  \draw[flow] (9.8,-1.35) -- (9.8,-3.05) -- node[below]{samples} (classical.east);
  \draw[flow] (classical.north) -- node[right]{new angles} (layers.south);
\end{tikzpicture}}
\caption{QAOA as a quantum circuit and a classical parameter loop.
Horizontal wires represent qubits, with omitted wires indicated by dots; operations act from left to right.
Hadamard gates (\(\mathsf H\)) prepare the initial superposition.
The wide evolution blocks denote compiled gate sequences, and measurement occurs only after all \(p\) problem--mixer pairs.
Every sampled assignment requires a fresh circuit execution.
TikZ code prepared with OpenAI Codex assistance.}\label{fig:qaoa-loop}
\end{figure}

\paragraph{Adiabatic computation and quantum annealing.}
Adiabatic quantum computation starts from an easily prepared ground state and slowly changes the Hamiltonian, aiming to follow its evolving ground state \citep{albash2018adiabatic}.
Quantum annealing commonly uses \(H(t)=A(t)H_D+B(t)H_C\), decreasing the driver weight \(A(t)\) while increasing the problem weight \(B(t)\) \citep{kadowaki1998quantum}.
The usual driver \(H_D=-\sum_i X_i\) is a \emph{transverse field}: \(X_i\) exchanges \(|0\rangle\) and \(|1\rangle\) on qubit \(i\), so this term couples assignments differing in one bit.
Its ground state is the uniform superposition of assignments, also prepared by the Hadamard gates in the circuit example of Figure~\ref{fig:qaoa-loop}.
``Transverse'' refers to its orientation perpendicular to the spin axis used for the diagonal Ising objective; its evolution redistributes amplitudes between assignments.
Figure~\ref{fig:adiabatic-sweep} plots the two lowest instantaneous energies against progress \(\tau=t/T_{\mathrm{ann}}\in[0,1]\), where \(T_{\mathrm{ann}}\) is the annealing duration.
For the nondegenerate levels shown, their vertical separation \(\Delta_{\mathrm{spec}}(\tau)=E_1(\tau)-E_0(\tau)\) is the spectral gap; it concerns the evolving Hamiltonian, rather than an incumbent's objective gap to optimality.
The adiabatic strategy is to remain in the lower level as its associated state changes.
Near the marked minimum gap, transitions out of that state can be difficult to suppress; a slower schedule may be required, with the sufficient runtime also depending on how the Hamiltonian changes and the allowed error \citep{albash2018adiabatic}.
Thus the minimum gap alone neither specifies a runtime nor establishes a computational advantage.
Finite temperature, control error, and restricted connectivity further limit the transfer of closed-system guarantees to implemented annealers \citep{hauke2020perspectives}.
Implemented annealing is therefore assessed as a sampling heuristic, with analog evolution realized continuously in hardware.
%


\paragraph{Variational algorithms.}
In the circuit-based variational algorithms considered here, a classical optimizer adjusts parameters of a gate-based circuit.
For each parameter choice, repeated executions estimate a training objective, often the expected problem cost \citep{cerezo2021vqa,bharti2022nisq}.
An \emph{ansatz} is the chosen family of states generated by the parameterized circuit; choosing it fixes the search space.
The variational quantum eigensolver (VQE) estimates the expected energy of a parameterized state and uses a classical optimizer to lower it, seeking the Hamiltonian's smallest eigenvalue and an associated state \citep{peruzzo2014variational}.
For a cost Hamiltonian diagonal in the computational basis, this energy is estimated by the mean objective value of the sampled assignments.
The quantum approximate optimization algorithm (QAOA) typically starts from an equal-amplitude superposition of all assignments and applies \(p\) layers, each a problem evolution followed by a mixing evolution \citep{farhi2014quantum,zhou2020qaoa}.
At layer \(\ell\), the problem step \(e^{-\mathrm{i}\gamma_\ell H_C}\) changes each assignment's phase according to its cost; alone, it leaves sampling probabilities unchanged and is not a descent step.
The mixing step \(e^{-\mathrm{i}\beta_\ell H_M}\), commonly with \(H_M=\sum_i X_i\), transfers amplitude between assignments, allowing interference to change their probabilities.
In Figure~\ref{fig:qaoa-loop}, the dashed enclosure repeats this pair \(p\) times without intermediate measurement; each wide block must be compiled into gates, so \(p\) is not the physical gate depth.
For fixed angles \((\beta_\ell,\gamma_\ell)_{\ell=1}^{p}\), fresh executions yield assignments whose costs define the training objective; the return arrow carries these samples to the classical optimizer, which selects the next set of \(2p\) angles.
This outer loop seeks angles that improve the training objective without guaranteeing descent or optimality; it is distinct from the mixer's transitions within one circuit (Section~\ref{sec:variational-orchestration}).
Warm-start QAOA maps a classical relaxation solution to initial qubit amplitudes and adapts the mixer \citep{egger2021warmstarting}.
Recursive QAOA uses measured correlations to fix relations between variables and reruns QAOA on the reduced problem; Section~\ref{sec:recursive-reduction} treats it as a reduction method.
%


\paragraph{Constraints inside quantum algorithms.}
Constraints that are not penalized can enter the quantum algorithm itself.
The feasibility-preserving generalization of QAOA starts with amplitudes only on feasible assignments and uses mixers whose ideal transitions stay within that set \citep{hadfield2019quantum}; implementation errors can violate this invariance.
Linear Ising penalties use single-spin biases rather than extra pairwise couplings and can reduce connectivity and coefficient-range demands, with theoretical, simulation, and annealing-hardware evidence for the studied constraints \citep{mirkarimi2024linear,mirkarimi2024experimental}.
Embedding Lagrangian duality in adiabatic computation can avoid slack-variable encodings; in a fault-tolerant circuit implementation, the proposed construction has a qubit count independent of the number of constraints, while solution-sampling probabilities still depend on the duality gap and complementary slackness \citep{gabbassov2025lagrangian}.
Variational methods can also keep multipliers or slack variables classical.
A transaction-settlement method updates circuit parameters and nonnegative slack variables, with sample-informed slack adjustments intended to escape poor local solutions \citep{braine2021transaction}.
A variational quantum eigensolver with constraints (VQEC) uses quantum estimates of costs, constraints, and gradients while a classical primal--dual method updates circuit parameters and Lagrange multipliers \citep{le2023vqec}.
Its expectation-constrained formulations can admit distributions containing infeasible assignments, and its nonconvex parameter updates have no general global-convergence guarantee; individual samples still require the original feasibility checks.
Quantum Zeno constraint enforcement instead repeatedly measures whether the state lies in the subspace spanned by feasible assignments, without identifying the full assignment.
Frequent checks suppress transitions out of that subspace in the ideal limit; finite measurement rates and implementation errors must be assessed in practice \citep{herman2023zeno,berwald2025zeno}.
Thus constraints can enter through objective penalties, dual multipliers, or restrictions on the quantum evolution.
Linear penalty constructions need not satisfy the assumptions of~(\ref{eq:exact-penalty}); their feasibility and optimum-preservation properties require separate, problem-specific analysis.
These operations can also appear within a decomposition.
Their mathematical exactness and the implemented solver's guarantees must be assessed separately.

\paragraph{Approximation guarantees and evidence for advantage.}
Approximation guarantees remain problem-specific.
For maximum-cut on unweighted 3-regular graphs, one-layer (\(p=1\)) QAOA guarantees an expected approximation ratio (cut weight divided by the optimum) of about 0.6924 \citep{farhi2014quantum}; classical semidefinite methods guarantee about 0.87856 on general graphs with nonnegative edge weights, including the unweighted case \citep{goemans1995improved}, and at least 0.9326 for unweighted graphs of maximum degree three \citep{halperin2004maxcut}.
E3LIN2 (Max-3XOR) maximizes satisfied three-variable parity equations; the early one-layer QAOA guarantee for its bounded-occurrence variant \citep{farhi2014bounded} stimulated a stronger classical guarantee \citep{barak2015beating}.
Numerical tests on random Max-\(k\)XOR, its \(k\)-variable generalization, find a classical mean-field method matching or outperforming QAOA on average in the tested regimes \citep{muller2025limitations}.
Neither the preceding approximation guarantees nor these numerical observations establish implemented system advantages.
More favorable evidence comes from several distinct regimes.
For low-autocorrelation binary sequences, noiseless simulations found better empirical scaling than the tested exact classical solvers, together with small hardware demonstrations \citep{shaydulin2024evidence}.
For random 8-SAT (Boolean satisfiability with eight literals per clause) near the satisfiability threshold, a fault-tolerant resource analysis combines QAOA with amplitude amplification, which increases the probability of measuring a satisfying assignment, and projects a runtime crossover under explicit assumptions on error correction and classical computational resources \citep{omanakuttan2025threshold}.
For maximum-cut, classically warm-started, regularized QAOA has favorable hardware quality comparisons under matched algorithmic restrictions and a projected fault-tolerant crossover against stronger unrestricted classical heuristics \citep{he2026regularized}.
A further analysis gives better cut-fraction guarantees (fraction of edges crossing the partition) than the compared classical local algorithms for maximum \(k\)-cut on specified regular graphs with girth (shortest cycle length) at least \(2p+2\), where \(p\) counts QAOA layers and vertices are partitioned into \(k\) sets to maximize crossing edges \citep{apte2026provable}.
These studies contribute empirical scaling evidence, resource projections, and comparisons of provable guarantees; their assumptions and the distinction from demonstrated end-to-end advantage must remain explicit.

\paragraph{Fault-tolerant primitives and decoded quantum interferometry.}
Fault-tolerant primitives offer prospective guarantees: resource estimates for optimization heuristics \citep{sanders2020compilation} and near-quadratic branch-and-bound speedups under search-tree and oracle-access assumptions \citep{montanaro2020quantum}.
Decoded quantum interferometry (DQI) uses interference and an associated classical decoding problem to bias sampling toward good solutions.
Its provable separation concerns folded optimal polynomial intersection (OPI), an approximate optimization problem over a finite field \(\mathbb{F}_q\) with \(q\) elements \citep{kramer2026provable}.
Polynomials of degree less than \(K\) are represented by coefficient vectors \(c\in\mathbb{F}_q^K\); the evaluation points are grouped into \(M\) blocks of size \(h\), and block \(i\) has an evaluation matrix \(B_i\in\mathbb{F}_q^{h\times K}\) and an acceptance set \(\mathcal{A}_i\subseteq\mathbb{F}_q^h\).
The problem is
\[
\max_{c\in\mathbb{F}_q^K}\ \frac1M\sum_{i=1}^{M}\mathbf{1}\bigl[B_i c\in\mathcal{A}_i\bigr].
\]
Here \(\mathbf{1}[\cdot]\) is the zero--one indicator, and \(B_i c\) uses field arithmetic.
The separation is not claimed for this general formulation.
It holds for a growing family of folded Reed--Solomon instances whose block alphabets \(\mathbb{F}_q^h\) have superpolynomially many elements, with each acceptance set drawn independently and uniformly from the subsets containing exactly half of its block alphabet and accessible only through a membership oracle \citep{kramer2026provable}.
A membership oracle tests whether a block value belongs to \(\mathcal{A}_i\); a coherent quantum query applies this test to a superposition without measuring the queried value.
For a family with code rate \(K/(Mh)\to R=0.3\), DQI makes one coherent quantum query to each block oracle and its expected score converges to \((1+\sqrt{R(2-R)})/2\approx0.857\).
Classical algorithms reach the threshold \((1+R)/2=0.65\) with linearly many classical membership queries, but exceeding it by any fixed amount with constant probability requires superpolynomially many \citep{kramer2026provable}.
Separate classical Markov-chain Monte Carlo tests attain DQI target scores on sparse maximum-XORSAT and ordinary, unfolded OPI, with exponential empirical scaling for the latter \citep{gilfuster2026mcmc}.
These experiments do not test the folded-OPI oracle family above; that separation also does not concern QUBO instances.
For the analyzed DQI construction on unweighted maximum-cut, any graph family supporting a constant asymptotic improvement over the random cut fraction \(1/2\) has enough structure for exact classical polynomial-time solution \citep{parekh2026limitations}.
This limits that proposed advantage on maximum-cut, a QUBO special case, without proving a limitation for every QUBO or every DQI variant.
Recent DQI screening criteria compare decoding assumptions with explicit classical baselines; candidate problems and finite-budget comparisons do not establish a general speedup \citep{mazumder2026discovering}.
Classical-output approximation guarantees must also account for parameter selection, preparation, measurement, decoding, and repetition \citep{hadfield2026approximation}; oracle-query results alone do not establish an implemented hybrid solver's wall-clock benefit.

\subsection{Hardware and physics-inspired platforms}\label{sec:analog-solvers}

Quantum annealers from D-Wave implement the transverse-field evolution on superconducting qubits connected by a fixed sparse coupling graph.
A problem whose graph does not match is minor-embedded: each logical variable becomes a chain of physical qubits, which consumes qubits and couplers \citep{boothby2020next}.
For example, the Chimera graph couples each qubit to at most six others, while Pegasus permits degree fifteen \citep{boothby2020next}.

Neutral-atom arrays hold individual atoms in optical traps and encode a binary choice in a ground state and a highly excited Rydberg state.
Strong interactions suppress simultaneous excitation of nearby atoms, providing a native interface for MIS on unit-disk graphs, whose vertices are atom positions and whose edges join atoms closer than a specified distance \citep{ebadi2022quantum}.
The excited atoms represent the selected vertices; finite interaction strengths and longer-range interactions affect how closely the physical evolution realizes this graph model.
Other graphs require the penalty formulation of Section~\ref{sec:encoding-costs} or additional encodings.
These devices return candidate configurations, whose objective and feasibility can be checked classically.

Superconducting circuits also implement gate-based processors at IBM and Google \citep{kim2023utility,harrigan2021qaoa}, while IonQ and Quantinuum use trapped ions \citep{wright2019benchmarking,moses2023racetrack}.
Neutral-atom platforms likewise support both analog Hamiltonian evolution and digital gates \citep{henriet2020neutral}.
Usable circuit width and depth depend on coherence, gate errors, connectivity, and compilation overhead \citep{montanez2025evaluating}.

Optical platforms span quantum and classical regimes.
Photonic quantum accelerators are being integrated into classical optimization pipelines \citep{slysz2026photonic,slysz2026corrigendum}.
QCi's Dirac-3 implements entropy computing with weak coherent light, photon counting, and electronic feedback; it optimizes discretized real variables and encodes combinatorial problems \citep{nguyen2025entropy}.
Its reported optical components are classically simulable in principle.
Other Ising solvers include optical coherent Ising machines \citep{mcmahon2016cim,honjo2021cim}, simulated bifurcation \citep{goto2019sbm}, and digital annealers \citep{aramon2019digital}; their common Ising objective does not imply common connectivity, precision, or cost.
Large-scale comparisons make these machines relevant baselines \citep{zeng2024performance}, but broad comparisons across applications and platforms remain incomplete.

The Quantum Optimization Benchmarking Library (\software{QOBLIB}) addresses this gap with ten problem classes, shared instances, and classical reference results, allowing application-level comparisons across formulations and devices \citep{koch2026qoblib}.
Application-derived circuits can also guide hardware design: a fixed-parameter optimization benchmark measures how circuit width and depth affect the persistence of a signal distinguishable from random sampling, exposing coherence and execution-rate bottlenecks across processors \citep{montanez2025evaluating}.
A complementary benchmark evaluates execution time and solution quality on optimization applications \citep{lubinski2024optimization}; parameter-setting costs can be assessed through the stochastic-solver benchmarking framework \citep{bernalneira2025benchmarking} discussed in Section~\ref{sec:performance-attribution}.
Hardware performance under these workloads remains distinct from solver advantage against capable classical algorithms on the original problem.

Table~\ref{tab:solver-interfaces} summarizes the interfaces.
Mathematical-programming baselines provide an established comparison for Ising processing units \citep{coffrin2019evaluating}.
Section~\ref{sec:performance-attribution} evaluates quality, resources, and attribution against the relevant classical alternatives.

%
\begin{table}[tbp]
\centering
\footnotesize
\setlength{\tabcolsep}{3pt}
\setlength{\tablewidth}{\dimexpr\linewidth-18pt\relax}
\begin{tabular}{@{}>{\RaggedRight\arraybackslash}p{0.225\tablewidth}>{\RaggedRight\arraybackslash}p{0.232\tablewidth}>{\RaggedRight\arraybackslash}p{0.262\tablewidth}>{\RaggedRight\arraybackslash}p{0.281\tablewidth}@{}}
\toprule
\textbf{Solver family} & \textbf{Problem class} & \textbf{Guarantees} & \textbf{Principal limitations} \\
\midrule
Quantum annealing / analog quantum devices &
Unconstrained QUBO / Ising (\ref{eq:qubo}); native MIS on restricted geometric graphs &
Heuristic; adiabatic guarantee only in idealized limits &
Qubit count and connectivity; noise and finite temperature; embedding overhead \\
Physics-inspired Ising machines &
Unconstrained QUBO / Ising &
Heuristic &
Restricted objective structure; coefficient precision; topology and coupling constraints; readout and control scaling; no broadly applicable guarantees \\
Gate-based variational (QAOA) &
Unconstrained QUBO; some constraints via mixers &
Special-case ratios; early examples matched or exceeded classically &
Circuit depth and noise; classical parameter optimization is nonconvex \\
Fault-tolerant primitives &
Broad, including search and branch and bound &
Conditional asymptotic speedups, e.g. near-quadratic branch and bound under search-tree and oracle-access assumptions &
Require fault-tolerant hardware not yet available at scale \\
\bottomrule
\end{tabular}
\caption{Quantum and physics-inspired solver families.
Near-term interfaces are narrow; fault-tolerant primitives are broader but require hardware not yet available at scale.
Section~\ref{sec:performance-attribution} assesses component benefits and system comparisons separately.}\label{tab:solver-interfaces}
\end{table}



%
\section{Hybrid quantum-classical decomposition}\label{sec:hybrid-methods}

%
\begin{figure}[tbp]
\centering
\scalebox{0.76}{\begin{tikzpicture}[
  node distance=1.1cm and 2.6cm,
  box/.style={rectangle, rounded corners, draw=black, thick, align=center, minimum height=1.2cm, minimum width=3.6cm, fill=blue!5},
  qbox/.style={rectangle, rounded corners, draw=black, thick, align=center, minimum height=1.2cm, minimum width=3.6cm, fill=orange!12},
  every node/.style={font=\small}
]
  \node[box] (master) {\textbf{Classical routine(s)}\\ coordination, verification,\\ remaining structure};
  \node[qbox, right=of master] (sub) {\textbf{Quantum / Ising routine}\\ sample low-energy states\\ of a QUBO \eqref{eq:qubo}};
  \draw[-{Latex[length=2.5mm]}, thick]
    (master.north east) to[bend left=22]
    node[above, midway, align=center, font=\footnotesize]{dispatched piece (family-dependent):\\ master, block, pricing, separation, sampling}
    (sub.north west);
  \draw[-{Latex[length=2.5mm]}, thick]
    (sub.south west) to[bend left=22]
    node[below, midway, align=center, font=\footnotesize]{returned information (family-dependent):\\ candidate assignments or family-specific oracle output}
    (master.south east);
\end{tikzpicture}}
\caption{Representative QUBO-based hybrid interface (Table~\ref{tab:hybrid-roles}).
Native independent-set sampling and constraint-preserving circuits need not use an unconstrained QUBO interface.
Classical checks turn candidates into feasible incumbents or valid witnesses.
In quantum-master Benders, classical recourse supplies cuts and a verified master relaxation supplies the lower bound; Section~\ref{sec:decomposition-interfaces} describes this allocation.}\label{fig:hybrid-interfaces}
\end{figure}
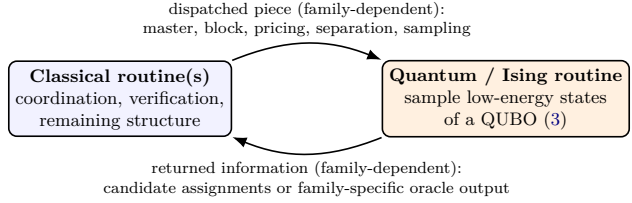


Figure~\ref{fig:hybrid-interfaces} separates the dispatched task from the returned information: a master assignment, pricing solution, block update, separation witness, or variational sample.
The return arrow need not carry a bound or optimality certificate; classical coordination must establish the family-specific property needed to continue or stop.
Tables~\ref{tab:hybrid-roles} and~\ref{tab:guarantee-requirements} summarize these interfaces and guarantee conditions.
Section~\ref{sec:performance-attribution} assesses their empirical value.
Following \citet{abbas2024challenges}, we cover quantum subproblems within classical decomposition and classical preparation, training, and postprocessing of quantum routines.

\subsection{Sampling and variational orchestration}\label{sec:variational-orchestration}

QAOA orchestration tunes the \(2p\) circuit parameters through continuous nonconvex search, typically without a global certificate \citep{farhi2014quantum,hadfield2019quantum,zhou2020qaoa,hao2025endtoend}.
The optimizer loop in Figure~\ref{fig:qaoa-loop} illustrates this potentially costly classical search \citep{nannicini2019performance}.
For suitable graph families and shallow depth, tensor-network calculations can select QAOA parameters classically before quantum sampling \citep{streif2020training}; training costs and the validity of transferring parameters between instances remain part of the assessment.
The training objective can also emphasize good samples: conditional value at risk (CVaR) averages a chosen lowest-cost fraction of outcomes \citep{barkoutsos2020cvar}.
The initial state, ansatz, and mixer (Section~\ref{sec:quantum-algorithms}) determine which assignments the circuit can reach.
A pairwise \(XY\) mixing step acts on two qubits: its evolution transfers amplitude between \(|01\rangle\) and \(|10\rangle\), leaving \(|00\rangle\) and \(|11\rangle\) unchanged.
Unlike a single-qubit \(X_i\) operation, this step preserves the number of ones and hence a fixed-cardinality constraint when starting from a feasible state \citep{hadfield2019quantum,fuchs2022constraint}.
Feasibility does not imply that every feasible assignment is reachable at a chosen depth: bounds on reachable assignments apply to fixed feasible basis states and schedules of disjoint pairwise \(XY\) mixing within fixed-cardinality variable groups \citep{ubale2026accessibility}.
For MIS (Section~\ref{sec:encoding-costs}), results establishing that suitable mixer controls can prepare an optimal state supply neither efficient parameter discovery nor an explicit circuit-depth bound \citep{tsvelikhovskiy2026expressive}.
Portfolio benchmarks exercise variational orchestration on structured instances \citep{brandhofer2022benchmarking}, while quantum-enhanced greedy optimization combines sampling with repeated variable fixing \citep{dupont2023greedy}, connecting orchestration to the recursive reductions of Section~\ref{sec:recursive-reduction}.

\paragraph{Sample-guided route construction.}
A preprint on fleet sizing with time windows constructs routes from sampled arc frequencies through classical path extraction and pruning on a directed acyclic graph \citep{makansi2025greedy}.
Its minimum-fleet guarantee requires an exact initial solution of the full discretized fleet-sizing model, an arc-frequency threshold in \((0,1)\), and multiple-longest-path extraction; heuristic sampling alone does not supply this guarantee.
This is orchestration of sampling and classical construction; it supplies routes rather than a pricing certificate.

\subsection{Graph and multilevel decomposition}\label{sec:graph-decomposition}

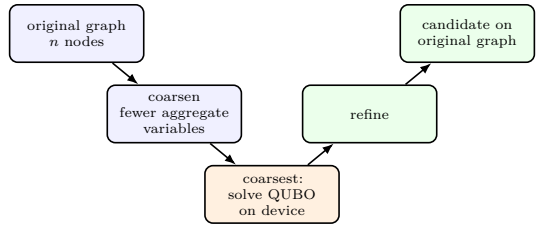
\begin{figure}[tbp]
\centering
\scalebox{0.76}{\begin{tikzpicture}[
  lvl/.style={rectangle, rounded corners, draw=black, thick, align=center, font=\scriptsize, minimum height=1.0cm, text width=2.1cm},
  every node/.style={font=\small}
]
  \node[lvl, fill=blue!6]   (a0) at (0,3)     {original graph\\ $n$ nodes};
  \node[lvl, fill=blue!6]   (a1) at (1.7,1.6) {coarsen\\ fewer aggregate\\ variables};
  \node[lvl, fill=orange!12](a2) at (3.4,0.2) {coarsest:\\ solve QUBO\\ on device};
  \node[lvl, fill=green!8]  (a3) at (5.1,1.6) {refine};
  \node[lvl, fill=green!8]  (a4) at (6.8,3)   {candidate on\\ original graph};
  \draw[-{Latex[length=2mm]}, thick] (a0)--(a1);
  \draw[-{Latex[length=2mm]}, thick] (a1)--(a2);
  \draw[-{Latex[length=2mm]}, thick] (a2)--(a3);
  \draw[-{Latex[length=2mm]}, thick] (a3)--(a4);
\end{tikzpicture}}
\caption{Multilevel QUBO coarsening: aggregate the graph, solve at the coarsest level, then interpolate and refine.
Coarsening is approximate.
Cut-set elimination \citep{ponce2025graph} and recursive splitting \citep{pelofske2021decomposition} have different preservation conditions; exactness requires the relevant component or leaf solves to be exact.}\label{fig:multilevel-coarsening}
\end{figure}

A QUBO is a weighted graph, making three approaches natural: multilevel coarsening, cut-set elimination, and exact divide-and-conquer.
Subgraph methods solve small components on quantum computers and stitch the results together \citep{shaydulin2019hybrid}; multilevel methods coarsen, solve, and refine in analogy to classical multigrid \citep{ushijima2021multilevel}, including a scheme for maximum-cut \citep{angone2023hybrid}.
The downward path in Figure~\ref{fig:multilevel-coarsening} replaces the graph by fewer aggregate variables.
Mapping a coarse assignment back supplies a candidate on the original graph; refinement seeks improvements under the original objective because aggregation need not preserve an optimum.
%


Cut-set elimination \citep{ponce2025graph} uses component value functions to construct a reduced, generally reweighted graph, exactly or approximately depending on the reduction and component solves.
Both component and reduced-problem solves are interfaces whose accuracy matters.
Approximation-loss bounds for classical sparsification and decomposition of nonnegative-weight maximum-cut target circuit resources and noise in trapped-ion compilations \citep{moondra2024promise}.
These guarantees control transformation error, not the quality of an arbitrary finite-depth QAOA output.

Divide-and-conquer algorithms split maximum-clique and minimum-vertex-cover instances on vertices until the resulting QUBOs fit the connectivity and qubit limits of D-Wave quantum annealers \citep{pelofske2019maxclique,pelofske2019vertexcover,pelofske2021decomposition}.
The early studies target the D-Wave 2X (1152 qubits; Chimera graph with 3360 couplers, each qubit coupled to at most six others) \citep{king2015timetotarget} and the D-Wave 2000Q (up to 2048 qubits; Chimera graph with 6016 couplers) \citep{pelofske2019maxclique,boothby2020next}.
These are nominal architecture counts: fabrication defects leave some qubits and couplers inactive, for example 1097 qubits and 3045 couplers on one D-Wave 2X chip \citep{king2015timetotarget}, and minor-embedding chains further reduce the usable problem size.
The maximum-clique study therefore splits subgraphs until they have at most 46 vertices on the 2X and 65 on the 2000Q \citep{pelofske2019maxclique}.
Complete recombination with optimal leaf solves is exact, but annealing itself does not guarantee those solves.
Classical branch-and-bound also provides the skeleton for hybrid knapsack and traveling-salesperson solvers \citep{sanavio2024hybrid} and for a QUBO study \citep{peng2025hybrid}.
Their samplers propose incumbents; node closure requires certified optimality, infeasibility, or valid bound-based pruning, not sampling failure.
Tree growth and coordination can remain the bottleneck \citep{haner2024solving}.
Software implementations of these interfaces are discussed in Section~\ref{sec:decomposition-software}.
The SPLIT framework coordinates parallel subproblems through effective local-field updates (Section~\ref{sec:ising-qubo}) and greedy sweeps, fixing outside-block variables when updating each block \citep{vandelli2026parallel}.
Its experiments use IBM \software{CPLEX}, a classical mathematical-programming solver, to solve the subproblems and provide a baseline for decomposition; quantum subsolvers remain prospective, and the heuristic has no general convergence proof.

\paragraph{Quantum preprocessing.}
Short QAOA circuits can also prepare input for a classical solver rather than solve the problem directly.
Quantum preconditioning transforms the input problem using correlations from a QAOA circuit, with circuit depth setting the level of preconditioning, and then passes the transformed instance to classical heuristics such as simulated annealing or the Burer--Monteiro method \citep{dupont2025preconditioning}.
A constrained variant gives a mixed-integer programming solver an objective built from QAOA correlations while retaining the original hard constraint and evaluating incumbents under the original objective \citep{ramesh2026constrainedpreconditioning}.
Relax-and-round solvers can likewise replace a classical relaxation with QAOA correlations before classical rounding; with depth-one circuits this matches its classical counterpart on several problem families, and a modified rounding recovers the Goemans--Williamson performance on certain graphs \citep{dupont2024relaxround}.
A self-consistent mean-field preprint instead splits an Ising problem into independent subproblems whose mutual influence is represented by an environment computed with a variational circuit \citep{dupont2026meanfield}.
Section~\ref{sec:performance-attribution} assesses the evidence, including circuit and transformation costs.

A portfolio pipeline applies random-matrix preprocessing, modified spectral clustering, and risk rebalancing, using an approximate decomposition that also benefits classical solvers \citep{acharya2025decomposition}.
An MIS pipeline combines exact graph reduction by lazy clique removals with compatibility testing and embedding on Rydberg arrays \citep{schuetz2025compilation}.

Symmetry reduction also applies to lattice problems.
A \emph{lattice} is the mathematical set \(\{Bz:z\in\mathbb{Z}^n\}\) generated by the integer combinations of the linearly independent columns of a basis matrix \(B\).
Its shortest vector problem is the integer quadratic program
\[
\min_{z\in\mathbb{Z}^n\setminus\{0\}}\ \|Bz\|_2^2=z^\top B^\top Bz.
\]
Approximate shortest-vector problems on structured lattices motivate hardness assumptions in post-quantum cryptography.
One method exploits the symmetries of those lattices to encode a reduced search space as a Hamiltonian and then applies VQE (Section~\ref{sec:quantum-algorithms}) \citep{schirman2026efficient}.
The algebraic reduction is exact within the selected subspace, while its qubit implementation also bounds the integer coefficients; neither restriction guarantees preservation of the global shortest vector.
The reduction also applies to classical lattice algorithms.

A multilevel QUBO study \citep{maciejewski2024multilevel} uses extended Noise-Directed Adaptive Remapping, which relabels bits so device noise favors previously found good solutions, and relax-and-round subsolvers \citep{dupont2024relaxround}.
A semidefinite-inspired pipeline likewise executes device-sized subcircuits within classical preprocessing and postprocessing \citep{dupont2025benchmarking}, with a companion circuit-execution study \citep{maciejewski2024design}.
A quantum local-search method solves and merges small MIS neighborhoods \citep{tomesh2022quantum}.
Shallow, fixed-parameter sampling can also steer classical search between energy basins \citep{zhu2026combinatorial}.
A recent preprint rotates an active neighborhood frontier under classical control while keeping circuit width bounded \citep{araujo2026hybrid}; it is not treated as settled evidence.

\paragraph{Quantum relaxations and readout.}
Quantum random-access optimization (QRAO) encodes up to three binary variables per qubit, maximizes the energy of a generally non-diagonal Hamiltonian, and rounds measurements in several bases into assignments \citep{fuller2024relaxations}.
For unweighted maximum-cut, the Hamiltonian's maximum energy bounds the optimal cut from above; a variational candidate energy need not.
A randomized rounding guarantee assumes that the prepared state has energy at least as high as the optimal cut value; variational preparation need not satisfy this condition.


\subsection{Quantum-informed recursive reduction}\label{sec:recursive-reduction}

Recursive QAOA estimates spin correlations from repeated measurements: for spins \(s_i,s_j\in\{-1,1\}\), the mean product \(\langle s_i s_j\rangle\) indicates whether the sampled values tend to agree or disagree \citep{bravyi2020obstacles}.
A strong positive or negative correlation suggests fixing \(s_i=s_j\) or \(s_i=-s_j\); substituting this relation eliminates a variable, after which QAOA is run again on the reduced instance.
This repeated solve--reduce cycle is the recursion, and the suggested relation need not hold in an optimum.
At QAOA depth \(p=1\), meaning one cost-and-mixing round per quantum call, an iterative MIS construction produces the same operations and selections as classical greedy MIS \citep{brady2024iterative}.
This is distinct from the number of recursive reductions; departing from that equivalence requires more QAOA rounds or modified selection rules.
Quantum-informed recursive optimization generalizes this approach: measured single-spin means and pair correlations guide problem-specific classical reductions, with backtracking to repair earlier decisions \citep{finzgar2024qiro}.
A constraint-aware shrinking and repair pipeline uses simulated variational routines without hardware and is discussed here as a mechanism example \citep{sharma2025adaptive}.
In qReduMIS, \emph{kernelization} repeatedly applies classical graph reductions that preserve an optimum and record how to reconstruct it from the remaining graph, called the kernel \citep{schuetz2026qredumis}.
For example, an isolated vertex can be selected and deleted because it belongs to every maximum independent set.
When no implemented exact reduction applies, sampled inclusion frequencies identify vertices to fix inside or outside the independent set; those heuristic decisions may enable further exact reductions, and the cycle repeats.

This family combines recursive orchestration with exact and heuristic reductions.
Kernelization preserves an optimum through known rules; a probability-guided fixing remains heuristic unless certified, exhaustively backtracked, or included in complete branch-and-reduce.
The quantum routine supplies correlations or candidate fixings, while classical routines supply feasibility checks, kernelization, reconstruction, and sometimes backtracking.
Section~\ref{sec:performance-attribution} examines the fixed-wrapper ablations that distinguish these contributions.

\subsection{Handling constraints: alternating direction, Benders, and column generation}\label{sec:decomposition-interfaces}

\paragraph{Block updates.}
ADMM splits variables into blocks coupled by consensus and alternates block and multiplier updates.
Its convex theory \citep{boyd2011distributed} does not transfer automatically to mixed-binary variants.
A mixed-binary ADMM architecture assigns a binary block to a quantum or Ising solver and continuous convex blocks to classical solvers \citep{gambella2020multiblock}.
Its variant-specific assumptions give stationary-point convergence for a soft-constrained augmented Lagrangian, without generally ensuring original feasibility or global optimality.
The stated theorems assume exact QUBO solves; the paper discusses extensions with vanishing, summable errors, not arbitrary sampling errors.
This architecture has also been applied to routing \citep{harwood2021formulating}.

Variational methods that keep slack variables or Lagrange multipliers classical, introduced in Section~\ref{sec:quantum-algorithms}, apply the same division of labor inside a single algorithm \citep{braine2021transaction,le2023vqec}.

\paragraph{Master proposals and cut selection.}
In most cited Benders schemes, the quantum routine handles the binary master while a classical linear or convex subproblem generates cuts.
Examples include mixed-integer programming \citep{chang2020hybrid}, QUBO-encoded masters \citep{zhao2022hybrid}, and centralized and distributed generalized Benders for unit commitment \citep{gao2025distributed}.
For transmission-network expansion, embedding preprocessing, conservative handling of cut constraints, and stopping rules have been adapted to heuristic annealers \citep{lopezbanos2026performance}.
A neutral-atom master uses tightened continuous-variable bounds to limit qubit count, retaining classical linear recourse \citep{naghmouchi2024mixed}.
This architecture also targets transmission switching and verification of power-system neural networks \citep{ellinas2024hybrid}.
Other discrete-continuous methods place quantum-annealed combinatorial solves inside classical search \citep{ajagekar2020quantum}.

Master-side cut selection in multi-cut Benders supplies a different QUBO interface \citep{paterakis2023hybrid}, also explored with gate-based QAOA \citep{devalk2026hybrid}.
Keeping only valid cuts keeps the master a valid relaxation, but permanently excluding needed cuts can prevent convergence.
Exact master and subproblem solves give the standard finite-convergence argument; safeguarded inexact variants can use other conditions \citep{rahmaniani2017benders}.
The quantum contribution in a master-sampling implementation is candidate generation: an encoded master assignment is decoded, checked against the current cuts, and passed to classical recourse.
Valid dual solutions or infeasibility certificates generate the corresponding cuts, while a feasible original solution updates the incumbent.
A sampled master objective is not a certified lower bound; termination requires a verified master bound or another valid stopping certificate, with a classical fallback when heuristic sampling fails to make progress.
The quantum substitution therefore changes the master-search procedure, while cut validity and final certification retain their classical proof obligations.

\paragraph{Pricing and restricted masters.}
In column generation, neutral-atom sampling supplies pricing solutions for minimum vertex coloring while the restricted master remains classical \citep{dasilvacoelho2023quantum}.
For fleet assignment, engineered embeddings and pulses generate columns, followed by diversifying postprocessing \citep{perron2026leveraging}.
Emulated neutral-atom pricing also targets fidelity-constrained routing, with classical path validation and restricted integer-master recovery \citep{naghmouchi2026neutral}.
A multi-agent pathfinding method reverses the roles: classical pricing generates paths and a quantum restricted master selects compatible sets \citep{gerlach2025hybrid}.
Its path-generation construction can accumulate a sufficient pool despite inexact master solutions; optimality of its branch-and-cut-and-price variant requires the specified restricted-master accuracy.
In general, certified termination requires exact pricing or a reliable no-improving-column certificate.
Failure of heuristic sampling to find an improving column is not such a certificate.
A proposed exact vehicle-routing branch-price-and-cut framework encodes pricing and separation as QUBOs, using every returned solution below a cost threshold \citep{wagner2026quantum}.
The experiments implement the master relaxation, pricing, and separation.

Dantzig--Wolfe branch-and-price for refinery scheduling keeps the master and continuous subproblem classical, while binary pricing goes to an annealer through \software{qbsolv} \citep{ossorio2022refinery}.
Dual prices induce columns and branching restores integrality \citep{barnhart1998branch}.
Pricing optimality had to be verified afterward by a classical global solver, preventing premature termination from missed improving columns.

\paragraph{Cutting planes and augmentation.}
A cutting-plane preprint combines a linear programming (LP) bound and graph reduction for maximum weighted independent set with neutral-atom samples and sample-informed classical odd-cycle separation \citep{perron2026iterative}.
For \(t\)-perfect graphs, variable bounds \(0\le x\le1\), edge inequalities, and odd-cycle inequalities describe the independent-set polytope.
Polynomial-time separation permits polynomial-time optimization with a suitable classical algorithm; it does not by itself bound the implemented cut-adding loop's iterations.
Another interface samples Graver-basis elements---integer kernel directions for improving a feasible solution while preserving its linear equalities---for augmentation in nonlinear integer programming \citep{alghassi2019graver}.
A global-optimality guarantee requires an objective class satisfying the relevant Graver optimality conditions (such as separable convex minimization) \citep{alghassi2019graver}, a complete sufficient test set, and valid augmentation and stopping; sampling does not certify completeness.

\paragraph{Commercial interfaces.}
Hosted services execute models as well as build them.
D-Wave \software{Leap} accepts binary and constrained quadratic models through hybrid solvers \citep{dwave_hybrid}, and nonlinear models through an interface with native constraints \citep{dwave_nonlinear_docs}; the open-source \software{dimod} and \software{Qiskit Optimization} libraries of Section~\ref{sec:encoding-costs} instead construct and convert models for a chosen solver.
Quantinuum \software{Nexus} instead exposes compilation, execution, and data management for user programs \citep{quantinuum_nexus_docs}; its \software{Guppy} documentation demonstrates a user-built QAOA loop for maximum-cut, with classical parameter optimization and a \software{Selene} emulator backend in the example \citep{quantinuum_guppy_qaoa_docs}.
These interfaces place different responsibilities on the user: submitting a constrained mathematical model differs from constructing its quantum circuit and classical optimization loop.
The patent record includes mixed-binary decomposition and discrete-continuous optimization \citep{gambella2019patent,you2023patent}, sample-guided variable fixing and decomposition into subproblems \citep{karimi2018patent}, and integration of quantum samples into classical branch-and-bound \citep{mcgeoch2024patent}.
Vendor documentation and patents describe available interfaces or disclosed designs; they do not by themselves establish a quantum contribution to performance.
Their dates and evidential status follow Appendix~\ref{app:survey-method}.

\subsection{Reformulation into continuous cones}\label{sec:conic-reformulations}

%
\begin{figure}[tbp]
\centering
\scalebox{0.76}{\begin{tikzpicture}[
  box/.style={rectangle, rounded corners, draw=black, thick, align=center, minimum height=1.6cm, text width=4.3cm},
  every node/.style={font=\small}
]
  \node[box, fill=blue!6] (m) {\textbf{Cutting-plane master (classical)}\\ relaxes the copositive dual;\\ polyhedral outer approximation\\ of the copositive cone};
  \node[box, fill=orange!12, right=3.1cm of m] (o) {\textbf{Ising / QUBO oracle}\\ witness search:\\ minimize $z^{\top} M(\eta)\, z$\\ over encoded $z \ge 0$};
  \draw[-{Latex[length=2.5mm]}, thick] (m.north east) to[bend left=18]
     node[above, midway, font=\scriptsize]{candidate $M(\eta)$} (o.north west);
  \draw[-{Latex[length=2.5mm]}, thick] (o.south west) to[bend left=18]
     node[below, midway, font=\scriptsize]{sampled bits encoding candidate $z$} (m.south east);
\end{tikzpicture}}
\caption{Copositive separation \citep{brown2024copositive}: classical verification converts a negative witness into a cut for the lifted formulation \citep{burer2009copositive}.
The matrix \(M(\eta)\) is determined by the dual master decision vector \(\eta\).
The exact-oracle guarantee requires continuous separation; finite-grid encoding and heuristic sampling introduce distinct limitations (Section~\ref{sec:conic-reformulations}).}\label{fig:copositive-separation}
\end{figure}
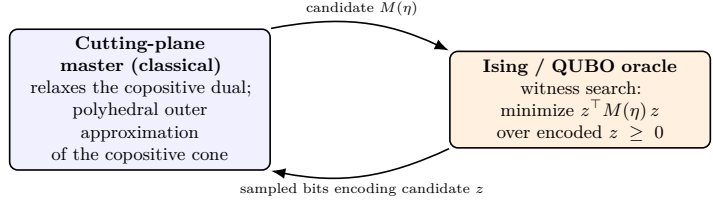


Lifting moves a problem into a higher-dimensional formulation.
Examples include Dantzig--Wolfe reformulation of linear programs \citep{dantzig1960decomposition}, moment and sum-of-squares semidefinite hierarchies \citep{lasserre2001global,parrilo2000structured}, and the exact \emph{completely positive} representation of mixed-binary quadratic programs \citep{burer2009copositive}.
A completely positive matrix has a factorization \(X=BB^\top\) with entrywise \(B\ge0\); a symmetric matrix \(M\) is copositive when \(z^\top Mz\ge0\) for every \(z\ge0\).
The completely positive reformulation relocates hardness into conic membership: non-copositivity testing is NP-complete \citep{murty1987some}, and strong membership of both the completely positive and copositive cones is NP-hard \citep{dickinson2014computational}.

A copositive framework derives a dual with the same optimal value as the mixed-binary primal under stated assumptions; equality with the dual supremum does not imply dual attainment \citep{brown2024copositive}.
In Figure~\ref{fig:copositive-separation}, the outgoing matrix \(M(\eta)\), determined by the dual master decision vector \(\eta\), defines the separation task; the return arrow carries a sampled encoding of a candidate vector \(z\).
Classically decoding it and checking \(z\ge0\) and \(z^\top M(\eta)z<0\) yields a valid separating cut for the continuous dual: one verified violation suffices.
With an exact continuous separation oracle and appropriate bounded localization and tolerance assumptions, cutting-plane variants such as the ellipsoid method approximate the dual value with polynomial oracle complexity and classical overhead.
The finite-bit implementation has a weaker interface: a negative grid witness still gives a valid cut, but even a certified zero grid minimum establishes only approximate nonnegativity when a discretization-error bound is available.
Failure of heuristic sampling to find a violation provides no such certificate.
The supplement derives a grid-error bound separately from the paper's exact-oracle result; exact reformulation, discretization, and finite-solver accuracy remain distinct.

Q-FW uses the same reformulation route with Frank--Wolfe: linear-minimization steps become QUBOs for an annealer \citep{yurtsever2022qfw}.
Q-FW replaces a single static exact-penalty weight with a prescribed augmented-Lagrangian or quadratic-penalty schedule.
Its stated convergence result assumes exact linear minimization, a nonempty feasible set, strong duality, and bounded dual iterates; heuristic annealer calls do not automatically satisfy these conditions, and finite-iteration rounding may require feasibility repair.

\subsection{Software for decomposition and reduction}\label{sec:decomposition-software}

Software exposes which variables a decomposition releases, how it treats their coupling to the remainder, and how it reconstructs a complete assignment.
Supplementary Table S7 catalogs 23 public implementations, distinguishing libraries, research prototypes, and educational code with inspected source versions and scope limits.
This targeted inventory establishes implemented interfaces, not comparative performance or independently validated correctness.

For block search, the archived \software{qbsolv} combines conditioned sub-QUBO updates with tabu search; its default subsolver is classical, with annealing optional \citep{booth2017partitioning}.
\software{dwave-hybrid} separates decomposers, samplers, and composers, supporting energy-impact, graph-traversal, and random neighborhoods \citep{dwave_hybrid_framework_docs}.
\software{QSplit} orchestrates split--solve--aggregate workflows through \software{StreamFlow}; current source offers graph, matrix, neighborhood, and macrovariable splitting, voting, and conditioned or consensus-based refinement \citep{bifulco2027qsplit,qsplit_software}.
\software{HADOF} instead substitutes sampled marginal means for exterior variables, so its subproblems differ from exact restrictions at a fixed incumbent \citep{Sankar2025HADOF,Sankar2026ParallelHADOF}.
The inspected \software{DQAOA} implementation proposes updates from principal submatrices and checks accepted bits against the full objective, without incorporating all omitted couplings into the local linear terms \citep{Kim2026DistributedQAOA}.
\software{Pasqal QUBO Solver} builds device-layout subproblems, requires nonnegative off-diagonal coefficients, and solves remaining variables classically \citep{pasqal_qubo_decomposition_docs}.
The registered \software{QUBODecomposition.jl} package provides connected-component solves and bounded conditioned sweeps through Julia optimizer interfaces \citep{qubodecomposition_software}.

Graph-specific implementations offer additional interfaces.
\software{ML-QLS} and \software{MLQAOA} coarsen, solve, project, and refine graph problems; the latter also transfers QAOA parameters using graph representations \citep{ushijima2021multilevel,bach2024mlqaoa}.
\software{QAOA-in-QAOA} merges maximum-cut subgraph assignments by optimizing whole-subgraph flips, while \software{Neural QAOA squared} learns partitioning and parameter initialization \citep{Zhou2023QAOAinQAOA,Zheng2026NeuralQAOA}.
These restricted merges are not arbitrary QUBO recombinations.
For recursive reduction, \software{OpenQAOA} implements correlation elimination and reconstruction; the \software{QIRO} reference code and \software{Qrisp} provide problem-specific replacement rules, whereas the public \software{qReduMIS} solver targets atom-position MIS with Rydberg backends \citep{sharma2022openqaoa,finzgar2024qiro,qrisp2026qiro,schuetz2026qredumis}.
\software{DASQA} supplies clique and vertex-cover branching with configurable leaf solvers; exactness requires exact leaves and valid pruning \citep{pelofske2021decomposition}.

Preprocessing is another reusable layer: \software{dwave-preprocessing} exposes roof-duality persistencies and component solves, \software{Qoolchain} combines reductions with branching, and \software{QUBOLite} supplies implication-based fixings and reconstructable partial assignments \citep{dwave_preprocessing_docs,orlandi2025qoolchain,mucke2025qubolite}.
For constrained interfaces, \software{Qiskit Optimization} implements mixed-binary ADMM and recursive correlation reduction; \software{HQC-Bend} delegates Benders masters, whereas \software{Logistiqs} keeps master and recourse solves classical and delegates cut selection \citep{qiskit_optimization_docs,zhao2025hqcbend,devalk2026hybrid}.
These interfaces have different proof obligations.
Even the word decomposition requires care: \software{Quafu-Qcover} decomposes QAOA expectation evaluation into local neighborhoods, rather than combining independently optimized blocks \citep{xu2024qcover}.
Across implementations, exact child solves alone do not certify coupled problems; comparisons must evaluate reconstructed assignments under the original objective and count preprocessing, training, repair, and coordination costs.
\subsection{Structure and guarantees as separate dimensions}\label{sec:structures-and-guarantees}

The structural vocabulary of Section~\ref{sec:classical-foundations} identifies what an algorithm does; its proof obligations are interface-specific.
Table~\ref{tab:hybrid-roles} summarizes candidate production and classical checks.
Table~\ref{tab:guarantee-requirements} gives the additional conditions needed for the stated guarantees, keeping coarsening, exact reduction, certified recursion, and heuristic fixing distinct.
%
\begin{table}[tbp]
\centering
\footnotesize
\setlength{\tabcolsep}{3pt}
\setlength{\tablewidth}{\dimexpr\linewidth-12pt\relax}
\begin{tabular}{@{}>{\RaggedRight\arraybackslash}p{0.24\tablewidth}>{\RaggedRight\arraybackslash}p{0.36\tablewidth}>{\RaggedRight\arraybackslash}p{0.40\tablewidth}@{}}
\toprule
\textbf{Family} & \textbf{Quantum task and returned object} & \textbf{Classical coordination and checks} \\
\midrule
Variational &
Sample states; assignments and expectation estimates &
Tune parameters; evaluate objectives and feasibility \\
Multilevel coarsening &
Sample coarse QUBOs; assignments &
Coarsen, interpolate, refine; evaluate original objective \\
Cut-set elimination &
Sample component and reduced QUBOs; assignments and values &
Construct value functions; reconstruct and check assignments \\
Recursive splitting / branch and bound &
Sample leaves; candidate incumbents &
Split and recombine; manage bounds and pruning in branch and bound \\
ADMM mixed-binary &
Sample binary block; candidate block assignment &
Solve continuous blocks; update multipliers; check residuals and original feasibility \\
Benders master &
Propose master assignments &
Solve recourse; verify incumbents, cuts, and master bounds \\
Benders cut selection &
Propose cut subsets &
Check cut validity; solve master and recourse; safeguard cut access \\
Column generation: quantum-pricing variant &
Propose pricing columns &
Check feasibility and reduced cost; solve restricted master; certify final pricing; branch \\
Copositive separation &
Propose witness vectors &
Check \(z\ge0\) and \(z^\top Mz<0\); form cuts; maintain outer approximation \\
Copositive Frank--Wolfe &
Sample linear-minimization QUBOs; assignments &
Update lifted iterate, penalty and duals; round and repair \\
Graver augmentation &
Propose test-set elements &
Check integer directions preserve equalities; verify feasible moves; augment \\
Quantum-informed fixing &
Estimate correlations; propose fixings &
Apply exact reduction rules; reconstruct, check feasibility, and backtrack \\
\bottomrule
\end{tabular}
\caption{Interface allocation for representative hybrid families.
Returned candidates support the listed checks; the additional requirements for guarantees appear in Table~\ref{tab:guarantee-requirements}, in the same family order.}\label{tab:hybrid-roles}
\end{table}


%
\begin{table}[tbp]
\centering
\scriptsize
\setlength{\tabcolsep}{3pt}
\setlength{\tablewidth}{\dimexpr\linewidth-18pt\relax}
\begin{tabular}{@{}>{\RaggedRight\arraybackslash}p{0.20\tablewidth}>{\RaggedRight\arraybackslash}p{0.19\tablewidth}>{\RaggedRight\arraybackslash}p{0.36\tablewidth}>{\RaggedRight\arraybackslash}p{0.25\tablewidth}@{}}
\toprule
\textbf{Family} & \textbf{Transformation status} & \textbf{Additional requirement} & \textbf{Resulting guarantee} \\
\midrule
Variational &
Orchestration &
Problem- and parameter-specific assumptions &
Special-case guarantees only \\
Multilevel coarsening &
Approximation &
Reduction-specific error analysis &
No general quality bound \\
Cut-set elimination &
Exact or approximate reduction &
Exact value functions and certified reduced solve &
Global optimum in exact variant \\
Recursive splitting / branch and bound &
Exact splitting &
Complete recombination; node closure by optimality, infeasibility, or valid bounds &
Certified optimum or infeasibility \\
ADMM mixed-binary &
Split; possibly softened constraints &
Variant-specific regularity and penalties; exact blocks or justified error control &
Conditional stationarity; original feasibility and optimality not implied \\
Benders master &
Exact partition &
Valid cuts, certified bounds, safeguarded progress &
Optimality when bounds meet \\
Benders cut selection &
Formulation unchanged &
Valid cuts; continued access to required cuts &
Bounds valid; convergence conditional \\
Column generation: quantum-pricing variant &
Exact reformulation &
Certified final pricing and master bounds; complete branching &
LP optimum; integer optimum with completed tree \\
Copositive separation &
Exact lift under stated assumptions &
Exact continuous separation and localization assumptions; grid error treated separately &
Dual-value approximation by cutting planes, polynomial relative to exact oracle; negative witness certifies violation \\
Copositive Frank--Wolfe &
Exact lift under stated assumptions &
Exact linear minimization, prescribed schedule, duality and bounded dual iterates &
Conditional objective and feasibility convergence \\
Graver augmentation &
Test-set framework &
Objective class satisfying Graver optimality conditions; complete sufficient test set; valid augmentation and stopping &
Global optimum for the covered objective class \\
Quantum-informed fixing &
Exact rules plus heuristic fixings &
Certified fixings or exhaustive backtracking &
Otherwise heuristic \\
\bottomrule
\end{tabular}
\caption{Guarantee requirements complementing Table~\ref{tab:hybrid-roles}.
Feasibility, objective values, and negative witnesses can be checked; heuristic samples alone do not certify optimality or the absence of improvement.
Conditions refer to the specific variants discussed in Section~\ref{sec:hybrid-methods}.}\label{tab:guarantee-requirements}
\end{table}


\FloatBarrier


%
\section{Attributing performance gains}\label{sec:performance-attribution}

A transformed pipeline can beat a monolithic quantum solver because of its classical wrapper, its quantum component, or their interaction.
Three questions distinguish these explanations:
\begin{itemize}
\item \emph{Component effect:} does substituting a quantum routine improve performance within a fixed classical wrapper at a stated quality or resource budget?
\item \emph{System effect:} does a separately optimized quantum-enabled pipeline outperform a separately optimized classical pipeline at comparable quality and resources, including all end-to-end costs?
\item \emph{Reduction--solver interaction:} does a reduction have a more favorable effect on the quantum solver than on a classical comparator on a common performance and resource scale?
\end{itemize}
A fixed wrapper may favor the solver for which it was designed; a system comparison permits each pipeline its own tuning.
For an interaction, let \(Y(r,\sigma)\) measure predeclared higher-is-better solution quality at a fixed common resource budget, where \(r\in\{0,1\}\) indicates absence or presence of a reduction and \(\sigma\in\{\mathrm{Q},\mathrm{C}\}\) selects the quantum or classical solver, respectively.
Define the simple effects \(\Delta_{\mathrm{Q}}=Y(1,\mathrm{Q})-Y(0,\mathrm{Q})\) and \(\Delta_{\mathrm{C}}=Y(1,\mathrm{C})-Y(0,\mathrm{C})\), and their contrast
\begin{equation}\label{eq:interaction-quality}
\Delta_{\mathrm{int}}=
\bigl[Y(1,\mathrm{Q})-Y(0,\mathrm{Q})\bigr]-\bigl[Y(1,\mathrm{C})-Y(0,\mathrm{C})\bigr].
\end{equation}
A positive contrast means a more favorable reduction effect for the quantum solver relative to the classical comparator.
An absolute quantum-solver improvement additionally requires \(\Delta_{\mathrm{Q}}>0\); report both simple effects alongside the contrast.
Neither condition implies system superiority.
This factorial contrast is not automatically causal: assignment, tuning, and resource comparability determine that interpretation.
Uncertainty must cover the declared contrast itself.
The alternative of cost at a fixed quality target, estimators, and bootstrap procedures appear in the supplementary methods.
Stochastic solvers add a further requirement: their parameters must be set, and the cost of setting them counts toward the comparison.
A benchmarking framework scores parameter-setting strategies through the statistics of a declared performance metric across resource budgets \citep{bernalneira2025benchmarking}; its open-source \software{stochastic-benchmark} package supports designing, evaluating, and visualizing such strategies \citep{stochastic_benchmark_repo}.

Effect type and evidence status are distinct.
A component comparison may use hardware, simulation, or nominal sample counts; a system comparison may report complete, partial, or projected timing.
Following established benchmarking cautions \citep{ronnow2014defining}, this survey counts end-to-end cost as preprocessing, parameter setting, mapping, queueing, communication, execution, decoding, verification, repair, classical coordination, and postprocessing where applicable.
A strong classical comparator is one identified within this survey's search, including relevant graphics processing unit (GPU) and physics-inspired methods at comparable quality.
Tables~\ref{tab:controlled-comparisons} and~\ref{tab:scaling-reassessments} retain each study's design, comparator, metric, and timing boundary; Appendix~\ref{app:survey-method} gives the survey scope and version policy.

\paragraph{Wrapper benefits and component effects.}
Controlled swaps distinguish useful quantum information from quantum superiority.
Within the same graph reduction, exact classical and QAOA component solves produce similar observed mean approximation ratios, near 0.96 versus about 0.75 for monolithic QAOA \citep{ponce2025graph}.
These observations support a reduction benefit, not formal equivalence or general solver insensitivity.
Both arms retain QAOA as the final solver, so neither is an entirely classical control.
In qReduMIS, quantum-informed reduction beats random selection and the wrapper improves on standalone optimized quantum annealing, yet SaReduMIS matches or exceeds the quantum arm on all 39 tested union-jack instances (site-diluted square grids with nearest and diagonal-neighbor edges) \citep{schuetz2026qredumis}.
SaReduMIS's observed success probability of one on this testbed is not a guarantee; 500 replicas and two simulated-annealing steps versus 1,000 quantum measurements are nominal counts, not matched time or energy.
The gate-based portfolio extension also supports a wrapper benefit, but retains QAOA in its comparator: 100 standalone shots and 20 hybrid trials with repeated calls use different budgets.
Cumulative-shot time-to-solution excludes classical and communication costs and does not uniformly favor the wrapper on hardware \citep{yalovetzky2026portfolio}.

Positive and negative component findings coexist.
In a Frank--Wolfe graph decomposition, simulated QAOA improves sparsity or error in the tested settings; IBM hardware beats simulated annealing at 50 and 70 nodes but loses at 100, without end-to-end timing \citep{pennington2025boosting}.
Partitioned QAOA extends this approach to 1,710 edge variables assembled from circuits of at most 50 qubits; with shared repair, it reports lower residual error than the tested annealing and random samplers, without repeat-based uncertainty or matched computational costs \citep{liaqat2026scaling}.
By contrast, specialized classical subproblem methods run faster than the quantum alternatives in pricing and separation experiments \citep{wagner2026quantum}.
These experiments implement the master relaxation and associated subproblems motivated by an exact branch-price-and-cut framework, not the complete integer solver.
Simulated annealing is the strongest sampler in the emulated cutting-plane wrapper \citep{perron2026iterative}.
In small power-grid decompositions, ideal variational sampling fails to improve on uniform sampling with the same continuous solver even when quantum evolution is charged zero time; this negative result is restricted to the 5--13-generator, convergence-selected set \citep{gaidai2026benchmarking}.

Classical coordination also determines whether the selected quantum task merits its cost.
In five toy hybrid pipelines, QAOA cut selection simulated classically with matrix product states consumes 99.3\% of wall time, although cut selection takes a small share of the classical baseline's runtime; a separate IBM experiment demonstrates only the component \citep{devalk2026hybrid}.
In emulated neutral-atom pricing, classical initialization narrows backend differences, and combined initialization and refinement bring the quantum arm's mean gaps below 1\% of \software{CPLEX} optima.
Neither wall-clock comparison nor exact final pricing is established \citep{naghmouchi2026neutral}.
In surveillance routing, a fixed-tour-pool comparison finds no systematic benefit of emulated neutral-atom selection over classical selection controls; the headline comparison against a full routing MIP projects QPU time from emulation, without hardware execution or complete pipeline costs \citep{lauwerier2026quantum}.
Varying quantum-call frequency without a matched zero-quantum arm and complete accounting \citep{holliday2026advanced} cannot isolate the marginal component contribution.

\paragraph{Preconditioning and co-design.}
For the preconditioners described in Section~\ref{sec:graph-decomposition}, accounting boundaries determine how much improvement can be attributed to the full pipeline.
On random 3-regular maximum-cut, the reported benefit includes correlation construction through emulation, while the hardware demonstration remains limited \citep{dupont2025preconditioning}.
For constrained bipartitioning, near-optimal incumbent discovery improves on 12--40 variables using classically evaluated correlations, whose construction cost is excluded \citep{ramesh2026constrainedpreconditioning}.
Main comparisons select each balance-penalty weight retrospectively to minimize that instance's time to the known near-optimal target; global and cross-validated choices retain benefits with smaller depth-one gains.
Rescoring on the original objective is essential because the surrogate need not preserve its optimum.

Positive one-in-three Boolean satisfiability (SAT) at clause-to-variable ratio \(m/n=0.626\) supplies favorable co-design evidence \citep{lu2026evidence}.
The restricting-space reduction improves the quantum methods relative to classical controls also given reduction-related information.
Fits through \(n=70\) favor the reported quantum scaling bases, but quantum cost is inverse success probability whereas classical cost counts conflicts and propagations.
The reduction also changes the search manifold, qubit count, Hamiltonian, and ansatz.
Thus this predominantly simulated evidence is neither a fixed-wrapper swap nor an instance of Equation~\eqref{eq:interaction-quality} on common resources.
A 13-qubit experiment reaches 22-variable instances, without establishing a hardware crossover under complete wall-clock accounting.

\paragraph{System performance and accounting boundaries.}
High solution quality, device-access time, and complete system speed answer different questions.
One superconducting study attains high approximation ratios, yet its projected run times trail the strongest classical heuristics, with no fixed-wrapper attribution \citep{dupont2025benchmarking}.
In another, experimental samples steer a fixed classical search and reach global minima 2.23 times more often than tuned simulated annealing on the 20 hardest of 4,000 nearest-neighbor Ising instances under a projected 0.4-megahertz normalization \citep{zhu2026combinatorial}.
The approximately 2.34-fold time-to-solution gain is fully projected against single-core simulated annealing, without GPU or planar-specific controls.
On trapped-ion hardware, classically warm-started QAOA with local search attains higher quality and probability of finding an optimal cut than implemented Goemans--Williamson \citep{goemans1995improved} and Halperin--Livnat--Zwick \citep{halperin2004maxcut} algorithms on five 96-node graphs, with local-search postprocessing also applied to classical outputs \citep{he2026regularized}.
Its runtime crossover against optimized rank-two Burer--Monteiro--Zhang heuristics \citep{burer2002ranktwo}, with classical local search and iterative refinement enabled, remains a fault-tolerant resource projection \citep{he2026regularized}.

Proprietary hybrids impose further attribution limits.
One large dense-instance comparison normalizes accuracy by the best observed objective, hides the quantum/classical division, and times quantum processing unit (QPU) access while excluding reading, queueing, and local communication \citep{kim2025quantum}.
Its reported quality and speed ratios therefore establish neither a component effect nor a complete system comparison.
Hybrid and simulated annealing reach the optima on the certified small instances within a 139-instance maximum-cut benchmark, whereas direct QPU execution does not; larger-instance quality favors simulated bifurcation and slower simulated annealing \citep{vodeb2026accuracy}.
Heterogeneous hardware and some literature-sourced timings preclude a common-environment efficiency ranking.
An audit of a portfolio service finds matches to all 54 certified head-to-head optima, with mean QPU access of 0.034 seconds in a five-second nominal budget \citep{lozano2026where}.
The revised preprint's matched-time tabu control has mean absolute objective gap 0.00078 to the service's objective over 36 instances and 108 seeded runs.
The small observable device-access share leaves the service's remaining work and causal contribution unresolved.

Selected comparisons nevertheless give restricted positive or mixed system evidence.
A run-time benefit is reported on specially constructed sparse instance classes \citep{tasseff2024emerging}.
Across 20 higher-order binary instances, the best energies over hybrid trials and settings match independently reported ground states on 14.
At the hybrid's mean solver-side runtime, several central processing unit (CPU) comparators fall short of these hybrid-defined targets on average, whereas enhanced parallel tempering and a GPU solver match or surpass them; orchestration and deployment costs remain excluded \citep{chandarana2026quest}.
Contextual classical controls further narrow interpretations.
GPU quantum-annealing-inspired samplers improve sampling-time comparisons with quantum hardware and complete-pipeline comparisons with classical solvers on related multi-objective benchmarks; the complete quantum-system comparison remains unresolved \citep{tao2026multiobjective}.
The largest industrial scheduling results in one study come from a classical Digital Annealer with cloud integration excluded \citep{safi2026reality}.
A product-state control is a useful candidate baseline, but best-of-ten cuts, trained schedules, and coloring excluded from timing prevent a matched advantage claim or general conclusion about entanglement \citep{bantysh2026disentangling}.

\paragraph{Scaling and reassessment.}
Scaling conclusions depend on the benchmark family and the classical frontier.
The claimed scaling advantage in time to reach a relative objective gap \(\varepsilon\), for targets of at least one percent on tailored high-precision spin glasses \citep{munoz2025scaling}, is closed by GPU simulated bifurcation \citep{pawlowski2026toward}; accessible sizes remain insufficient for robust asymptotic inference.
Likewise, the reported superlinear speedup over optimized simulated annealing on hard native unit-disk instances \citep{ebadi2022quantum} is qualified by a classical reassessment that solves structured quasiplanar instances with thousands of vertices in minutes and obtains competitive simulated annealing results \citep{andrist2023hardness}.
On harder families constructed through density and treewidth, quantum runs are slower on the small instances compared exactly; classical scaling suggests roughly 1,000 atoms at kilohertz repetition as a necessary step toward advantage, not an observed crossover \citep{cazals2026identifying}.
Optimal adiabatic schedules are found to depend principally on independent-set instance difficulty rather than size \citep{perseguers2025hardness}.
For gate-based methods, a QAOA scaling claim derives from idealized noiseless analysis and simulation; separate small-scale hardware demonstrations do not establish hardware scaling \citep{shaydulin2024evidence}.
These distinct forms of evidence require distinct resource interpretations.

The oracle model offers a complementary accounting check.
The polynomial classical cutting-plane loop \citep{brown2024copositive} does not remove the difficult separation requirement.
Holding the outer algorithm fixed isolates the oracle contribution; changes in iteration count and classical overhead, as well as oracle cost, require separate accounting.
Polynomial overhead can itself remain costly in dimension, tolerance, and conditioning.

The analyzed corpus supports conditional component benefits, useful wrapper improvements, and favorable co-design evidence, alongside negative comparisons and classical reassessments.
It does not establish robust, transferable system superiority under full end-to-end accounting against the strongest relevant comparable classical baselines considered.
Here robustness concerns repeated instances, seeds, hardware runs, and settings with uncertainty; transfer concerns independently generated or held-out instances.
Those conclusions are scoped to the purposive survey, and do not negate its component findings.
Declaring the effect, quality target, and timing boundary makes future comparisons interpretable; \software{QOBLIB} provides classically hard problem classes and solution records for such comparisons \citep{koch2026qoblib}.

\FloatBarrier
%
\begingroup
\footnotesize
\setlength{\tabcolsep}{3pt}
\setlength{\tablewidth}{\dimexpr\linewidth-18pt\relax}
\begin{longtable}{@{}>{\RaggedRight\arraybackslash}p{0.255\tablewidth}>{\RaggedRight\arraybackslash}p{0.205\tablewidth}>{\RaggedRight\arraybackslash}p{0.375\tablewidth}>{\RaggedRight\arraybackslash}p{0.165\tablewidth}@{}}
\caption{Controlled and black-box comparisons underlying Section~\ref{sec:performance-attribution}.
The supported inference is restricted to the stated design, comparator, and metric; a component or wrapper benefit does not imply system superiority.
Bracketed labels link to the supplementary audit, with publication-status checks dated \SurveyVerificationDate{}.
Shared terminology: CPU, central processing unit; MIP, mixed-integer programming; MPS, matrix-product-state simulation; SK, Sherrington--Kirkpatrick; TTS, time-to-solution.
Table~\ref{tab:scaling-reassessments} continues with scaling claims and reassessments.}\label{tab:controlled-comparisons}\\
\toprule
\textbf{Study and setting} & \textbf{Quantum / classical roles} & \textbf{Design, baseline, metric, and timing boundary} & \textbf{Supported\newline inference} \\
\midrule\endfirsthead
\multicolumn{4}{l}{\normalsize Table~\thetable\ continued.}\\
\toprule
\textbf{Study and setting} & \textbf{Quantum / classical roles} & \textbf{Design, baseline, metric, and timing boundary} & \textbf{Supported\newline inference} \\
\midrule\endhead
\multicolumn{4}{r}{\scriptsize continues on the next page}\\\endfoot
\bottomrule\endlastfoot
\textup{[E1]} \citep{ponce2025graph}: 25 maximum-cut instances, 100 vertices, 3-regular; depth-one QAOA via closed-form expectations; ten reduced instances (12--16 qubits) on trapped-ion hardware &
Component solves and final QAOA / graph reduction &
Reduction fixed; exact classical vs. QAOA component solver, with final QAOA in both arms.
Similar observed mean ratios near 0.96 vs. monolithic QAOA near 0.75; certified optima, no formal equivalence or timing claim &
Positive wrapper benefit; limited component sensitivity \\
\addlinespace[3pt]
\textup{[E2]} \citep{schuetz2026qredumis}: 39 union-jack instances, 7--231 nodes; Rydberg hardware; version details in Appendix~\ref{app:survey-method} &
Kernel samples / clique-based reduction, selection, reconstruction &
Wrapper fixed; quantum vs. random vs. simulated annealing.
Quantum samples improve on random; wrapper improves on standalone quantum annealing.
SaReduMIS matches or exceeds quantum on all 39 instances at nominal counts (500 replicas, two steps vs. 1,000 measurements), not matched time &
Positive wrapper benefit; no quantum component superiority \\
\addlinespace[3pt]
\textup{[E3]} \citep{dupont2025benchmarking}: 3-regular maximum-cut to thousands of variables; superconducting processor, small subcircuits &
Sampler / relax-and-round preprocessing and postprocessing &
Whole quantum pipeline vs. commercial MIP, simulated annealing, and Burer--Monteiro.
Ratios near 99\% against classical references; projected 100-qubit run time behind strongest heuristics; no fixed-wrapper substitution &
Negative projected system-speed comparison; attribution unresolved \\
\addlinespace[3pt]
\textup{[E4]} \citep{kim2025quantum}: over 50 large dense instances; annealer and proprietary hybrid layer &
Annealing / undisclosed internal decomposition &
Black-box solvers vs. classical heuristics.
Relative accuracy uses best observed objective.
Timing counts QPU programming and sampling; excludes reading, queueing, and local-machine communication &
Unresolved attribution; incomplete system timing \\
\addlinespace[3pt]
\textup{[E5]} \citep{tasseff2024emerging}: large sparse instances; quantum annealer &
Direct annealing / parameter setting and evaluation &
Hardware vs. strong classical heuristics; authors report run-time benefit only on specially constructed instance classes &
Limited positive system comparison \\
\addlinespace[3pt]
\textup{[E6]} \citep{vodeb2026accuracy}: 139 maximum-cut instances, 100--10,000 nodes; small-instance certified subset &
Direct and hybrid annealing / simulated annealing and bifurcation &
Black-box comparison.
Hybrid and simulated annealing reach certified optima where available; direct QPU does not.
Larger-instance quality favors bifurcation and slower simulated annealing.
Heterogeneous hardware, some literature timings &
Mixed system-quality comparison; attribution unresolved \\
\addlinespace[3pt]
\textup{[E7]} \citep{lozano2026where}: portfolio instances, 10--640 assets; proprietary service; preprint &
Hybrid service / exact, annealing, and tabu anchors &
Service matches all 54 certified head-to-head optima.
Mean QPU access 0.034 s of a 5 s nominal budget; remaining time unresolved.
Version 3's matched-wall-clock tabu control: 36 instances, 108 seeded runs, mean absolute objective gap 0.00078 to the service objective &
Competitive interface outputs; component attribution unresolved \\
\addlinespace[3pt]
\textup{[E8]} \citep{zhu2026combinatorial}: 20 hardest of 4,000 nearest-neighbor Ising instances; 104 qubits &
Shallow sampling steers basin transitions / fixed classical search &
Experimental samples reach global minima 2.23\(\times\) more often than tuned simulated annealing under projected 0.4 MHz throughput normalization.
About 2.34\(\times\) time-to-solution gain fully projected vs. single-core simulated annealing; no GPU or planar-specific baseline &
Limited positive sampling evidence; system timing projected \\
\addlinespace[3pt]
\textup{[E9]} \citep{pennington2025boosting}: complete and bipartite graphs, 7--25 qubits, simulated; heavy-hex graphs, 52--111 qubits, simulated and IBM Kingston; preprint &
QAOA matching samples / fully corrective Frank--Wolfe decomposition &
Fixed wrapper; QAOA vs. simulated annealing vs. random.
Simulated QAOA sparser on small graphs and leads on heavy-hex error.
Hardware beats simulated annealing at 50 and 70 nodes, loses at 100; no end-to-end timing &
Positive component benefit in some settings; mixed hardware result \\
\addlinespace[3pt]
\textup{[E10]} \citep{dupont2025preconditioning}: SK, 3-regular maximum-cut, grid energy; emulation and Ankaa-3 test &
QAOA correlations / preconditioning then simulated annealing or Burer--Monteiro &
Same solver on original vs. preconditioned input.
Faster convergence across families; benefit survives preconditioning cost for 3-regular maximum-cut via emulation.
200-instance averages, standard errors; limited hardware demonstration &
Positive preconditioning benefit in emulation; hardware system effect unresolved \\
\addlinespace[3pt]
\textup{[E11]} \citep{wagner2026quantum}: vehicle routing; master-relaxation experiments &
Annealing on pricing and separation / classical master relaxation and exact pricing fallback &
Pricing and separation motivated by an exact branch-price-and-cut framework.
Specialized classical methods run faster than the quantum alternatives; complete integer solver not implemented &
Negative quantum component comparison \\
\addlinespace[3pt]
\textup{[E12]} \citep{perron2026iterative}: weighted independent set to 300 vertices; noiseless emulation; preprint &
Analog candidates and statistics / LP bounds, reduction, odd-cycle cuts &
Fixed cutting-plane wrapper and equal sampling budgets; comparisons with direct analog protocols and greedy baselines.
Within 5--10\% of optimum at 300 vertices; simulated annealing strongest sampler.
Polynomial solvability on \(t\)-perfect graphs does not itself bound this loop's iterations; no end-to-end timing &
Wrapper benefit; no quantum component superiority \\
\addlinespace[3pt]
\textup{[E13]} \citep{devalk2026hybrid}: routing; ten permutations of one 20-customer instance; five 5-customer toys; preprint &
QAOA cut selection / Benders master and dual subproblems &
Classical bottleneck profile; MPS-simulation toy pipeline spends 99.3\% of wall time in cut selection.
IBM demonstration is separate; no hardware pipeline speedup &
Negative result for tested scale \\
\addlinespace[3pt]
\textup{[E14]} \citep{naghmouchi2026neutral}: routing; two topologies, up to 30 nodes; neutral-atom emulation; preprint &
Pricing samples / LP master, reductions, repair and integer recovery &
Shared wrapper versus simulated annealing; warm-start/refinement ablations.
Reinforced quantum arm's mean gaps below 1\% to \software{CPLEX}; 95\% intervals.
No timing or exact-pricing certificate &
Limited component sensitivity in simulation \\
\addlinespace[3pt]
\textup{[E15]} \citep{ramesh2026constrainedpreconditioning}: balanced bipartitioning, 12--40 variables; classical correlation evaluation; preprint &
QAOA correlations / modified-objective \software{Gurobi} with hard balance &
Original-objective incumbent quality improves sooner.
Correlation cost excluded; main penalties selected with oracle information, alternative choices also tested.
No original-optimum guarantee for surrogate &
Positive preconditioning benefit; system effect unresolved \\
\addlinespace[3pt]
\textup{[E16]} \citep{yalovetzky2026portfolio}: four market graphs up to 225 assets; 78-qubit kernels; 73 emulator graphs; preprint &
Gate-based samples / recursive reduction, fixing and repair &
Same first kernel versus standalone QAOA; unequal raw-success budgets.
Cumulative-shot TTS excludes classical/communication time; no classical sampler substitution in hardware comparison &
Positive wrapper benefit; component/system superiority unresolved \\
\addlinespace[3pt]
\textup{[E17]} \citep{liaqat2026scaling}: 1,710 edge variables; partitions up to 50 qubits; IBM Heron; preprint &
QAOA samples / partitioning, repair, matching and Frank--Wolfe updates &
Shared wrapper/repair; QAOA reports lower residual than tested annealing/random samplers.
Unmatched costs; no repeat-based uncertainty.
Packing times use different stopping rounds; no end-to-end timing &
Lower residual in tested settings; no matched-cost or system advantage \\
\addlinespace[3pt]
\textup{[E18]} \citep{lauwerier2026quantum}: tour selection; 10/15/20 vertices; one risk-map date; neutral-atom emulation; preprint &
Independent-set samples / tour generation, conflict graph and weighted selection &
Fixed candidate pools and weights; random, greedy, annealing and exact controls.
Quantum matches restricted optimum only at 15 vertices.
Classical selection faster than projected 100-shot QPU cost; no hardware pipeline timing &
No systematic component superiority; quantum timing projected \\
\end{longtable}\endgroup


%
\begin{table}[tbp]
\centering
\footnotesize
\setlength{\tabcolsep}{3pt}
\setlength{\tablewidth}{\dimexpr\linewidth-18pt\relax}
\begin{tabular}{@{}>{\RaggedRight\arraybackslash}p{0.255\tablewidth}>{\RaggedRight\arraybackslash}p{0.205\tablewidth}>{\RaggedRight\arraybackslash}p{0.375\tablewidth}>{\RaggedRight\arraybackslash}p{0.165\tablewidth}@{}}
\toprule
\textbf{Study and setting} & \textbf{Quantum / classical roles} & \textbf{Design, baseline, metric, and timing boundary} & \textbf{Supported\newline inference} \\
\midrule
\textup{[A1]} \citep{munoz2025scaling}: tailored 2D spin glasses, high-precision couplings; over 1,300 logical qubits &
Error-suppressed annealing / calibration and analysis &
Time-to-\(\varepsilon\) scaling for low-energy sampling (target relative objective gap \(\varepsilon\ge1\%\)) vs. parallel tempering with isoenergetic cluster moves; affirmative scaling claim on this family &
Limited positive system-scaling comparison \\
\addlinespace[3pt]
\textup{[A2]} \citep{pawlowski2026toward}: same benchmark family &
No quantum routine / GPU simulated bifurcation &
Same metric; classical method closes reported scaling gap; accessible sizes insufficient for robust asymptotic inference &
Classical reassessment removes observed scaling separation \\
\addlinespace[3pt]
\textup{[A3]} \citep{chandarana2026quest}: 20 higher-order binary instances; superconducting processor; preprint &
Hybrid solver / pre- and postprocessing within budget &
Matched solver-side time, including preprocessing, quantum execution, postprocessing; orchestration and deployment excluded.
Best hybrid energies over trials/settings match 14 ground states; at mean hybrid runtime, several CPU comparators fall short on average, while enhanced parallel tempering and a GPU solver match or surpass these targets &
Mixed system comparison at solver-side budgets \\
\addlinespace[3pt]
\textup{[A4]} \citep{ebadi2022quantum}: native unit-disk independent set, 39--289 atoms &
Analog evolution / pulse optimization and readout &
Hardware vs. optimized simulated annealing; superlinear speedup on hardest native instances in deep-circuit regime &
Limited positive system comparison \\
\addlinespace[3pt]
\textup{[A5]} \citep{andrist2023hardness}: unit-disk family reassessed &
No quantum routine / exact and heuristic classical solvers &
Structured quasiplanar instances with thousands of nodes solved to optimality in minutes; less constrained annealing competitive on studied family &
Classical competitiveness on reassessed family \\
\addlinespace[3pt]
\textup{[A6]} \citep{cazals2026identifying}: hard native unit-disk instances, varied density and treewidth &
Small-scale neutral-atom evolution / exact and heuristic baselines &
Quantum slower on small instances compared exactly.
Classical scaling suggests about 1,000 atoms at kilohertz repetition as a necessary step toward advantage &
Negative small-scale system comparison; larger-scale outcome unresolved \\
\addlinespace[3pt]
\textup{[A7]} \citep{lu2026evidence}: positive one-in-three SAT, clause-to-variable ratio \(m/n=0.626\); fits to \(n=70\); 13-qubit hardware, 22-variable instances &
Enhanced QAOA and adiabatic solvers / restricting-space reduction, with information shared to classical controls &
Original vs. reduced quantum vs. informed classical solvers.
Scaling about \(1.0077^n\), \(1.0070^n\) vs. \(1.0128^n\).
Quantum cost: inverse success probability; classical: conflict and propagation counts.
Predominantly simulated co-design evidence, not a common-resource interaction contrast or matched wall clock &
Favorable co-design evidence; common-resource interaction not established \\
\bottomrule
\end{tabular}
\caption{Continuation of Table~\ref{tab:controlled-comparisons}: scaling claims, classical reassessments, and hardness-aware benchmark designs.
Positive comparisons retain their instance, metric, and baseline restrictions; co-design evidence does not establish common-resource interaction or system superiority.}\label{tab:scaling-reassessments}
\end{table}



\FloatBarrier
%
\section{Design lessons and outlook}\label{sec:outlook}

The surveyed methods suggest three concrete questions for constraint programming and operations research.

\paragraph{Certified progress.}
When should sampling yield to classical certification?
Benders needs lower bounds, column generation needs final pricing, and copositive separation must distinguish a negative witness from certified absence.
Schedule certification adaptively, comparing its cost with progress from additional samples.
Logic-based Benders and inference duals offer a language for specifying the proofs required at these interfaces.

\paragraph{Constraint-aware reduction.}
Which reductions preserve feasibility or an optimum, and which need repair, backtracking, or certified implications?
Constraint-preserving mixers, classical masters, Lagrangian interfaces, and projection methods offer alternatives to large static penalties; their comparison should include encoding cost, accessibility of useful feasible states, and the consequences of implementation error.

\paragraph{Cost-aware interface selection.}
Is the quantum task still a bottleneck after preprocessing, communication, verification, and repeated calls?
Fixed-wrapper substitutions, independently tuned complete systems, and four-cell reduction comparisons answer complementary questions.
Specify the returned object, target, comparator, and accounting boundary before interpreting a gain.

\paragraph{Community contributions.}

Constraint programming and operations research can contribute established empirical discipline.
The DIMACS implementation challenges \citep{johnson1996dimacs}, the international SAT solver competitions \citep{jarvisalo2012sat}, and the MiniZinc Challenge \citep{stuckey2014minizinc} compare solvers on shared instances under declared resource limits.
Curated libraries such as \software{MIPLIB} \citep{gleixner2021miplib} and \software{QPLIB} \citep{furini2019qplib} supply instances with documented reference solutions and bounds.
\software{QOBLIB} \citep{koch2026qoblib} extends this effort to quantum optimization.
These communities can contribute certified benchmarks, competition-style evaluations with complete resource accounting, and exact, bounding, and inference methods that certify feasibility and bound the quality of sampled candidates.
Beyond benchmarking, these communities can help design the quantum computation itself: qubit placement and routing, gate scheduling, variational parameter optimization, and allocation of sampling and classical computation budgets.
Statistical modeling and uncertainty quantification can support device calibration and control, while application expertise can identify useful structure and realistic decision requirements \citep{bernalneira2025whatcanquantum}.


\FloatBarrier
%
\appendix
\begingroup\small
\section{Survey method}\label{app:survey-method}

This purposive survey uses topic queries, citation follow-up, and primary records to select decomposition, certification, and attribution studies through \SurveyEligibilityDate{}; publication status was checked on \SurveyVerificationDate{}.

The supplementary audit appended to this preprint gives methods, access limits, and Tables S1--S7; the source download and ancillary files provide source and audit data.
These cover 30 studies (25 main-table entries), four historical exclusions, 82 discovery decisions, and an inventory of 23 software implementations.
The survey is neither exhaustive nor an experimental reproduction.

\paragraph{Author involvement.}
The author coauthored cited work on classical optimization and quantum interfaces \citep{kronqvist2019convex,kronqvist2026fifty,bernal2022perspectives,quintero2022characterization,bernalneira2025whatcanquantum}, hybrid methods \citep{brown2024copositive,maciejewski2024multilevel,harwood2021formulating,peng2025hybrid,ramesh2026constrainedpreconditioning,makansi2025greedy}, and benchmarking \citep{dupont2025benchmarking,lubinski2024optimization,bernalneira2025benchmarking,koch2026qoblib,montanez2025evaluating}.
The author's software contributions include \software{QUBO.jl} \citep{xavier2026qubojl}, \software{stochastic-benchmark} \citep{stochastic_benchmark_repo}, and \software{QUBODecomposition.jl} \citep{qubodecomposition_software}; the cited CPAIOR tutorial is also the author's work \citep{bernal2024cpaior}.
The same evidentiary criteria apply to these works as to the other sources.
\endgroup


%
\section*{Statements and Declarations}
\paragraph{Competing interests.}
The author declares no financial competing interests.
Cited work coauthored by the author is identified in the author-involvement disclosure in Appendix~\ref{app:survey-method}.

\paragraph{Author contributions.}
David E. Bernal Neira is the sole author and is responsible for the survey's conception, interpretation, and manuscript.

\paragraph{Data availability.}
The study records, exclusions, discovery decisions, source-version checks, and field definitions accompany this preprint in the appended supplementary audit and ancillary files; the latter contain the machine-readable audit and provenance files.

\paragraph{AI assistance.}
OpenAI ChatGPT and Codex assisted with literature and software discovery, drafting, structural revision, technical and citation checks, copyediting, and preparation of LaTeX/TikZ source and diagrams.
Anthropic Claude provided additional source-based critique of the technical explanations and figures.
This assistance does not constitute independent peer review or verification of the cited experiments; the author retains responsibility for the article's content, interpretations, and conclusions.


%
\section*{Acknowledgments}
The author thanks Carleton Coffrin for the invitation to speak at the CPAIOR 2024 Master Class, ``Quantum Computing for CP, AI, and OR, and vice-versa.''
The author also thanks Mark Wallace for the invitation to contribute this article to the associated special issue of \emph{Constraints}.

This material is based upon work supported by the Center for Quantum Technologies under the Industry-University Cooperative Research Center Program at the US National Science Foundation under Grant No. 2224960.
The author also acknowledges support from the Quantum Seed Fund of Purdue University's Colleges of Engineering and Science and from the Davidson School of Chemical Engineering at Purdue University.


\par


{\small\RaggedRight%

}
\end{document}


\maketitle
\noindent Supplement to the manuscript prepared for \emph{Constraints}, CPAIOR 2024 Master Class special issue.
%
Tables~S1 and~S2 audit the 25 main evidence-table studies. Table~S3 preserves 4 historical exclusion records. Table~S4 records empirical source versions and publication-status checks for all 30 study records, including 5 contextual records. Table~S5 records 82 discovery decisions from the targeted literature extension. Table~S6 identifies source versions and passages for all 30 studies, distinguishing retained inspections from selective new checks. A contextual-record section and the complete data dictionary follow the tables.


The evidence-eligibility cutoff is now \SurveyEligibilityDate.
The original cutoff \SurveyOriginalEligibilityDate\ and audit snapshot \SurveyPriorVerificationDate\ remain historical provenance.
Targeted discovery covers \SurveySearchStartDate--\SurveyEligibilityDate.
Publication-status checks were refreshed on \SurveyVerificationDate; this does not imply that every empirical source was reread.
New empirical audits and selective primary-text rereads identify the claims and versions inspected.
The purposive update does not claim exhaustive literature coverage.
The source-location index in Table~\ref{tab:s6} records the passages behind the main comparisons, including targeted corrections checked on October 6, 2026; these do not extend the eligibility cutoff.

Component effects, system effects, and reduction--solver interactions support different inferences.
Each row retains its design, comparator, target quality, and timing limits.
A wrapper benefit or a favorable interaction can coexist with a faster classical system.
A negative or unresolved component comparison remains informative; the audit does not impose one blanket advantage checklist on every study.

The CSV files \nolinkurl{supplementary_studies.csv}, \nolinkurl{supplementary_exclusion_log.csv}, \nolinkurl{supplementary/literature-search.csv}, and \nolinkurl{supplementary/claim-locations.csv} are the source of truth.
The tables and dictionary are rendered deterministically after identifier, field-coverage, inference-label, and source-location checks.
The source download includes the standalone LaTeX source (\texttt{supplementary.tex}) and these records as ancillary files; the ancillary \texttt{README.txt} explains compilation and packaged filenames, and \texttt{SOURCE\_MANIFEST.json} records the original paths and source checksums.
Historical artifact labels are retained unless an explicit correction is recorded; blank artifact identifiers mean that no identifier was captured.
A status check does not certify the underlying empirical claims or archival availability.
No archival DOI has been assigned to this manuscript package by this update.
%
\section{Survey methods and statistical definitions}\label{supp:methods}

\subsection{Search scope, update history, and evidence records}

This is a narrative, purposive survey: studies were selected for their bearing on decomposition and attribution, rather than through exhaustive screening.
Conclusions apply to the studies analyzed.
The search combined Google Scholar and arXiv queries with publisher browsing across APS, Springer, Nature Portfolio, ACM, IEEE, Elsevier, and Taylor \& Francis.
Query families included ``hybrid quantum classical optimization'', ``Benders quantum'', ``column generation quantum annealing'', ``QUBO decomposition'', ``graph coarsening QUBO'', ``copositive optimization quantum'', and ``quantum annealing benchmark''.
The accompanying study audit records search routes and inclusion reasons; the exclusion log records close studies omitted from the evidence synthesis.
These records document a purposive selection, not an exhaustive systematic search.

The original corpus used an eligibility cutoff of \SurveyOriginalEligibilityDate{} and an audit snapshot dated \SurveyPriorVerificationDate{}.
A targeted update searched work first released or substantively revised from \SurveySearchStartDate{} through \SurveyEligibilityDate{}, extending the current eligibility cutoff to the latter date.
Date-filtered topic queries and primary arXiv and publisher records covered decomposition, constrained optimization, preprocessing, controlled sampler comparisons, classical reassessments, and theoretical optimization guarantees.
The update screened new studies as well as journal versions of older preprints; its queries, candidate decisions, and primary URLs accompany the supplement.
This remains a purposive synthesis, with selections based on methodological relevance and extractable evidence rather than exhaustive coverage.

Publication status was checked on \SurveyVerificationDate{}.
First release, inspected empirical version, and current publication status are distinct fields.
The October 4 pass retrieved arXiv lists ending with the October 2 listing; later releases may not then have been exposed.
The October 12 refresh used date-qualified web searches over the same topic families, primary arXiv listings through October 12, and primary publication records; it did not repeat a Google Scholar search.
The incremental window is October 5--12, including papers submitted earlier but first announced in this interval.
Submission and listing dates are retained separately, and an October journal issue does not by itself establish an October first release.
The update therefore documents what was found by the stated searches, not a complete census through the cutoff.
Earlier empirical records are retained unless a specific reread is recorded; refreshing metadata alone does not revalidate every result.

The SaReduMIS comparison in \citet{schuetz2026qredumis}, first included in \href{https://arxiv.org/html/2503.12551v2}{arXiv version 2}, and the expanded tabu control in \citet{lozano2026where}, from version 3, postdate the original June cutoff but lie inside the extended scope.
The graph-decomposition experiments of \citet{pennington2025boosting} are read from version 2.
Version-specific citations identify inspected evidence when a later version could not be read.
Patents and vendor documents describe the commercial record and are not treated as empirical evidence.

Selective October 6, 2026 claim checks were recorded independently of the publication-status refresh.
For [E11], \href{https://arxiv.org/html/2412.15665v2}{Wagner and Liers, arXiv version 2}, Section~3, distinguishes the proposed complete framework from the implemented master relaxation; Section~5 and Figures~5--10 report pricing and separation quality and timing.
For [A7], the journal version of \citet{lu2026evidence}, Table~1 and \emph{Scalings of time complexities}, gives the reported scaling bases and their distinct quantum and classical cost measures; Figure~5 and \emph{Experiment on superconducting quantum processors} locate the 13-qubit, 22-variable demonstration.
The CSV records these local rereads separately from the corpus eligibility and publication-status dates.

A study enters the evidence tables when its comparative design bears directly on component, system, or interaction effects.
Mechanism papers and contextual counterpoints may be discussed without a main-table row.
Labels [E1]--[E18] and [A1]--[A7] link the main article's controlled-comparison and reassessment tables to \texttt{supplementary\_studies.csv}.
Supplementary Tables S1 and S2 report the associated scale, resource accounting, baseline, uncertainty, artifact, and certificate fields; S3 reproduces all four historical exclusion records; S4 records source versions and publication-status checks for all 30 study records, including five contextual studies.
S5 records the extension's candidate decisions, including mechanism papers without evidence-table rows.
The supplement and \texttt{supplementary\_data\_dictionary.md} define the fields.
The publication-status check log provides the primary URLs and the resolution of each check.
These files accompany the article; no archival deposit identifier is claimed before one exists.

Targeted October 7--12 checks preceded the global refresh and clarified background, theoretical guarantees, benchmark targets, and comparator roles.
The October 12 refresh then checked all 28 existing study records and four historical exclusions, plus 27 other preprints (two overlap the exclusions) and three additional published or correction-sensitive references.
The record of each check states its outcome and access limits; no new version or correction was located among the 32 existing study/exclusion records.
The bibliography now cites the published VQEC journal article and Peng conference paper; these metadata upgrades do not establish stronger algorithmic or empirical guarantees.
Newly selected comparisons have separate full-text extractions, while the prior empirical inspections, exclusion decisions, and 50 October 4 discovery records remain historical provenance.
The ancillary files include the incremental queries and decisions, status-check reports, and source-location history in Table S6.

\subsection{Comparison design and resource accounting}

The main article's attribution definitions follow established benchmarking cautions \citep{ronnow2014defining}.
A \emph{strong classical comparator} means a relevant method identified within this survey's scope at comparable target quality, including competitive GPU or physics-inspired implementations where applicable.
It is a search-relative designation, not a claim of exhaustive baseline coverage.
\emph{End-to-end time} includes preprocessing, reformulation, embedding or compilation, queueing, communication, device execution or search, decoding, verification, and postprocessing whenever applicable.
Solver-side pipeline time, device execution time, projected time, and nominal operation or sample counts are distinguished from that complete wall clock.

The main article defines component and system effects.
A fixed-wrapper swap estimates a component effect but may disadvantage a method whose wrapper was tuned to its original solver; system comparisons allow each pipeline its own tuning.
To quantify a reduction--solver interaction, let \(Y(r,\sigma)\) be a predeclared higher-is-better solution-quality statistic at one fixed resource budget, with \(r\in\{0,1\}\) indicating absence or presence of a reduction and \(\sigma\in\{\mathrm{Q},\mathrm{C}\}\) selecting the quantum or classical solver, respectively.
All four cells use the same budget and resource-accounting boundary; achieved quality is measured, not fixed.
The simple effects are \(\Delta_{\mathrm{Q}}=Y(1,\mathrm{Q})-Y(0,\mathrm{Q})\) and \(\Delta_{\mathrm{C}}=Y(1,\mathrm{C})-Y(0,\mathrm{C})\); the main article defines their contrast as \(\Delta_{\mathrm{int}}=\Delta_{\mathrm{Q}}-\Delta_{\mathrm{C}}\).
A positive contrast means a more favorable relative effect for the quantum solver, while an absolute improvement additionally requires \(\Delta_{\mathrm{Q}}>0\).
For example, \(\Delta_{\mathrm{Q}}=-0.1\) and \(\Delta_{\mathrm{C}}=-0.3\) give a positive contrast of \(0.2\) although both solvers worsen.
Report both simple effects with the contrast.

Alternatively, at a fixed quality target, let \(T(r,\sigma)\) measure lower-is-better cost in the same resource units and accounting boundary for all four cells; cost is measured, not fixed.
Use

\begin{equation}\label{eq:interaction-cost}
\Delta^{\mathrm{cost}}_{\mathrm{int}}=
\bigl[T(0,\mathrm{Q})-T(1,\mathrm{Q})\bigr]-\bigl[T(0,\mathrm{C})-T(1,\mathrm{C})\bigr].
\end{equation}

A positive cost contrast means a more favorable relative cost effect for the quantum solver.
An absolute quantum cost reduction additionally requires \(T(0,\mathrm{Q})-T(1,\mathrm{Q})>0\); report both solvers' cost reductions alongside the contrast.
Neither contrast sign implies that the quantum-enabled system outperforms the classical system.
These are factorial contrasts, not automatically causal estimates: assignment, tuning, and resource comparability determine that interpretation.
The statistic may be an expectation over randomness, an instance average, a median, or a per-instance quantity, but these choices are not interchangeable for heavy-tailed distributions.
With cell estimates \(\widehat\mu^Y_{r,\sigma}\) of the quality statistic \(Y\), estimate

\begin{equation}\label{eq:interaction-estimator}
\widehat\Delta_{\mathrm{int}}=
\bigl(\widehat\mu^Y_{1,\mathrm{Q}}-\widehat\mu^Y_{0,\mathrm{Q}}\bigr)
-\bigl(\widehat\mu^Y_{1,\mathrm{C}}-\widehat\mu^Y_{0,\mathrm{C}}\bigr).
\end{equation}

For cell estimates \(\widehat\mu^T_{r,\sigma}\) of the cost statistic \(T\), reverse the differences:
\[
\widehat\Delta^{\mathrm{cost}}_{\mathrm{int}}=
\bigl(\widehat\mu^T_{0,\mathrm{Q}}-\widehat\mu^T_{1,\mathrm{Q}}\bigr)
-\bigl(\widehat\mu^T_{0,\mathrm{C}}-\widehat\mu^T_{1,\mathrm{C}}\bigr).
\]
Report uncertainty for the declared contrast itself over both instances and solver randomness.
A paired bootstrap resamples four-cell instance records and recomputes the declared statistic and contrast; resampling precomputed instance contrasts suffices for a mean contrast, but not generally for a difference of medians.
With stochastic runs nested in instances, a hierarchical bootstrap also resamples runs within each selected instance, preserving any pairing across cells.
Mixed-effects models are another option.
Time-to-target and gap-at-time profiles \citep{lubinski2024optimization}, together with stochastic-solver scoring \citep{bernalneira2025benchmarking}, make the quality--cost trade explicit.

A benefit is \emph{robust} when it replicates across instances, random seeds, hardware runs, and solver settings with stated uncertainty; it is \emph{transferable} when it persists across generators or independently generated or held-out instances.
Practical scale concerns classical difficulty at the target quality, not variable count alone.
These criteria qualify each effect separately.

\subsection{Finite-grid certificates for copositive separation}\label{supp:copositive-grid}

The main article distinguishes the exact continuous separation oracle of the copositive framework \citep{brown2024copositive} from a finite-bit implementation.
The error estimate below is derived here to explain that distinction; it is not attributed as a theorem of that paper.
For a real symmetric \(d\times d\) matrix \(M\), degree-two homogeneity allows any nonnegative vector with a negative quadratic value to be rescaled into the unit box.
Define
\begin{equation}\label{eq:copositive-box-minimum}
\varphi(M)=\min_{z\in[0,1]^d}z^\top Mz,
\qquad
\varphi_N(M)=\min_{z\in G_N^d}z^\top Mz,
\quad G_N=\{0,1/N,\ldots,1\}.
\end{equation}
Here the positive integer \(N\) is the number of grid intervals per coordinate.
The grid lies in the box and both contain zero, so \(\varphi(M)\le\varphi_N(M)\le0\); copositivity is equivalent to \(\varphi(M)=0\).
For \(N=2^b-1\), a uniform \(b\)-bit encoding uses \(db\) logical bits before embedding or error-control overhead.
Let \(z^\star\) minimize the continuous problem and round each coordinate to its nearest grid value, giving \(\widehat z\) with \(\|\widehat z-z^\star\|_\infty\le1/(2N)\).
Because both vectors lie in the unit box,
\[
\left|\widehat z_i\widehat z_j-z_i^\star z_j^\star\right|
\le |\widehat z_i|\,|\widehat z_j-z_j^\star|
   +|z_j^\star|\,|\widehat z_i-z_i^\star|
\le N^{-1}.
\]
It follows that
\begin{equation}\label{eq:copositive-discretization-error}
0\le\varphi_N(M)-\varphi(M)
\le\varepsilon_N,
\qquad
\varepsilon_N=\frac{1}{N}\sum_{i,j}|M_{ij}|.
\end{equation}
Thus a certified grid optimum \(\varphi_N(M)=0\) implies only
\begin{equation}\label{eq:copositive-grid-bound}
-\varepsilon_N\le\varphi(M)\le0.
\end{equation}
For positive \(\varepsilon_N\), this is approximate nonnegativity on the box, not exact copositivity.
A verified negative grid witness disproves copositivity directly, whereas unsuccessful heuristic sampling supplies neither a certified grid optimum nor this bound.


\clearpage
\begin{landscape}
%
\begingroup
\fontsize{7.5}{9.0}\selectfont
\setlength{\tabcolsep}{3.5pt}
\setlength{\AuditTableWidth}{\dimexpr\linewidth-35pt\relax}
\begin{longtable}{@{}>{\raggedright\arraybackslash}p{0.105\AuditTableWidth}>{\raggedright\arraybackslash}p{0.19\AuditTableWidth}>{\raggedright\arraybackslash}p{0.2\AuditTableWidth}>{\raggedright\arraybackslash}p{0.17\AuditTableWidth}>{\raggedright\arraybackslash}p{0.15\AuditTableWidth}>{\raggedright\arraybackslash}p{0.185\AuditTableWidth}@{}}
\caption{Comparison-study audit, rows [E1]--[E18]. Each supported inference is restricted to its design, resource boundary, and comparator; component, system, and interaction effects are assessed separately.}\label{tab:s1}\\
\toprule
\textbf{Study} & \textbf{Scale and implementation} & \textbf{Resource boundary and target} & \textbf{Baseline and stopping rule} & \textbf{Uncertainty and artifacts} & \textbf{Supported inference and certificate} \\
\midrule\endfirsthead
\multicolumn{6}{l}{\normalsize Table~\thetable\ continued.}\\
\toprule
\textbf{Study} & \textbf{Scale and implementation} & \textbf{Resource boundary and target} & \textbf{Baseline and stopping rule} & \textbf{Uncertainty and artifacts} & \textbf{Supported inference and certificate} \\
\midrule\endhead
\multicolumn{6}{r}{Continues on the next page.}\\\endfoot
\bottomrule\endlastfoot
\textbf{[E1]}\par Ponce et al., 2025 &
\textbf{Scale:} 25 instances; BFGS with 100 restarts per parameter search; 500 shots per hardware instance\par\smallskip \textbf{Hardware / software:} trapped-ion hardware at 500 shots per reduced instance; BFGS with 100 restarts per parameter search &
\textbf{Timing:} no timing claimed\par\smallskip \textbf{Target:} ratio of expected cut value to certified optimum &
Monolithic QAOA and exact component solver; expectation-based evaluation in all arms &
\textbf{Uncertainty:} per-instance ratios tabulated\par\smallskip \textbf{Availability:} not found by documented search\par No artifact identifier captured. &
\textbf{Inference:} Positive wrapper benefit; limited component sensitivity\par\smallskip \textbf{Certificate:} none; ratios vs certified optima \\
\addlinespace[4pt]
\textbf{[E2]}\par Schuetz et al. (2026) &
\textbf{Scale:} 39 instances, 7-231 nodes; 20 classical repetitions; 1,000 shots per call; simulated-annealing arm from a post-cutoff version\par\smallskip \textbf{Hardware / software:} Rydberg hardware; 1000 shots per QPU call; 20 classical repetitions &
\textbf{Timing:} nominal sample count; device cycle-time model reported; no matched end-to-end wall clock\par\smallskip \textbf{Target:} success probability of finding a maximum independent set &
Simulated annealing in the same reduction rule at nominal equal counts; random-informed selection; optimized quantum annealing &
\textbf{Uncertainty:} 90\% confidence intervals via bootstrapping\par\smallskip \textbf{Availability:} source-code repository\par \href{https://github.com/jpmorganchase/qReduMIS}{\nolinkurl{github.com/jpmorganchase/qReduMIS}} (commit 43feb750, inspected July 26, 2026) &
\textbf{Inference:} Positive wrapper benefit; no quantum component superiority\par\smallskip \textbf{Certificate:} No certificate that a sampler-suggested fixing is safe or globally optimal; final independent-set feasibility and cardinality are checkable; exact kernelization rules are independently verifiable \\
\addlinespace[4pt]
\textbf{[E3]}\par Dupont et al., 2025 &
\textbf{Scale:} Thousands of instances; sizes to thousands of variables\par\smallskip \textbf{Hardware / software:} superconducting processor with small subcircuits per call &
\textbf{Timing:} projected time (authors' run-time analysis); device execution time per sampling call\par\smallskip \textbf{Target:} approximation ratio vs classically established references &
Commercial MIP solver, simulated annealing, Burer-Monteiro under the authors' run-time model &
\textbf{Uncertainty:} ensemble averages\par\smallskip \textbf{Availability:} archival dataset\par \href{https://doi.org/10.5281/zenodo.11061744}{\nolinkurl{doi:10.5281/zenodo.11061744}} &
\textbf{Inference:} Negative projected system-speed comparison; attribution unresolved\par\smallskip \textbf{Certificate:} none; heuristic comparison \\
\addlinespace[4pt]
\textbf{[E4]}\par Kim et al., 2025 &
\textbf{Scale:} Over 50 instances (the paper's own count)\par\smallskip \textbf{Hardware / software:} annealer plus proprietary hybrid layer &
\textbf{Timing:} device execution time (QPU access: programming, sampling); problem reading, queueing, and machine-annealer communication excluded\par\smallskip \textbf{Target:} relative accuracy vs best observed objective &
Classical heuristics as black boxes; per-solver settings &
\textbf{Uncertainty:} not found by this survey's search of the study text\par\smallskip \textbf{Availability:} not found by documented search\par No artifact identifier captured. &
\textbf{Inference:} Unresolved attribution; incomplete system timing\par\smallskip \textbf{Certificate:} none; relative to best observed \\
\addlinespace[4pt]
\textbf{[E5]}\par Tasseff et al., 2024 &
\textbf{Scale:} Large sparse instances across classes\par\smallskip \textbf{Hardware / software:} quantum annealer &
\textbf{Timing:} solver-side pipeline time (authors' run-time protocol)\par\smallskip \textbf{Target:} time-to-target style comparison &
Strong classical heuristics; benefits only on specially constructed instance classes, per the authors &
\textbf{Uncertainty:} reported across instance classes\par\smallskip \textbf{Availability:} source-code repository\par No artifact identifier captured. &
\textbf{Inference:} Limited positive system comparison\par\smallskip \textbf{Certificate:} none \\
\addlinespace[4pt]
\textbf{[E6]}\par Vodeb et al., 2026 &
\textbf{Scale:} 139 instances, 100-10,000 nodes\par\smallskip \textbf{Hardware / software:} fast-annealing quantum processor; commercial hybrid solver; simulated-bifurcation machine &
\textbf{Timing:} full wall clock across heterogeneous hardware; some simulated-bifurcation records taken from the literature rather than rerun\par\smallskip \textbf{Target:} Certified Max-Cut optima through 251 vertices (250 original QUBO variables); best-known solution quality above. &
Two simulated-annealing speeds and simulated bifurcation vs the fast-annealing QPU and the hybrid solver &
\textbf{Uncertainty:} Not extracted from the study by this survey's audit\par\smallskip \textbf{Availability:} archival dataset\par \href{https://doi.org/10.5281/zenodo.14290290}{\nolinkurl{doi:10.5281/zenodo.14290290}} &
\textbf{Inference:} Mixed system-quality comparison; attribution unresolved\par\smallskip \textbf{Certificate:} none \\
\addlinespace[4pt]
\textbf{[E7]}\par Lozano, 2026 &
\textbf{Scale:} Portfolios with 10-640 assets; 54 certified head-to-head cases; expanded Tabu control on 36 instances with 108 seeded five-second CPU runs.\par\smallskip \textbf{Hardware / software:} TabuSampler 1.7.0 on Apple M3 Max; proprietary D-Wave hybrid service. Immutable service identifier not re-extracted. &
\textbf{Timing:} Service-reported time and SDK fields; five-second nominal service budget and five-second CPU Tabu runs.\par\smallskip \textbf{Target:} Agreement with Gurobi optima on the certified subset; absolute mean-variance objective gap to hybrid CQM for Tabu controls. &
Gurobi optimality anchors; matched five-second CPU Tabu control, TabuSampler 1.7.0 on Apple M3 Max. Mean absolute gap to CQM 0.00078; maximum 0.00798. &
\textbf{Uncertainty:} 108 Tabu runs: three fixed algorithm seeds per each of 36 instances; mean and maximum absolute gaps reported.\par\smallskip \textbf{Availability:} source-code repository\par No artifact identifier captured. &
\textbf{Inference:} Competitive interface outputs; component attribution unresolved\par\smallskip \textbf{Certificate:} External Gurobi certificates anchor the certified subset; no certificate of causal QPU contribution. \\
\addlinespace[4pt]
\textbf{[E8]}\par Zhu et al., 2026 &
\textbf{Scale:} 20 hardest of 4,000 two-dimensional nearest-neighbor Ising instances; 104 qubits\par\smallskip \textbf{Hardware / software:} 104 superconducting qubits; software versions not extracted by this audit &
\textbf{Timing:} hardware-generated samples evaluated at a sample count derived from projected 0.4 MHz controller throughput (nominal 40 ms budget); time to solution about 2.34x fully projected against single-core simulated annealing; no measured end-to-end timing\par\smallskip \textbf{Target:} frequency of sampling the global minimum &
Tuned single-core simulated annealing; QAOA with local search; no GPU baseline &
\textbf{Uncertainty:} standard deviations over the 20 instances reported\par\smallskip \textbf{Availability:} not found by documented search\par No artifact identifier captured. &
\textbf{Inference:} Limited positive sampling evidence; system timing projected\par\smallskip \textbf{Certificate:} Sampled assignment, energy, and success event are independently checkable; no certificate of global optimality or of the projected timing advantage \\
\addlinespace[4pt]
\textbf{[E9]}\par Pennington et al. (2025) &
\textbf{Scale:} Complete graphs, 6-9 nodes (9-25 qubits), and bipartite graphs, 6-10 nodes (7-24 qubits), state-vector simulation; heavy-hex graphs, 50-100 nodes (52-111 qubits), matrix-product-state simulation and IBM Kingston hardware; preprint (v2)\par\smallskip \textbf{Hardware / software:} Qiskit Aer state-vector simulator (complete graphs 9 to 25 qubits; bipartite graphs 7 to 24 qubits); Qiskit Aer matrix-product-state simulator and IBM Kingston hardware (heavy-hexagonal graphs 52 to 111 qubits); Qiskit Aer 0.16.0; single-layer QAOA with fixed parameters (-0.5, 0.5) on heavy-hexagonal graphs &
\textbf{Timing:} no measured end-to-end timing; sampler-level comparison\par\smallskip \textbf{Target:} decomposition sparsity; approximation error &
Simulated annealing and random sampling, same wrapper; matrix-product-state simulation and superconducting hardware on heavy-hex graphs &
\textbf{Uncertainty:} 10 instances per graph size and topology; reported plots show maximum and minimum extent across instances\par\smallskip \textbf{Availability:} source-code repository\par \href{https://github.com/victor-ibm/QAOA-graph-decomposition}{\nolinkurl{github.com/victor-ibm/QAOA-graph-decomposition}} &
\textbf{Inference:} Positive component benefit in some settings; mixed hardware result\par\smallskip \textbf{Certificate:} sampled matchings checkable for feasibility and objective; no optimality certificate \\
\addlinespace[4pt]
\textbf{[E10]}\par Dupont et al. (2025) &
\textbf{Scale:} Sherrington-Kirkpatrick, random 3-regular maximum-cut, and a grid energy problem; circuit emulation; Rigetti Ankaa-3 hardware test\par\smallskip \textbf{Hardware / software:} circuit emulation of light-cone-induced QAOA at depths p=1 and p=2; Rigetti Ankaa-3 superconducting hardware for the experimental test &
\textbf{Timing:} downstream solver run time to convergence; preconditioning cost accounted in the emulation comparison\par\smallskip \textbf{Target:} convergence speed of the downstream classical solver; solution quality &
Same classical solvers (simulated annealing, Burer-Monteiro) on original vs. preconditioned input &
\textbf{Uncertainty:} averages over 200 randomly generated instances with standard error of the mean for the main random 3-regular maximum-cut experiments\par\smallskip \textbf{Availability:} not found by documented search\par No artifact identifier captured. &
\textbf{Inference:} Positive preconditioning benefit in emulation; hardware system effect unresolved\par\smallskip \textbf{Certificate:} preconditioned objective and downstream solution quality are checkable; no certificate that the quantum-derived preconditioner is optimal or that it yields an end-to-end advantage \\
\addlinespace[4pt]
\textbf{[E11]}\par Wagner and Liers (2026) &
\textbf{Scale:} Vehicle-routing instances; pricing and separation subproblems as QUBOs\par\smallskip \textbf{Hardware / software:} quantum annealer (proof-of-concept study); exact device, solver versions, and parameter protocol not extracted from the study by this survey's audit &
\textbf{Timing:} pricing wall-clock includes model building, embedding, sampling, unembedding and solution construction; master-relaxation and separation timings reported separately; no complete integer-solver end-to-end comparison\par\smallskip \textbf{Target:} negative-reduced-cost routes and exact pricing-call reduction; residual relaxation integrality gap after separation &
Simulated annealing, exact integer-programming pricing and specialized classical separation methods in master-relaxation experiments; classical methods run on a standard laptop &
\textbf{Uncertainty:} Not extracted from the study by this survey's audit\par\smallskip \textbf{Availability:} not found by documented search\par No artifact identifier captured. &
\textbf{Inference:} Negative quantum component comparison\par\smallskip \textbf{Certificate:} Exact pricing certifies the master-relaxation optimum; returned routes and cuts are classically checked. Complete integer optimality additionally requires the proposed branch-price-and-cut tree, which is not fully implemented in these experiments \\
\addlinespace[4pt]
\textbf{[E12]}\par Perron et al. (2026) &
\textbf{Scale:} Maximum weighted independent set, instances to 300 vertices; noiseless emulation; preprint\par\smallskip \textbf{Hardware / software:} noiseless emulation of the neutral-atom sampler &
\textbf{Timing:} no end-to-end timing claim; equal sampling budgets across samplers\par\smallskip \textbf{Target:} solution quality relative to optimality; upper bound from the relaxed linear program &
Direct analog neutral-atom protocols, greedy baselines, and simulated annealing at equal sample counts within the same cutting-plane wrapper &
\textbf{Uncertainty:} Not extracted from the study by this survey's audit\par\smallskip \textbf{Availability:} not found by documented search\par No artifact identifier captured. &
\textbf{Inference:} Wrapper benefit; no quantum component superiority\par\smallskip \textbf{Certificate:} primal independent sets independently checkable; classical polynomial solvability on t-perfect graphs does not by itself establish a polynomial iteration bound for the implemented cut-addition loop; upper bound from the relaxed linear program \\
\addlinespace[4pt]
\textbf{[E13]}\par de Valk et al., 2026 &
\textbf{Scale:} 10 permutations of one 20-customer/4-vehicle instance; five 5-customer/3-vehicle toys; 1000 shots\par\smallskip \textbf{Hardware / software:} MPS-JuliQAOA, bond dimension 64, AMD EPYC 9554P; IBM ibm\_strasbourg; Qiskit; Cbc; HiGHS &
\textbf{Timing:} End-to-end wall time and phase breakdown; hardware component timing only\par\smallskip \textbf{Target:} (UB-LB)/UB; feasibility and cut-cover validity &
Benchmark 1.5 h / 2500 iterations; toys 4 h / 1000 iterations; QAOA COBYLA max 1000 iterations &
\textbf{Uncertainty:} Benchmark medians/interquartile ranges; toy runtime 2900±800 s\par\smallskip \textbf{Availability:} source-code repository\par \href{https://github.com/camilledevalk/hqc-mcms-benders}{\nolinkurl{https://github.com/camilledevalk/hqc-mcms-benders}} &
\textbf{Inference:} Negative result for tested scale\par\smallskip \textbf{Certificate:} Classical bounds; benchmark gap unclosed after 1.5 h; no new convergence theorem \\
\addlinespace[4pt]
\textbf{[E14]}\par Naghmouchi et al., 2026 &
\textbf{Scale:} Two topologies up to 30 nodes / 44 edges and 12 nodes / 23 edges; 25 configurations x 20 instances each; 50 demands\par\smallskip \textbf{Hardware / software:} Pasqal QUBO Solver; state-vector/MPS emulation (MPS above 25 logical variables); CPLEX; no quantum hardware. &
\textbf{Timing:} No wall-clock timing; iterations and quality only\par\smallskip \textbf{Target:} Mean percentage gap to compact-ILP CPLEX optimum. With initialization and refinement, CG-Q means are below 1\%; CG-SA means are 1-3\% or 0-6\% by topology. &
SA tuning/sample-budget parity not established; heuristic failure stops pricing &
\textbf{Uncertainty:} Means and 95\% confidence intervals\par\smallskip \textbf{Availability:} not found by documented search\par No artifact identifier captured. &
\textbf{Inference:} Limited component sensitivity in simulation\par\smallskip \textbf{Certificate:} Heuristic pricing and restricted integer solution; no full-LP/integer certificate \\
\addlinespace[4pt]
\textbf{[E15]}\par Ramesh et al., 2026 &
\textbf{Scale:} 50 instances per size; n=12-40; p=1,2,3 with deeper emulation limited to smaller sizes; edge weights uniform in [0,2].\par\smallskip \textbf{Hardware / software:} Classical analytic p=1 formula; exact state-vector p=2,3; Gurobi 13.0, 8 threads, AMD EPYC 7643; BFGS with 100 starts at n<=20 and angle transfer thereafter. &
\textbf{Timing:} MIP solver time only; correlation construction and quantum parameter optimization excluded. No hardware timing.\par\smallskip \textbf{Target:} Time to an original-objective incumbent within 1\% of the baseline optimum; exact recovery also reported. &
Same Gurobi environment; presolve on/off; relative gap tolerance 1e-6; no wall-time cap. Main scaling uses per-instance oracle penalty; Appendix B checks fixed and held-out selections. &
\textbf{Uncertainty:} Means and standard errors over 50 instances per size; exact-recovery failures retained through penalized average time.\par\smallskip \textbf{Availability:} source-code repository\par \href{https://github.com/SECQUOIA/Quantum_Preconditioning_Constrained_Optimization}{\nolinkurl{https://github.com/SECQUOIA/Quantum_Preconditioning_Constrained_Optimization}} &
\textbf{Inference:} Positive preconditioning benefit; system effect unresolved\par\smallskip \textbf{Certificate:} MIP certificate applies to the transformed objective; finite-depth surrogate solution need not optimize the original objective. Hard balance remains exact. \\
\addlinespace[4pt]
\textbf{[E16]}\par Yalovetzky et al., 2026 &
\textbf{Scale:} Four market graphs with 85, 89, 98, 225 assets and first kernels of 49, 45, 68, 78 qubits; 73 emulator graphs with kernels <=22.\par\smallskip \textbf{Hardware / software:} Quantinuum 98-qubit Helios, p=2; H2-1 noisy state-vector emulator, p=2,6; pretrained/scaled parameters; OR-Tools reference. &
\textbf{Timing:} Cumulative quantum shots, not wall clock. Hardware: 100 standalone shots versus 20 hybrid trials, 10 shots/call, at most 5 calls/trial. Emulator: 500 standalone shots versus 20 hybrid trials, 5 shots/call.\par\smallskip \textbf{Target:} MIS success against OR-Tools reference; independent-set size divided by MIS size; 99\%-target optimal shot-TTS. &
Standalone QAOA uses the identical first kernel, depth and repair; cumulative-shot TTS normalizes unequal raw-success budgets. Separate SA comparator; no classical solver swapped into the full wrapper. &
\textbf{Uncertainty:} 95\% bootstrap intervals from 10,000 resamples; emulator slope standard errors; infinite standalone TTS values excluded from slope fit.\par\smallskip \textbf{Availability:} input-data repository\par \href{https://github.com/AchilleasKn/PortfolioOptimization}{\nolinkurl{https://github.com/AchilleasKn/PortfolioOptimization}} &
\textbf{Inference:} Positive wrapper benefit; component/system superiority unresolved\par\smallskip \textbf{Certificate:} Exact simplicial reductions; sampled frozen-node fixing heuristic. Reference optimum is external; a fully reduced kernel need not give original optimum. \\
\addlinespace[4pt]
\textbf{[E17]}\par Liaqat et al., 2026 &
\textbf{Scale:} 120 vertices/132 edges and PEGASE1354 vertices/1710 edge variables; PEGASE uses 41 partitions with bound 50 qubits, packing four 34/38-qubit circuits into 144 of 156 physical qubits; principal trajectories lack independent-repeat counts\par\smallskip \textbf{Hardware / software:} IBM ibm\_kingston and ibm\_marrakesh, 156-qubit Heron; Qiskit optimization level 3; Qoro Divi; Maestro MPS (chi=128 for PEGASE); exact software version pins not specified in checked passages &
\textbf{Timing:} Packing comparison reports approximately 3× hardware-execution reduction, at different stopping rounds/final errors; no comprehensive end-to-end accounting or system speedup versus classical solvers\par\smallskip \textbf{Target:} Normalized squared Frobenius residual ||D-D*||\_F\textasciicircum{}2/n\textasciicircum{}2 versus retained decomposition length; threshold 1e-6 &
QAOA 25,000 shots versus SA 1,500 reads × 1,000 sweeps per partition/round; top five distinct repaired candidates; residual<=1e-6, relative improvement<=1e-3, or 20 rounds. CPU specification/matched complete cost not established. &
\textbf{Uncertainty:} No repeated-instance/seed uncertainty estimates identified for the principal reported convergence trajectories; 25,000 shots are within-run samples, not independent-instance replications\par\smallskip \textbf{Availability:} available on request\par No artifact identifier captured. &
\textbf{Inference:} Lower residual in tested settings; no matched-cost or system advantage\par\smallskip \textbf{Certificate:} Repair verifies matching feasibility; normalized residual is measurable; no global minimum-cardinality certificate or quantum-advantage guarantee \\
\addlinespace[4pt]
\textbf{[E18]}\par Lauwerier et al., 2026 &
\textbf{Scale:} 10,15,20-node/qubit conflict graphs; Table 3 uses September30,2021 instance;100 emulated output samples; broader daily/monthly risk-map evaluations do not replicate the controlled comparison\par\smallskip \textbf{Hardware / software:} Pasqal maximum-independent-set 0.3.1 and EMU-SV; dual AMD EPYC 7643 (96 cores),500 GB RAM; HiGHS; no physical quantum device &
\textbf{Timing:} Classical selector wall time versus assumed 100-shot/1-Hz QPU time; all quantum work emulated; headline 120 s is not measured end-to-end hardware time\par\smallskip \textbf{Target:} Selected tour reward; feasibility; relative gap to restricted-selection LP bound, which is not the exact restricted optimum or full routing optimum &
Full routing HiGHS comparison allows 3600 s; restricted MILP/heuristic selector timings in Table 3 are distinct. Classical selector CPU/tuning/sample parity not fully specified; default quantum-library configuration,100 samples. &
\textbf{Uncertainty:} Monthly means/minima/maxima for pipeline outputs; no independent-repeat uncertainty for Table 3 fixed-date sampler comparison\par\smallskip \textbf{Availability:} source-code repository\par \href{https://github.com/QUASAR-UdS/TOP/}{\nolinkurl{https://github.com/QUASAR-UdS/TOP/}} &
\textbf{Inference:} No systematic component superiority; quantum timing projected\par\smallskip \textbf{Certificate:} Feasible selections checked; exact/brute-force controls certify restricted candidate-pool optima; LP gap can remain positive at that optimum; no full routing certificate or finite-time candidate-generation guarantee \\
\addlinespace[4pt]
\end{longtable}
\endgroup

\clearpage
\begingroup
\fontsize{7.5}{9.0}\selectfont
\setlength{\tabcolsep}{3.5pt}
\setlength{\AuditTableWidth}{\dimexpr\linewidth-35pt\relax}
\begin{longtable}{@{}>{\raggedright\arraybackslash}p{0.105\AuditTableWidth}>{\raggedright\arraybackslash}p{0.19\AuditTableWidth}>{\raggedright\arraybackslash}p{0.2\AuditTableWidth}>{\raggedright\arraybackslash}p{0.17\AuditTableWidth}>{\raggedright\arraybackslash}p{0.15\AuditTableWidth}>{\raggedright\arraybackslash}p{0.185\AuditTableWidth}@{}}
\caption{Affirmative claims and classical reassessments, rows [A1]--[A7]. Each supported inference is restricted to its design, resource boundary, and comparator; component, system, and interaction effects are assessed separately.}\label{tab:s2}\\
\toprule
\textbf{Study} & \textbf{Scale and implementation} & \textbf{Resource boundary and target} & \textbf{Baseline and stopping rule} & \textbf{Uncertainty and artifacts} & \textbf{Supported inference and certificate} \\
\midrule\endfirsthead
\multicolumn{6}{l}{\normalsize Table~\thetable\ continued.}\\
\toprule
\textbf{Study} & \textbf{Scale and implementation} & \textbf{Resource boundary and target} & \textbf{Baseline and stopping rule} & \textbf{Uncertainty and artifacts} & \textbf{Supported inference and certificate} \\
\midrule\endhead
\multicolumn{6}{r}{Continues on the next page.}\\\endfoot
\bottomrule\endlastfoot
\textbf{[A1]}\par Munoz-Bauza and Lidar, 2025 &
\textbf{Scale:} 2D lattices; over 1,300 logical qubits under quantum annealing correction; scaling fits across sizes\par\smallskip \textbf{Hardware / software:} quantum annealing correction, over 1300 logical qubits &
\textbf{Timing:} device execution time (anneal-time-based time-to-epsilon protocol)\par\smallskip \textbf{Target:} time-to-epsilon at optimality gaps >= 1\% &
Parallel tempering with isoenergetic cluster moves at matched targets &
\textbf{Uncertainty:} scaling fits with stated uncertainties\par\smallskip \textbf{Availability:} archival dataset\par No artifact identifier captured. &
\textbf{Inference:} Limited positive system-scaling comparison\par\smallskip \textbf{Certificate:} none; sampling metric \\
\addlinespace[4pt]
\textbf{[A2]}\par Pawlowski et al., 2026 &
\textbf{Scale:} Same benchmark family; GPU implementation\par\smallskip \textbf{Hardware / software:} GPU implementation &
\textbf{Timing:} solver-side pipeline time (GPU wall clock per the paper)\par\smallskip \textbf{Target:} same metric as the reassessed study &
Same time-to-epsilon style target as the reassessed study &
\textbf{Uncertainty:} scaling fits\par\smallskip \textbf{Availability:} not found by documented search\par No artifact identifier captured. &
\textbf{Inference:} Classical reassessment removes observed scaling separation\par\smallskip \textbf{Certificate:} none \\
\addlinespace[4pt]
\textbf{[A3]}\par Chandarana et al., 2026 &
\textbf{Scale:} 20 instances; repetitions within each budget\par\smallskip \textbf{Hardware / software:} superconducting processor; baselines on 128 vCPUs and 8 GPUs &
\textbf{Timing:} solver-side pipeline time (preprocessing, quantum execution, postprocessing); orchestration and deployment latencies excluded\par\smallskip \textbf{Target:} Target is the best energy observed across hybrid trials and repetition-delay settings for each instance; 14 of 20 targets match independently reported ground states. &
Enhanced parallel tempering and a GPU solver on 128 vCPUs and 8 GPUs at the same budget &
\textbf{Uncertainty:} reported per instance\par\smallskip \textbf{Availability:} not found by documented search\par No artifact identifier captured. &
\textbf{Inference:} Mixed system comparison at solver-side budgets\par\smallskip \textbf{Certificate:} none \\
\addlinespace[4pt]
\textbf{[A4]}\par Ebadi et al., 2022 &
\textbf{Scale:} 39-289 atoms; hardest-instance subset emphasized\par\smallskip \textbf{Hardware / software:} neutral-atom array, 39-289 atoms &
\textbf{Timing:} device execution time (hardware evolution) vs nominal operation count (annealing sweeps)\par\smallskip \textbf{Target:} probability of finding exact maximum independent sets &
Optimized simulated annealing &
\textbf{Uncertainty:} reported across instance classes\par\smallskip \textbf{Availability:} archival dataset\par \href{https://doi.org/10.5281/zenodo.6462687}{\nolinkurl{doi:10.5281/zenodo.6462687}} &
\textbf{Inference:} Limited positive system comparison\par\smallskip \textbf{Certificate:} none \\
\addlinespace[4pt]
\textbf{[A5]}\par Andrist et al., 2023 &
\textbf{Scale:} Instances to thousands of nodes; scaling analysis\par\smallskip \textbf{Hardware / software:} commodity classical hardware &
\textbf{Timing:} solver-side pipeline time (classical run time on commodity hardware)\par\smallskip \textbf{Target:} same maximum-independent-set success metric &
Exact and heuristic solvers, including a less constrained simulated annealing &
\textbf{Uncertainty:} scaling fits reported\par\smallskip \textbf{Availability:} source-code repository\par No artifact identifier captured. &
\textbf{Inference:} Classical competitiveness on reassessed family\par\smallskip \textbf{Certificate:} none \\
\addlinespace[4pt]
\textbf{[A6]}\par Cazals et al., 2026 &
\textbf{Scale:} Small instances compared exactly; larger classical-only scaling runs; hardness set by density and treewidth\par\smallskip \textbf{Hardware / software:} neutral-atom platform on small instances; modern classical solvers &
\textbf{Timing:} solver-side classical run time; projected time for the device from repetition rate\par\smallskip \textbf{Target:} exact-solution time on hardness-controlled instances &
Modern exact and heuristic classical solvers on deliberately hardened native instances &
\textbf{Uncertainty:} reported across the constructed families\par\smallskip \textbf{Availability:} not found by documented search\par No artifact identifier captured. &
\textbf{Inference:} Negative small-scale system comparison; larger-scale outcome unresolved\par\smallskip \textbf{Certificate:} none \\
\addlinespace[4pt]
\textbf{[A7]}\par Lu et al., 2026 &
\textbf{Scale:} Scaling fits to n = 70; simulations to 100-150 variables; 13-qubit hardware on 22-variable instances\par\smallskip \textbf{Hardware / software:} 13-qubit superconducting hardware; 40-layer QAOA; 150-layer quantum adiabatic &
\textbf{Timing:} inverse-success-probability cost at fixed depth (quantum) vs conflict and propagation counters, a nominal operation count (classical); no matched wall clock\par\smallskip \textbf{Target:} satisfying assignment; fitted scaling base uses inverse success probability for quantum solvers and conflicts/propagations for classical solvers &
Selected classical baseline also given reduction-related information; fixed depths (40-layer QAOA, 150-layer adiabatic) &
\textbf{Uncertainty:} confidence intervals reported with the scaling fits\par\smallskip \textbf{Availability:} archival dataset\par \href{https://doi.org/10.5281/zenodo.19970881}{\nolinkurl{doi:10.5281/zenodo.19970881}} &
\textbf{Inference:} Favorable co-design evidence; common-resource interaction not established\par\smallskip \textbf{Certificate:} none \\
\addlinespace[4pt]
\end{longtable}
\endgroup

\clearpage
\begingroup
\fontsize{7.5}{9.0}\selectfont
\setlength{\tabcolsep}{3.5pt}
\setlength{\AuditTableWidth}{\dimexpr\linewidth-21pt\relax}
\begin{longtable}{@{}>{\raggedright\arraybackslash}p{0.23\AuditTableWidth}>{\raggedright\arraybackslash}p{0.22\AuditTableWidth}>{\raggedright\arraybackslash}p{0.29\AuditTableWidth}>{\raggedright\arraybackslash}p{0.26\AuditTableWidth}@{}}
\caption{All 4 historical exclusion records. These are purposive selection decisions, not a requirement that every included study demonstrate an end-to-end quantum advantage. Current discovery decisions are recorded separately in Table~S5.}\label{tab:s3}\\
\toprule
\textbf{Study and identifier} & \textbf{Why considered} & \textbf{Recorded exclusion and disposition} & \textbf{Publication-status check} \\
\midrule\endfirsthead
\multicolumn{4}{l}{\normalsize Table~\thetable\ continued.}\\
\toprule
\textbf{Study and identifier} & \textbf{Why considered} & \textbf{Recorded exclusion and disposition} & \textbf{Publication-status check} \\
\midrule\endhead
\multicolumn{4}{r}{Continues on the next page.}\\\endfoot
\bottomrule\endlastfoot
\textbf{[X1]} Mendes de Araujo and de Abreu Faria, Hybrid Quantum Neighborhood Selection\par\smallskip \href{https://arxiv.org/abs/2606.29647}{\nolinkurl{arXiv:2606.29647}}\par\smallskip \textbf{First version:} June 28, 2026 &
frontier-size, warm-start, filtering, and rotation ablations relevant to reduction-solver interaction; eligible by first-version date &
very recent single-group simulation preprint; mentioned in Section 4.2 of the main article and treated cautiously per review guidance; excluded from the Section 5 evidence synthesis\par\smallskip \textbf{Current disposition:} Historical reason retained as provenance. The separate July HQNS companion is eligible under the expanded cutoff but remains screen-only: execution/timing boundaries and budget matching are not established, Section V-H does not name the backend, and no matched frontier classical control was identified. Date eligibility is no longer the reason for omission. &
\textbf{Current status:} Preprint; no journal publication identified in checked primary records\par\smallskip \textbf{Check:} verified primary record (2026-10-12)\par\smallskip October 12 metadata refresh: arXiv title/authors/submission history rechecked; latest v1 (2026-06-28) unchanged. No journal reference appears on the record, and the recorded title search located no matching journal publication. This bounded search is not proof of nonpublication. Existing empirical source versions and audit dates are retained; no empirical reread. The separately identified HQNS companion remains arXiv:2607.07336v1; its existing screen-only disposition is unchanged.\par \url{https://arxiv.org/abs/2606.29647}\par \url{https://arxiv.org/abs/2607.07336} \\
\addlinespace[4pt]
\textbf{[X2]} Lozano, companion direct-QPU pipeline study\par\smallskip \href{https://arxiv.org/abs/2605.17628}{\nolinkurl{arXiv:2605.17628}}\par\smallskip \textbf{First version:} May 2026 &
companion to the operational audit [E7] on the same instance family &
formulation and direct-QPU pipeline study rather than a hybrid decomposition audit; the companion audit [E7] carries the evidence relevant to this survey\par\smallskip \textbf{Current disposition:} Historical omission retained after the targeted update: this formulation/direct-QPU companion adds no distinct decomposition comparison beyond the operational audit [E7]. This is a purposive coverage decision. &
\textbf{Current status:} preprint; no journal version located in this status refresh\par\smallskip \textbf{Check:} verified primary record (2026-10-12)\par\smallskip October 12 metadata refresh: arXiv title/authors/submission history rechecked; latest v2 (2026-05-31) unchanged. No journal reference appears on the record, and the recorded title search located no matching journal publication. This bounded search is not proof of nonpublication. Existing empirical source versions and audit dates are retained; no empirical reread.\par \url{https://arxiv.org/abs/2605.17628} \\
\addlinespace[4pt]
\textbf{[X3]} Sharma and Lau, adaptive constraint-aware graph shrinking with verification and repair\par\smallskip \href{https://arxiv.org/abs/2506.14250}{\nolinkurl{arXiv:2506.14250}}\par\smallskip \textbf{First version:} June 2025 &
classical constraint-aware shrinking and repair pipeline in the reduction family, directly relevant to the attribution thesis &
simulation-only variational experiments without a hardware demonstration; no comparative design isolating a marginal quantum contribution under the audit fields; noted in Section 4 of the main article\par\smallskip \textbf{Current disposition:} Retained as a contextual mechanism citation without a main evidence row; the inspected comparison does not isolate a marginal quantum contribution. Simulation alone is not a general bar to inclusion. &
\textbf{Current status:} preprint; no journal version located in this status refresh\par\smallskip \textbf{Check:} verified primary record (2026-10-12)\par\smallskip October 12 metadata refresh: arXiv title/authors/submission history rechecked; latest v1 (2025-06-17) unchanged. No journal reference appears on the record, and the recorded title search located no matching journal publication. This bounded search is not proof of nonpublication. Existing empirical source versions and audit dates are retained; no empirical reread.\par \url{https://arxiv.org/abs/2506.14250} \\
\addlinespace[4pt]
\textbf{[X4]} Holliday, Osaba, and Luu, hybrid quantum tabu search for capacitated vehicle routing with varied quantum-invocation frequency\par\smallskip \href{https://doi.org/10.1007/s42484-026-00375-8}{\nolinkurl{doi:10.1007/s42484-026-00375-8}}\par\smallskip \textbf{First version:} January 2025 (arXiv); journal 2026 &
quantum-invocation-frequency ablation inside a classical tabu wrapper, relevant to component-frequency attribution &
varies quantum-invocation frequency and compares hybrid variants, a frequency ablation relevant to component attribution, but lacks a matched zero-quantum fixed-wrapper arm and complete resource accounting, so it does not estimate a marginal quantum contribution; noted in Section 5 of the main article\par\smallskip \textbf{Current disposition:} Frequency-ablation findings remain discussed contextually; unmatched zero-quantum controls and incomplete resource accounting limit their marginal-component interpretation. No general exclusion of system or interaction evidence is implied. &
\textbf{Current status:} Published journal article\par\smallskip \textbf{Check:} verified primary record (2026-10-12)\par\smallskip October 12 metadata refresh: published status and identifier retained; Quantum Machine Intelligence 8, 37. Crossmark reports the document current with no listed update. Existing empirical source versions and audit dates are retained; no empirical reread.\par \url{https://link.springer.com/article/10.1007/s42484-026-00375-8}\par \url{https://crossmark.crossref.org/dialog/?doi=10.1007/s42484-026-00375-8}\par \url{https://arxiv.org/abs/2501.12652} \\
\addlinespace[4pt]
\end{longtable}
\endgroup

\clearpage
\begingroup
\fontsize{7.5}{9.0}\selectfont
\setlength{\tabcolsep}{3.5pt}
\setlength{\AuditTableWidth}{\dimexpr\linewidth-28pt\relax}
\begin{longtable}{@{}>{\raggedright\arraybackslash}p{0.1\AuditTableWidth}>{\raggedright\arraybackslash}p{0.2\AuditTableWidth}>{\raggedright\arraybackslash}p{0.19\AuditTableWidth}>{\raggedright\arraybackslash}p{0.21\AuditTableWidth}>{\raggedright\arraybackslash}p{0.3\AuditTableWidth}@{}}
\caption{Version and publication-status audit for all 30 study records. Original cutoff: \SurveyOriginalEligibilityDate; historical audit: \SurveyPriorVerificationDate; expanded eligibility: \SurveyEligibilityDate; current checks: \SurveyVerificationDate. A metadata refresh does not imply empirical rereading; new empirical audits and selective rereads are identified explicitly.}\label{tab:s4}\\
\toprule
\textbf{Study} & \textbf{Empirical source inspected} & \textbf{Historical publication status} & \textbf{Publication record at refresh} & \textbf{Evidence and limits of the check} \\
\midrule\endfirsthead
\multicolumn{5}{l}{\normalsize Table~\thetable\ continued.}\\
\toprule
\textbf{Study} & \textbf{Empirical source inspected} & \textbf{Historical publication status} & \textbf{Publication record at refresh} & \textbf{Evidence and limits of the check} \\
\midrule\endhead
\multicolumn{5}{r}{Continues on the next page.}\\\endfoot
\bottomrule\endlastfoot
\textbf{[E1]}\par Ponce et al., 2025 &
\textbf{Source:} journal version\par\smallskip \textbf{Audit venue:} Quantum Information Processing\par\smallskip \href{https://doi.org/10.1007/s11128-025-04675-z}{\nolinkurl{doi:10.1007/s11128-025-04675-z}} &
\textbf{At original cutoff:} journal article\par\smallskip \textbf{At historical audit:} journal article &
\textbf{Current record:} Published journal article\par Quantum Information Processing\par \href{https://doi.org/10.1007/s11128-025-04675-z}{\nolinkurl{doi:10.1007/s11128-025-04675-z}} &
\textbf{Check result:} verified primary record\par\smallskip \textbf{Check date:} 2026-10-12\par\smallskip October 12 metadata refresh: published status and identifier retained; Quantum Information Processing 24, 60. Crossmark reports the document current with no listed update. Existing empirical source versions and audit dates are retained; no empirical reread. Publisher abstract and bibliographic record accessible; article body is subscription-limited. No empirical reread.\par \url{https://link.springer.com/article/10.1007/s11128-025-04675-z}\par \url{https://crossmark.crossref.org/dialog/?doi=10.1007/s11128-025-04675-z}\par \url{https://api.crossref.org/works/10.1007/s11128-025-04675-z} \\
\addlinespace[4pt]
\textbf{[E2]}\par Schuetz et al. (2026) &
\textbf{Source:} arXiv v2 (2026-07-01), with central SaReduMIS claims selectively rechecked against journal version of record (2026-09-10) on 2026-10-04\par\smallskip \textbf{Audit venue:} Physical Review Research\par\smallskip \href{https://doi.org/10.1103/nh3b-1wv5}{\nolinkurl{doi:10.1103/nh3b-1wv5}} &
\textbf{At original cutoff:} preprint (v1, March 16, 2025)\par\smallskip \textbf{At historical audit:} accepted article, Physical Review Research (accepted July 20, 2026) &
\textbf{Current record:} published journal article\par Physical Review Research\par \href{https://doi.org/10.1103/nh3b-1wv5}{\nolinkurl{doi:10.1103/nh3b-1wv5}}\par \textbf{Latest observed:} journal version of record (2026-09-10); latest arXiv v3 (2026-07-10) &
\textbf{Check result:} verified primary record\par\smallskip \textbf{Check date:} 2026-10-12\par\smallskip October 12 metadata refresh: published status and identifier retained; Physical Review Research 8, 033296. Crossmark reports the document current with no listed update. Existing empirical source versions and audit dates are retained; no empirical reread.\par \url{https://journals.aps.org/prresearch/abstract/10.1103/nh3b-1wv5}\par \url{https://crossmark.crossref.org/dialog/?doi=10.1103/nh3b-1wv5}\par \url{https://arxiv.org/abs/2503.12551}\par \textbf{Empirical inspection (2026-10-04):} Selective verification against journal PDF pp. 6 and 8; historical v2 provenance retained, published support added. No full audit of all numerical fields or code artifact. Source: https://journals.aps.org/prresearch/pdf/10.1103/nh3b-1wv5 \\
\addlinespace[4pt]
\textbf{[E3]}\par Dupont et al., 2025 &
\textbf{Source:} journal version\par\smallskip \textbf{Audit venue:} journal article (superconducting hardware benchmark)\par\smallskip \href{https://doi.org/10.1103/PhysRevApplied.23.014045}{\nolinkurl{doi:10.1103/PhysRevApplied.23.014045}} &
\textbf{At original cutoff:} journal article\par\smallskip \textbf{At historical audit:} journal article &
\textbf{Current record:} Published journal article\par Physical Review Applied\par \href{https://doi.org/10.1103/PhysRevApplied.23.014045}{\nolinkurl{doi:10.1103/PhysRevApplied.23.014045}} &
\textbf{Check result:} verified primary record\par\smallskip \textbf{Check date:} 2026-10-12\par\smallskip October 12 metadata refresh: published status and identifier retained; Physical Review Applied 23, 014045. Crossmark reports the document current with no listed update. Existing empirical source versions and audit dates are retained; no empirical reread. Publisher bibliographic record accessible; article body is subscription-limited. No empirical reread.\par \url{https://journals.aps.org/prapplied/abstract/10.1103/PhysRevApplied.23.014045}\par \url{https://crossmark.crossref.org/dialog/?doi=10.1103/PhysRevApplied.23.014045}\par \url{https://api.crossref.org/works/10.1103/PhysRevApplied.23.014045} \\
\addlinespace[4pt]
\textbf{[E4]}\par Kim et al., 2025 &
\textbf{Source:} journal version\par\smallskip \textbf{Audit venue:} journal article (annealer application)\par\smallskip \href{https://doi.org/10.1038/s41534-025-01020-1}{\nolinkurl{doi:10.1038/s41534-025-01020-1}} &
\textbf{At original cutoff:} journal article\par\smallskip \textbf{At historical audit:} journal article &
\textbf{Current record:} Published journal article\par npj Quantum Information\par \href{https://doi.org/10.1038/s41534-025-01020-1}{\nolinkurl{doi:10.1038/s41534-025-01020-1}} &
\textbf{Check result:} verified primary record\par\smallskip \textbf{Check date:} 2026-10-12\par\smallskip October 12 metadata refresh: published status and identifier retained; npj Quantum Information 11, 77. Crossmark reports the document current with no listed update. Existing empirical source versions and audit dates are retained; no empirical reread. Publisher HTML returned an internal fetch error; publisher PDF page 1 and Crossmark were accessible. Metadata and update status verified, not empirical results.\par \url{https://www.nature.com/articles/s41534-025-01020-1.pdf}\par \url{https://crossmark.crossref.org/dialog/?doi=10.1038/s41534-025-01020-1} \\
\addlinespace[4pt]
\textbf{[E5]}\par Tasseff et al., 2024 &
\textbf{Source:} journal version\par\smallskip \textbf{Audit venue:} journal article\par\smallskip \href{https://doi.org/10.1007/s10732-024-09530-5}{\nolinkurl{doi:10.1007/s10732-024-09530-5}} &
\textbf{At original cutoff:} journal article\par\smallskip \textbf{At historical audit:} journal article &
\textbf{Current record:} Published journal article\par Journal of Heuristics\par \href{https://doi.org/10.1007/s10732-024-09530-5}{\nolinkurl{doi:10.1007/s10732-024-09530-5}} &
\textbf{Check result:} verified primary record\par\smallskip \textbf{Check date:} 2026-10-12\par\smallskip October 12 metadata refresh: published status and identifier retained; Journal of Heuristics 30, 325--358. Crossmark reports the document current with no listed update. Existing empirical source versions and audit dates are retained; no empirical reread.\par \url{https://link.springer.com/article/10.1007/s10732-024-09530-5}\par \url{https://crossmark.crossref.org/dialog/?doi=10.1007/s10732-024-09530-5}\par \url{https://api.crossref.org/works/10.1007/s10732-024-09530-5} \\
\addlinespace[4pt]
\textbf{[E6]}\par Vodeb et al., 2026 &
\textbf{Source:} journal version (July 3, 2026)\par\smallskip \textbf{Audit venue:} Quantum Information Processing\par\smallskip \href{https://doi.org/10.1007/s11128-026-05263-5}{\nolinkurl{doi:10.1007/s11128-026-05263-5}} &
\textbf{At original cutoff:} preprint (arXiv:2412.07460, December 2024)\par\smallskip \textbf{At historical audit:} journal article (July 3, 2026) &
\textbf{Current record:} published journal article\par Quantum Information Processing\par \href{https://doi.org/10.1007/s11128-026-05263-5}{\nolinkurl{doi:10.1007/s11128-026-05263-5}}\par \textbf{Latest observed:} journal version of record (2026-07-03); arXiv v1 (2024-12-10) &
\textbf{Check result:} verified primary record\par\smallskip \textbf{Check date:} 2026-10-12\par\smallskip October 12 metadata refresh: published status and identifier retained; Quantum Information Processing 25, 240. Crossmark reports the document current with no listed update. Existing empirical source versions and audit dates are retained; no empirical reread.\par \url{https://link.springer.com/article/10.1007/s11128-026-05263-5}\par \url{https://crossmark.crossref.org/dialog/?doi=10.1007/s11128-026-05263-5}\par \url{https://arxiv.org/abs/2412.07460}\par \url{https://api.crossref.org/works/10.1007/s11128-026-05263-5}\par \textbf{Empirical inspection (2026-10-04):} Selectively reread Sections 4.1, 4.2 and 5, resolving QUBO-variable versus Max-Cut-vertex counts and confirming heterogeneous timing caution. No full re-audit of all instance results. Source: https://link.springer.com/article/10.1007/s11128-026-05263-5 \\
\addlinespace[4pt]
\textbf{[E7]}\par Lozano, 2026 &
\textbf{Source:} arXiv v3 (2026-07-12); post-cutoff expanded Tabu control selectively rechecked on 2026-10-04\par\smallskip \textbf{Audit venue:} arXiv\par\smallskip \href{https://arxiv.org/abs/2605.17623}{\nolinkurl{arXiv:2605.17623}} &
\textbf{At original cutoff:} preprint (v1 May 17, 2026; v2 June 2, 2026)\par\smallskip \textbf{At historical audit:} preprint &
\textbf{Current record:} preprint; no journal version located in this status refresh\par arXiv\par \href{https://arxiv.org/abs/2605.17623}{\nolinkurl{arXiv:2605.17623}}\par \textbf{Latest observed:} v3 (2026-07-12) &
\textbf{Check result:} verified primary record\par\smallskip \textbf{Check date:} 2026-10-12\par\smallskip October 12 metadata refresh: arXiv title/authors/submission history rechecked; latest v3 (2026-07-12) unchanged. No journal reference appears on the record, and the recorded title search located no matching journal publication. This bounded search is not proof of nonpublication. Existing empirical source versions and audit dates are retained; no empirical reread.\par \url{https://arxiv.org/abs/2605.17623}\par \textbf{Empirical inspection (2026-10-04):} Selectively checked v3 abstract and matched-budget subsection/Table 8 against v2 abstract. Update version provenance and label post-cutoff expanded control; do not claim all service or artifact fields reaudited. Source: https://arxiv.org/html/2605.17623v3 \\
\addlinespace[4pt]
\textbf{[E8]}\par Zhu et al., 2026 &
\textbf{Source:} journal version (March 2, 2026)\par\smallskip \textbf{Audit venue:} National Science Review\par\smallskip \href{https://doi.org/10.1093/nsr/nwag124}{\nolinkurl{doi:10.1093/nsr/nwag124}} &
\textbf{At original cutoff:} journal article (March 2, 2026, inside the cutoff)\par\smallskip \textbf{At historical audit:} journal article &
\textbf{Current record:} Published journal article\par National Science Review\par \href{https://doi.org/10.1093/nsr/nwag124}{\nolinkurl{doi:10.1093/nsr/nwag124}}\par \textbf{Latest observed:} 2509.11535v1 (2025-09-15) &
\textbf{Check result:} verified primary record\par\smallskip \textbf{Check date:} 2026-10-12\par\smallskip October 12 metadata refresh: published status and identifier retained; National Science Review 13(9), nwag124. Crossmark has no data; the recorded DOI correction search located no correction. Existing empirical source versions and audit dates are retained; no empirical reread. Publisher HTML redirected to an inaccessible Silverchair endpoint. Publisher-deposited Crossref metadata and arXiv DOI link confirm publication; Crossmark has no data. DOI/title correction search located no correction, but absence is not proved.\par \url{https://api.crossref.org/works/10.1093/nsr/nwag124}\par \url{https://crossmark.crossref.org/dialog/?doi=10.1093/nsr/nwag124}\par \url{https://arxiv.org/abs/2509.11535} \\
\addlinespace[4pt]
\textbf{[E9]}\par Pennington et al. (2025) &
\textbf{Source:} arXiv v2 (November 14, 2025)\par\smallskip \textbf{Audit venue:} arXiv\par\smallskip \href{https://arxiv.org/abs/2509.10657}{\nolinkurl{arXiv:2509.10657}} &
\textbf{At original cutoff:} preprint (v2, November 14, 2025)\par\smallskip \textbf{At historical audit:} preprint; no peer-reviewed version located as of July 26, 2026 &
\textbf{Current record:} preprint; no journal version located in this status refresh\par arXiv\par \href{https://arxiv.org/abs/2509.10657}{\nolinkurl{arXiv:2509.10657}}\par \textbf{Latest observed:} v2 (2025-11-14) &
\textbf{Check result:} verified primary record\par\smallskip \textbf{Check date:} 2026-10-12\par\smallskip October 12 metadata refresh: arXiv title/authors/submission history rechecked; latest v2 (2025-11-14) unchanged. No journal reference appears on the record, and the recorded title search located no matching journal publication. This bounded search is not proof of nonpublication. Existing empirical source versions and audit dates are retained; no empirical reread.\par \url{https://arxiv.org/abs/2509.10657} \\
\addlinespace[4pt]
\textbf{[E10]}\par Dupont et al. (2025) &
\textbf{Source:} journal version (October 3, 2025)\par\smallskip \textbf{Audit venue:} Physical Review Applied\par\smallskip \href{https://doi.org/10.1103/9prw-684p}{\nolinkurl{doi:10.1103/9prw-684p}} &
\textbf{At original cutoff:} journal article (October 3, 2025, inside the cutoff)\par\smallskip \textbf{At historical audit:} journal article &
\textbf{Current record:} Published journal article\par Physical Review Applied\par \href{https://doi.org/10.1103/9prw-684p}{\nolinkurl{doi:10.1103/9prw-684p}}\par \textbf{Latest observed:} 2502.18570v2 (2025-10-07) &
\textbf{Check result:} verified primary record\par\smallskip \textbf{Check date:} 2026-10-12\par\smallskip October 12 metadata refresh: published status and identifier retained; Physical Review Applied 24, 044013. Crossmark reports the document current with no listed update. Existing empirical source versions and audit dates are retained; no empirical reread. Publisher abstract request timed out; arXiv journal reference and Crossmark confirm the published record. Crossref metadata was successfully retrieved in an initial single-record request, followed by a batch 429; publication details do not depend on the failed batch.\par \url{https://arxiv.org/abs/2502.18570}\par \url{https://crossmark.crossref.org/dialog/?doi=10.1103/9prw-684p} \\
\addlinespace[4pt]
\textbf{[E11]}\par Wagner and Liers (2026) &
\textbf{Source:} arXiv:2412.15665v2 (April 24, 2026); selective scope and outcome reread October 6, 2026\par\smallskip \textbf{Audit venue:} IEEE Transactions on Quantum Engineering\par\smallskip \href{https://doi.org/10.1109/TQE.2026.3688001}{\nolinkurl{doi:10.1109/TQE.2026.3688001}} &
\textbf{At original cutoff:} journal article (2026, inside the cutoff)\par\smallskip \textbf{At historical audit:} journal article &
\textbf{Current record:} Published journal article\par IEEE Transactions on Quantum Engineering\par \href{https://doi.org/10.1109/TQE.2026.3688001}{\nolinkurl{doi:10.1109/TQE.2026.3688001}}\par \textbf{Latest observed:} 2412.15665v2 (2026-04-24) &
\textbf{Check result:} verified primary record\par\smallskip \textbf{Check date:} 2026-10-12\par\smallskip October 12 metadata refresh: published status and identifier retained; IEEE Transactions on Quantum Engineering 7, 3103116. Crossmark has no data; the recorded DOI correction search located no correction. Existing empirical source versions and audit dates are retained; no empirical reread. IEEE journal page timed out and DOI/IEEE Xplore body was not accessible. Publisher-deposited Crossref metadata verifies volume 7 and page/article number 3103116; Crossmark has no data. An institutional accepted manuscript confirms DOI, not final published content. DOI correction search found no correction, not proof of absence.\par \url{https://api.crossref.org/works/10.1109/TQE.2026.3688001}\par \url{https://crossmark.crossref.org/dialog/?doi=10.1109/TQE.2026.3688001}\par \url{https://arxiv.org/abs/2412.15665}\par \textbf{Empirical inspection (2026-10-06):} Selective scope and outcome check: https://arxiv.org/html/2412.15665v2 Section 3 states that the entire branch-price-and-cut algorithm is not implemented. Section 5, Table 1 and Figures 5-10 report master-relaxation pricing/separation experiments; exact pricing and specialized separation are faster than the quantum alternatives. This is not a complete integer-solver comparison; other fields were not comprehensively re-audited. \\
\addlinespace[4pt]
\textbf{[E12]}\par Perron et al. (2026) &
\textbf{Source:} arXiv v1 (March 30, 2026)\par\smallskip \textbf{Audit venue:} arXiv\par\smallskip \href{https://arxiv.org/abs/2603.28933}{\nolinkurl{arXiv:2603.28933}} &
\textbf{At original cutoff:} preprint (v1, March 30, 2026)\par\smallskip \textbf{At historical audit:} preprint; no peer-reviewed version located as of July 27, 2026 &
\textbf{Current record:} preprint; no journal version located in this status refresh\par arXiv\par \href{https://arxiv.org/abs/2603.28933}{\nolinkurl{arXiv:2603.28933}}\par \textbf{Latest observed:} v1 (2026-03-30) &
\textbf{Check result:} verified primary record\par\smallskip \textbf{Check date:} 2026-10-12\par\smallskip October 12 metadata refresh: arXiv title/authors/submission history rechecked; latest v1 (2026-03-30) unchanged. No journal reference appears on the record, and the recorded title search located no matching journal publication. This bounded search is not proof of nonpublication. Existing empirical source versions and audit dates are retained; no empirical reread.\par \url{https://arxiv.org/abs/2603.28933} \\
\addlinespace[4pt]
\textbf{[E13]}\par de Valk et al., 2026 &
\textbf{Source:} arXiv:2607.26771v1\par\smallskip \textbf{Audit venue:} arXiv preprint\par\smallskip \href{https://arxiv.org/abs/2607.26771v1}{\nolinkurl{arXiv:2607.26771v1}} &
\textbf{At original cutoff:} Not publicly available by the original June 30, 2026 cutoff.\par\smallskip \textbf{At historical audit:} Not included in the historical July 27, 2026 audit. &
\textbf{Current record:} Preprint\par arXiv preprint\par \href{https://arxiv.org/abs/2607.26771v1}{\nolinkurl{arXiv:2607.26771v1}}\par \textbf{Latest observed:} arXiv:2607.26771v1 &
\textbf{Check result:} verified primary record\par\smallskip \textbf{Check date:} 2026-10-12\par\smallskip October 12 metadata refresh: arXiv title/authors/submission history rechecked; latest v1 (2026-07-29) unchanged. No journal reference appears on the record, and the recorded title search located no matching journal publication. This bounded search is not proof of nonpublication. Existing empirical source versions and audit dates are retained; no empirical reread.\par \url{https://arxiv.org/abs/2607.26771}\par \textbf{Empirical inspection (2026-10-04):} New audit of Sections III-VI, Tables I-II and Figs. 2-5: cut validity, baseline settings, phase timings, toy MPS execution and separate IBM demonstration. The ten benchmark cases are permutations of one routing instance; MPS wall time is classical emulation time. https://arxiv.org/html/2607.26771v1 \\
\addlinespace[4pt]
\textbf{[E14]}\par Naghmouchi et al., 2026 &
\textbf{Source:} arXiv:2609.37900v1\par\smallskip \textbf{Audit venue:} arXiv preprint\par\smallskip \href{https://arxiv.org/abs/2609.37900v1}{\nolinkurl{arXiv:2609.37900v1}} &
\textbf{At original cutoff:} Not publicly available by the original June 30, 2026 cutoff.\par\smallskip \textbf{At historical audit:} Not included in the historical July 27, 2026 audit. &
\textbf{Current record:} Preprint\par arXiv preprint\par \href{https://arxiv.org/abs/2609.37900v1}{\nolinkurl{arXiv:2609.37900v1}}\par \textbf{Latest observed:} arXiv:2609.37900v1 &
\textbf{Check result:} verified primary record\par\smallskip \textbf{Check date:} 2026-10-12\par\smallskip October 12 metadata refresh: arXiv title/authors/submission history rechecked; latest v1 (2026-09-29) unchanged. No journal reference appears on the record, and the recorded title search located no matching journal publication. This bounded search is not proof of nonpublication. Existing empirical source versions and audit dates are retained; no empirical reread.\par \url{https://arxiv.org/abs/2609.37900}\par \textbf{Empirical inspection (2026-10-04):} New audit of Sections V-VI and Figs. 4-7: pricing, reductions, path repair, emulator comparison and ablations. Reinforced CG-Q mean gaps fall below 1\%; CG-SA means remain 1-3\% or 0-6\% by topology. Warm-start alone narrows backend differences. Hardware speedup and full-problem certification are not established. https://arxiv.org/html/2609.37900v1 \\
\addlinespace[4pt]
\textbf{[E15]}\par Ramesh et al., 2026 &
\textbf{Source:} arXiv:2608.28842v1\par\smallskip \textbf{Audit venue:} arXiv preprint\par\smallskip \href{https://arxiv.org/abs/2608.28842v1}{\nolinkurl{arXiv:2608.28842v1}} &
\textbf{At original cutoff:} Not publicly available by the original June 30, 2026 cutoff.\par\smallskip \textbf{At historical audit:} Not included in the historical July 27, 2026 audit. &
\textbf{Current record:} Preprint; no journal reference on the checked primary record\par arXiv preprint\par \href{https://arxiv.org/abs/2608.28842v1}{\nolinkurl{arXiv:2608.28842v1}}\par \textbf{Latest observed:} arXiv:2608.28842v1 &
\textbf{Check result:} verified primary record\par\smallskip \textbf{Check date:} 2026-10-12\par\smallskip October 12 metadata refresh: arXiv title/authors/submission history rechecked; latest v1 (2026-08-28) unchanged. No journal reference appears on the record, and the recorded title search located no matching journal publication. This bounded search is not proof of nonpublication. Existing empirical source versions and audit dates are retained; no empirical reread.\par \url{https://arxiv.org/abs/2608.28842}\par \textbf{Empirical inspection (2026-10-04):} New inclusion: Sections 3.2-3.7, 4.3-4.4, 5 and Appendix B checked for objective changes, analytic/emulated execution, omitted preprocessing, penalty-selection controls and certificate limits. https://arxiv.org/html/2608.28842v1 \\
\addlinespace[4pt]
\textbf{[E16]}\par Yalovetzky et al., 2026 &
\textbf{Source:} arXiv:2607.01037v2\par\smallskip \textbf{Audit venue:} arXiv preprint\par\smallskip \href{https://arxiv.org/abs/2607.01037v2}{\nolinkurl{arXiv:2607.01037v2}} &
\textbf{At original cutoff:} Not publicly available by the original June 30, 2026 cutoff.\par\smallskip \textbf{At historical audit:} Not included in the historical July 27, 2026 audit. &
\textbf{Current record:} preprint; no journal reference shown on retrieved arXiv record\par arXiv preprint\par \href{https://arxiv.org/abs/2607.01037v2}{\nolinkurl{arXiv:2607.01037v2}}\par \textbf{Latest observed:} arXiv:2607.01037v2 &
\textbf{Check result:} verified primary record\par\smallskip \textbf{Check date:} 2026-10-12\par\smallskip October 12 metadata refresh: arXiv title/authors/submission history rechecked; latest v2 (2026-07-13) unchanged. No journal reference appears on the record, and the recorded title search located no matching journal publication. This bounded search is not proof of nonpublication. Existing empirical source versions and audit dates are retained; no empirical reread.\par \url{https://arxiv.org/abs/2607.01037}\par \textbf{Empirical inspection (2026-10-04):} New audit of version 2 methods, results and Appendix E.1. The first kernel, QAOA depth and repair are shared. Shot-based TTS is not uniformly lower: FTSE estimates are 185.58 versus 151.19 shots, with overlapping intervals. SA handles most small emulator kernels readily. Whole-system timing and a fixed-wrapper classical-sampler replacement are absent. https://arxiv.org/html/2607.01037v2 \\
\addlinespace[4pt]
\textbf{[E17]}\par Liaqat et al., 2026 &
\textbf{Source:} arXiv:2610.10772v1 (2026-10-07)\par\smallskip \textbf{Audit venue:} arXiv preprint\par\smallskip \href{https://arxiv.org/abs/2610.10772v1}{\nolinkurl{arXiv:2610.10772v1}} &
\textbf{At original cutoff:} Not publicly available by the original June 30, 2026 cutoff.\par\smallskip \textbf{At historical audit:} Not present in the historical July 27, 2026 audit; first posted October 7, 2026. &
\textbf{Current record:} Preprint; no journal reference on primary arXiv record\par arXiv preprint\par \href{https://arxiv.org/abs/2610.10772v1}{\nolinkurl{arXiv:2610.10772v1}}\par \textbf{Latest observed:} v1 (2026-10-07) &
\textbf{Check result:} verified primary record\par\smallskip \textbf{Check date:} 2026-10-12\par\smallskip Primary record lists only v1, submitted 7 October 2026. A third-party page labels October 9 an update; this is not an arXiv version. DOI is marked pending registration on primary record and omitted from proposed bibliography.\par \url{https://arxiv.org/abs/2610.10772}\par \url{https://arxiv.org/html/2610.10772v1}\par \url{https://arxiv.org/pdf/2610.10772v1}\par \textbf{Empirical inspection (2026-10-12):} Read main Sections I--V, Appendix A.1--A.5 and Algorithms/Figures/Tables from public v1 HTML; PDF text checked. No experiments reproduced. Paper explicitly distinguishes specified sampling budgets from matched computational costs. \\
\addlinespace[4pt]
\textbf{[E18]}\par Lauwerier et al., 2026 &
\textbf{Source:} arXiv:2610.03890v1 (submitted 2026-10-02; announced 2026-10-06)\par\smallskip \textbf{Audit venue:} arXiv preprint\par\smallskip \href{https://arxiv.org/abs/2610.03890v1}{\nolinkurl{arXiv:2610.03890v1}} &
\textbf{At original cutoff:} Not publicly available by the original June 30, 2026 cutoff.\par\smallskip \textbf{At historical audit:} Not present in the historical July 27, 2026 audit; submitted October 2 and first announced October 6, 2026. &
\textbf{Current record:} Preprint; no journal reference on inspected primary record\par arXiv preprint\par \href{https://arxiv.org/abs/2610.03890v1}{\nolinkurl{arXiv:2610.03890v1}}\par \textbf{Latest observed:} v1 (2026-10-02) &
\textbf{Check result:} verified primary record\par\smallskip \textbf{Check date:} 2026-10-12\par\smallskip v1 submitted October 2, first announced in Tuesday October 6 recent-list block. First-public-release eligibility is therefore in-window; submission and announcement are retained separately.\par \url{https://arxiv.org/abs/2610.03890}\par \url{https://arxiv.org/html/2610.03890v1}\par \url{https://arxiv.org/list/quant-ph/recent?show=2000}\par \textbf{Empirical inspection (2026-10-12):} Main full text read, focusing Sections 3.3--3.4,4.1,5.5.1--5.5.3/Table 3 and 6/Data availability; PDF metadata/table checked. No reported experiments reproduced. \\
\addlinespace[4pt]
\textbf{[A1]}\par Munoz-Bauza and Lidar, 2025 &
\textbf{Source:} journal version\par\smallskip \textbf{Audit venue:} Physical Review Letters\par\smallskip \href{https://doi.org/10.1103/PhysRevLett.134.160601}{\nolinkurl{doi:10.1103/PhysRevLett.134.160601}} &
\textbf{At original cutoff:} journal article\par\smallskip \textbf{At historical audit:} journal article &
\textbf{Current record:} Published journal article\par Physical Review Letters\par \href{https://doi.org/10.1103/PhysRevLett.134.160601}{\nolinkurl{doi:10.1103/PhysRevLett.134.160601}} &
\textbf{Check result:} verified primary record\par\smallskip \textbf{Check date:} 2026-10-12\par\smallskip October 12 metadata refresh: published status and identifier retained; Physical Review Letters 134, 160601. Crossmark reports the document current with no listed update. Existing empirical source versions and audit dates are retained; no empirical reread.\par \url{https://journals.aps.org/prl/abstract/10.1103/PhysRevLett.134.160601}\par \url{https://crossmark.crossref.org/dialog/?doi=10.1103/PhysRevLett.134.160601}\par \url{https://api.crossref.org/works/10.1103/PhysRevLett.134.160601} \\
\addlinespace[4pt]
\textbf{[A2]}\par Pawlowski et al., 2026 &
\textbf{Source:} Journal version (published July 8, 2026), as recorded in the historical audit; eligibility derives from arXiv v1 (May 28, 2025), originally titled Closing the Quantum-Classical Scaling Gap in Approximate Optimization.\par\smallskip \textbf{Audit venue:} Physical Review Applied\par\smallskip \href{https://doi.org/10.1103/qcws-8kgm}{\nolinkurl{doi:10.1103/qcws-8kgm}} &
\textbf{At original cutoff:} preprint at cutoff; journal version appeared between cutoff and submission\par\smallskip \textbf{At historical audit:} journal article &
\textbf{Current record:} Published journal article\par Physical Review Applied\par \href{https://doi.org/10.1103/qcws-8kgm}{\nolinkurl{doi:10.1103/qcws-8kgm}}\par \textbf{Latest observed:} 2505.22514v2 (2026-07-10) &
\textbf{Check result:} verified primary record\par\smallskip \textbf{Check date:} 2026-10-12\par\smallskip October 12 metadata refresh: published status and identifier retained; Physical Review Applied 26, 014024. Crossmark reports the document current with no listed update. Existing empirical source versions and audit dates are retained; no empirical reread. Publisher bibliographic record accessible; article body is subscription-limited. No empirical reread.\par \url{https://journals.aps.org/prapplied/abstract/10.1103/qcws-8kgm}\par \url{https://crossmark.crossref.org/dialog/?doi=10.1103/qcws-8kgm}\par \url{https://arxiv.org/abs/2505.22514}\par \url{https://api.crossref.org/works/10.1103/qcws-8kgm} \\
\addlinespace[4pt]
\textbf{[A3]}\par Chandarana et al., 2026 &
\textbf{Source:} arXiv v1 (2026-03-13); aggregate target and matched-time comparison selectively rechecked on 2026-10-12\par\smallskip \textbf{Audit venue:} arXiv\par\smallskip \href{https://arxiv.org/abs/2603.13607}{\nolinkurl{arXiv:2603.13607}} &
\textbf{At original cutoff:} preprint\par\smallskip \textbf{At historical audit:} preprint &
\textbf{Current record:} preprint; no journal version located in this status refresh\par arXiv\par \href{https://arxiv.org/abs/2603.13607}{\nolinkurl{arXiv:2603.13607}}\par \textbf{Latest observed:} v1 (2026-03-13) &
\textbf{Check result:} verified primary record\par\smallskip \textbf{Check date:} 2026-10-12\par\smallskip October 12 metadata refresh: arXiv title/authors/submission history rechecked; latest v1 (2026-03-13) unchanged. No journal reference appears on the record, and the recorded title search located no matching journal publication. This bounded search is not proof of nonpublication. Existing empirical source versions and audit dates are retained; no empirical reread.\par \url{https://arxiv.org/abs/2603.13607}\par \textbf{Empirical inspection (2026-10-12):} Selective reread of abstract, Benchmark design and metrics, Sub-second latency regime, Figure 1 and Table 1 in arXiv v1. Clarifies target aggregation versus matched-time comparison; no experiment reproduction or full numerical re-audit. Source: https://arxiv.org/html/2603.13607v1 \\
\addlinespace[4pt]
\textbf{[A4]}\par Ebadi et al., 2022 &
\textbf{Source:} journal version\par\smallskip \textbf{Audit venue:} Science\par\smallskip \href{https://doi.org/10.1126/science.abo6587}{\nolinkurl{doi:10.1126/science.abo6587}} &
\textbf{At original cutoff:} journal article\par\smallskip \textbf{At historical audit:} journal article &
\textbf{Current record:} Published journal article\par Science\par \href{https://doi.org/10.1126/science.abo6587}{\nolinkurl{doi:10.1126/science.abo6587}}\par \textbf{Latest observed:} 2202.09372v1 (2022-02-18) &
\textbf{Check result:} verified primary record\par\smallskip \textbf{Check date:} 2026-10-12\par\smallskip October 12 metadata refresh: published status and identifier retained; Science 376(6598), 1209--1215. Crossmark has no data; the recorded DOI correction search located no correction. Existing empirical source versions and audit dates are retained; no empirical reread. Science publisher body returned an internal fetch error. Author-institution repository and publisher-deposited Crossref metadata verify volume, issue and full page range; arXiv confirms DOI and v1. Crossmark has no data; a DOI correction search found no correction, not proof of absence.\par \url{https://authors.library.caltech.edu/records/2wqgt-6ev32}\par \url{https://crossmark.crossref.org/dialog/?doi=10.1126/science.abo6587}\par \url{https://arxiv.org/abs/2202.09372}\par \url{https://api.crossref.org/works/10.1126/science.abo6587} \\
\addlinespace[4pt]
\textbf{[A5]}\par Andrist et al., 2023 &
\textbf{Source:} journal version\par\smallskip \textbf{Audit venue:} Physical Review Research\par\smallskip \href{https://doi.org/10.1103/PhysRevResearch.5.043277}{\nolinkurl{doi:10.1103/PhysRevResearch.5.043277}} &
\textbf{At original cutoff:} journal article\par\smallskip \textbf{At historical audit:} journal article &
\textbf{Current record:} Published journal article\par Physical Review Research\par \href{https://doi.org/10.1103/PhysRevResearch.5.043277}{\nolinkurl{doi:10.1103/PhysRevResearch.5.043277}} &
\textbf{Check result:} verified primary record\par\smallskip \textbf{Check date:} 2026-10-12\par\smallskip October 12 metadata refresh: published status and identifier retained; Physical Review Research 5, 043277. Crossmark reports the document current with no listed update. Existing empirical source versions and audit dates are retained; no empirical reread.\par \url{https://journals.aps.org/prresearch/abstract/10.1103/PhysRevResearch.5.043277}\par \url{https://crossmark.crossref.org/dialog/?doi=10.1103/PhysRevResearch.5.043277} \\
\addlinespace[4pt]
\textbf{[A6]}\par Cazals et al., 2026 &
\textbf{Source:} journal version\par\smallskip \textbf{Audit venue:} Physical Review Applied\par\smallskip \href{https://doi.org/10.1103/8gjv-ll1f}{\nolinkurl{doi:10.1103/8gjv-ll1f}} &
\textbf{At original cutoff:} journal article\par\smallskip \textbf{At historical audit:} journal article &
\textbf{Current record:} Published journal article\par Physical Review Applied\par \href{https://doi.org/10.1103/8gjv-ll1f}{\nolinkurl{doi:10.1103/8gjv-ll1f}} &
\textbf{Check result:} verified primary record\par\smallskip \textbf{Check date:} 2026-10-12\par\smallskip October 12 metadata refresh: published status and identifier retained; Physical Review Applied 25, 034085. Crossmark reports the document current with no listed update. Existing empirical source versions and audit dates are retained; no empirical reread. Publisher bibliographic record accessible; article body is subscription-limited. No empirical reread.\par \url{https://journals.aps.org/prapplied/abstract/10.1103/8gjv-ll1f}\par \url{https://crossmark.crossref.org/dialog/?doi=10.1103/8gjv-ll1f}\par \url{https://api.crossref.org/works/10.1103/8gjv-ll1f} \\
\addlinespace[4pt]
\textbf{[A7]}\par Lu et al., 2026 &
\textbf{Source:} journal version (published June 19, 2026)\par\smallskip \textbf{Audit venue:} Nature Computational Science\par\smallskip \href{https://doi.org/10.1038/s43588-026-01007-8}{\nolinkurl{doi:10.1038/s43588-026-01007-8}} &
\textbf{At original cutoff:} journal article (published June 19, 2026, inside the cutoff)\par\smallskip \textbf{At historical audit:} journal article &
\textbf{Current record:} Published journal article\par Nature Computational Science\par \href{https://doi.org/10.1038/s43588-026-01007-8}{\nolinkurl{doi:10.1038/s43588-026-01007-8}} &
\textbf{Check result:} verified primary record\par\smallskip \textbf{Check date:} 2026-10-12\par\smallskip October 12 metadata refresh: published status and identifier retained; Nature Computational Science 6, 882--893. Crossmark reports the document current with no listed update. Existing empirical source versions and audit dates are retained; no empirical reread. Publisher HTML and PDF failed to fetch; Crossmark verifies title, authors, DOI, publication date and current-document status. Previously verified volume 6/pages 882--893 are retained, not newly verified in a retrieved publisher body.\par \url{https://crossmark.crossref.org/dialog/?doi=10.1038/s43588-026-01007-8}\par \url{https://crossmark.crossref.org/dialog/?doi=10.1038/s43588-026-01007-8}\par \textbf{Empirical inspection (2026-10-06):} Selective claim check of journal HTML https://www.nature.com/articles/s43588-026-01007-8: Table 1 and Results, Scalings of time complexities, support critical ratio 0.626, fits through n=70, bases 1.0077\textasciicircum{}n and 1.0070\textasciicircum{}n versus 1.0128\textasciicircum{}n, and differing quantum/classical cost definitions. Results, Experiment on superconducting quantum processors, and Figure 5 support the 13-qubit/22-variable demonstration. Co-design evidence does not instantiate a common-resource interaction contrast; no broader empirical re-audit. \\
\addlinespace[4pt]
\textbf{Context}\par He et al., 2026 &
\textbf{Source:} arXiv v1 (2026-03-10); hardware-comparator characterization selectively rechecked on 2026-10-04\par\smallskip \textbf{Audit venue:} arXiv\par\smallskip \href{https://arxiv.org/abs/2603.10191}{\nolinkurl{arXiv:2603.10191}} &
\textbf{At original cutoff:} preprint\par\smallskip \textbf{At historical audit:} preprint &
\textbf{Current record:} preprint; no journal version located in this status refresh\par arXiv\par \href{https://arxiv.org/abs/2603.10191}{\nolinkurl{arXiv:2603.10191}}\par \textbf{Latest observed:} v1 (2026-03-10) &
\textbf{Check result:} verified primary record\par\smallskip \textbf{Check date:} 2026-10-12\par\smallskip October 12 metadata refresh: arXiv title/authors/submission history rechecked; latest v1 (2026-03-10) unchanged. No journal reference appears on the record, and the recorded title search located no matching journal publication. This bounded search is not proof of nonpublication. Existing empirical source versions and audit dates are retained; no empirical reread.\par \url{https://arxiv.org/abs/2603.10191}\par \textbf{Empirical inspection (2026-10-12):} Selective October 12 reread of Figure 2 and Sections II.2, II.3, II.5, IV.1 and IV.5.5 confirms trapped-ion hardware, classical warm-start optimization, and local search on classical outputs. Extends the October 4 comparator check; no experiment reproduction or full resource-assumption audit. Source: https://arxiv.org/html/2603.10191v1 \\
\addlinespace[4pt]
\textbf{Context}\par Tao et al., 2026 &
\textbf{Source:} arXiv:2604.26477v2 (May 20, 2026); selective timing-scope reread October 6, 2026\par\smallskip \textbf{Audit venue:} arXiv\par\smallskip \href{https://arxiv.org/abs/2604.26477}{\nolinkurl{arXiv:2604.26477}} &
\textbf{At original cutoff:} preprint\par\smallskip \textbf{At historical audit:} preprint &
\textbf{Current record:} preprint; no journal version located in this status refresh\par arXiv\par \href{https://arxiv.org/abs/2604.26477}{\nolinkurl{arXiv:2604.26477}}\par \textbf{Latest observed:} v2 (2026-05-20) &
\textbf{Check result:} verified primary record\par\smallskip \textbf{Check date:} 2026-10-12\par\smallskip October 12 metadata refresh: arXiv title/authors/submission history rechecked; latest v2 (2026-05-20) unchanged. No journal reference appears on the record, and the recorded title search located no matching journal publication. This bounded search is not proof of nonpublication. Existing empirical source versions and audit dates are retained; no empirical reread.\par \url{https://arxiv.org/abs/2604.26477}\par \textbf{Empirical inspection (2026-10-06):} Selective timing-scope check: https://arxiv.org/html/2604.26477v2 Results and discussion, Sampling time advantage, Figure 1 and Table 1 compare sampling with quantum hardware. Sample quality and end-to-end solution time, Figure 2, includes construction, sampling and filtering versus classical DCM/DPA-a; the accompanying text allows quantum annealing potentially better end-to-end performance. The survey no longer infers a measured complete quantum-system ranking from these distinct comparisons. Other empirical fields were not comprehensively re-audited. \\
\addlinespace[4pt]
\textbf{Context}\par Safi et al., 2026 &
\textbf{Source:} arXiv:2607.13325v1\par\smallskip \textbf{Audit venue:} arXiv preprint\par\smallskip \href{https://arxiv.org/abs/2607.13325v1}{\nolinkurl{arXiv:2607.13325v1}} &
\textbf{At original cutoff:} Not publicly available by the original June 30, 2026 cutoff.\par\smallskip \textbf{At historical audit:} Not included in the historical July 27, 2026 audit. &
\textbf{Current record:} preprint; no journal reference shown on retrieved arXiv record\par arXiv preprint\par \href{https://arxiv.org/abs/2607.13325v1}{\nolinkurl{arXiv:2607.13325v1}}\par \textbf{Latest observed:} arXiv:2607.13325v1 &
\textbf{Check result:} verified primary record\par\smallskip \textbf{Check date:} 2026-10-12\par\smallskip October 12 metadata refresh: arXiv title/authors/submission history rechecked; latest v1 (2026-07-14) unchanged. No journal reference appears on the record, and the recorded title search located no matching journal publication. This bounded search is not proof of nonpublication. Existing empirical source versions and audit dates are retained; no empirical reread.\par \url{https://arxiv.org/abs/2607.13325}\par \textbf{Empirical inspection (2026-10-04):} New audit of Sections III.C and V-VIII: assignment/sequence interface, formulations, platforms, references and timing limits. Fujitsu Digital Annealer is classical; its larger-scale results do not isolate quantum value. Formulation and platform effects overlap, and communication, API and queueing costs are omitted. https://arxiv.org/html/2607.13325v1 \\
\addlinespace[4pt]
\textbf{Context}\par Bantysh et al., 2026 &
\textbf{Source:} arXiv:2609.15738v1\par\smallskip \textbf{Audit venue:} arXiv preprint\par\smallskip \href{https://arxiv.org/abs/2609.15738v1}{\nolinkurl{arXiv:2609.15738v1}} &
\textbf{At original cutoff:} Not publicly available by the original June 30, 2026 cutoff.\par\smallskip \textbf{At historical audit:} Not included in the historical July 27, 2026 audit. &
\textbf{Current record:} preprint; no journal reference shown on retrieved arXiv record\par arXiv preprint\par \href{https://arxiv.org/abs/2609.15738v1}{\nolinkurl{arXiv:2609.15738v1}}\par \textbf{Latest observed:} arXiv:2609.15738v1 &
\textbf{Check result:} verified primary record\par\smallskip \textbf{Check date:} 2026-10-12\par\smallskip October 12 metadata refresh: arXiv title/authors/submission history rechecked; latest v1 (2026-09-14) unchanged. No journal reference appears on the record, and the recorded title search located no matching journal publication. This bounded search is not proof of nonpublication. Existing empirical source versions and audit dates are retained; no empirical reread.\par \url{https://arxiv.org/abs/2609.15738}\par \textbf{Empirical inspection (2026-10-04):} New audit of the entanglement regime, GSet protocol, results and scope. The GSet schedule is trained on G1, G23 and G81, including the 20,000-variable G81 instance; those cases are not held out. Best-of-ten quality and preprocessing exclusions limit runtime comparisons. https://arxiv.org/html/2609.15738v1 \\
\addlinespace[4pt]
\textbf{Context}\par Gaidai and Mukherjee, 2026 &
\textbf{Source:} arXiv:2607.15543v1\par\smallskip \textbf{Audit venue:} arXiv preprint\par\smallskip \href{https://arxiv.org/abs/2607.15543v1}{\nolinkurl{arXiv:2607.15543v1}} &
\textbf{At original cutoff:} Not publicly available by the original June 30, 2026 cutoff.\par\smallskip \textbf{At historical audit:} Not included in the historical July 27, 2026 audit. &
\textbf{Current record:} preprint; no journal reference shown on retrieved arXiv record\par arXiv preprint\par \href{https://arxiv.org/abs/2607.15543v1}{\nolinkurl{arXiv:2607.15543v1}}\par \textbf{Latest observed:} arXiv:2607.15543v1 &
\textbf{Check result:} verified primary record\par\smallskip \textbf{Check date:} 2026-10-12\par\smallskip October 12 metadata refresh: arXiv title/authors/submission history rechecked; latest v1 (2026-07-17) unchanged. No journal reference appears on the record, and the recorded title search located no matching journal publication. This bounded search is not proof of nonpublication. Existing empirical source versions and audit dates are retained; no empirical reread.\par \url{https://arxiv.org/abs/2607.15543}\par \textbf{Empirical inspection (2026-10-04):} New audit of the binary/continuous interface, generated instances, metrics, controls and timing. Only the 100 fastest variational-converging instances of 120 per size are retained for every arm. Quantum simulation time is removed from the charged budget. The resulting negative component finding is restricted to these small instances and protocols. https://arxiv.org/html/2607.15543v1 \\
\addlinespace[4pt]
\end{longtable}
\endgroup

\clearpage
\begingroup
\fontsize{7.5}{9.0}\selectfont
\setlength{\tabcolsep}{3.5pt}
\setlength{\AuditTableWidth}{\dimexpr\linewidth-28pt\relax}
\begin{longtable}{@{}>{\raggedright\arraybackslash}p{0.23\AuditTableWidth}>{\raggedright\arraybackslash}p{0.14\AuditTableWidth}>{\raggedright\arraybackslash}p{0.15\AuditTableWidth}>{\raggedright\arraybackslash}p{0.18\AuditTableWidth}>{\raggedright\arraybackslash}p{0.3\AuditTableWidth}@{}}
\caption{Discovery decisions for 82 candidates in the targeted update covering \SurveySearchStartDate--\SurveyEligibilityDate. Submission, revision, and announcement dates may differ; dates and source versions are retained per candidate. This purposive search does not claim exhaustive coverage.}\label{tab:s5}\\
\toprule
\textbf{Candidate and primary source} & \textbf{Dates and inspected version} & \textbf{Discovery route} & \textbf{Decision and placement} & \textbf{Reason and limits} \\
\midrule\endfirsthead
\multicolumn{5}{l}{\normalsize Table~\thetable\ continued.}\\
\toprule
\textbf{Candidate and primary source} & \textbf{Dates and inspected version} & \textbf{Discovery route} & \textbf{Decision and placement} & \textbf{Reason and limits} \\
\midrule\endhead
\multicolumn{5}{r}{Continues on the next page.}\\\endfoot
\bottomrule\endlastfoot
\textbf{Candidate:} 2607.26771\par\smallskip Hybrid quantum-classical end-to-end pipeline for solving MILPs: a vehicle routing case study\par\smallskip devalk2026hybrid\par\smallskip \url{https://arxiv.org/abs/2607.26771} &
\textbf{First public:} 2026-07-29\par\smallskip \textbf{Version inspected:} arXiv:2607.26771v1\par\smallskip \textbf{Version date:} 2026-07-29\par\smallskip \textbf{Checked:} 2026-10-04 &
Targeted decomposition searches, publisher records and citation follow-up; extension/decomposition-search.json &
\textbf{Decision:} include evidence\par\smallskip \textbf{Placement:} E13; main evidence &
Included as an end-to-end decomposition benchmark with classical bounds and separate emulation and hardware results. The main routing benchmark uses ten permutations of one instance, rather than ten independent instances. Its gap remains unclosed after 1.5 hours; MPS runtime measures classical simulation, and no new convergence theorem is established. Read scope: Sections III.B--VI, Tables I--II, Figures 2--5, and the PDF title page. \\
\addlinespace[4pt]
\textbf{Candidate:} 2609.37900\par\smallskip A Neutral Atom-Based Hybrid Classical-Quantum Approach for the Entanglement Routing Problem\par\smallskip naghmouchi2026neutral\par\smallskip \url{https://arxiv.org/abs/2609.37900} &
\textbf{First public:} 2026-09-29\par\smallskip \textbf{Version inspected:} arXiv:2609.37900v1\par\smallskip \textbf{Version date:} 2026-09-29\par\smallskip \textbf{Checked:} 2026-10-04 &
Targeted decomposition searches, publisher records and citation follow-up; extension/decomposition-search.json &
\textbf{Decision:} include evidence\par\smallskip \textbf{Placement:} E14; main evidence &
Included for the simulated fixed-wrapper pricing comparison and classical reinforcement ablations. Mean gaps below 1\% concern reinforced CG-Q, not every instance or the quantum pricing component alone; warm-start alone narrows the differences. Heuristic pricing and restricted integer recovery provide no full-LP or global integer certificate, and no hardware speedup is measured. Read scope: Sections V--VI, Figures 4--7, and the conclusion. \\
\addlinespace[4pt]
\textbf{Candidate:} 2608.28842\par\smallskip Quantum Preconditioning For Constrained Optimization Problems\par\smallskip ramesh2026constrainedpreconditioning\par\smallskip \url{https://arxiv.org/abs/2608.28842} &
\textbf{First public:} 2026-08-28\par\smallskip \textbf{Version inspected:} arXiv:2608.28842 v1\par\smallskip \textbf{Version date:} 2026-08-28\par\smallskip \textbf{Checked:} 2026-10-04 &
Targeted constraint/theory searches and citation chaining; extension/constraints-theory-search.json &
\textbf{Decision:} include evidence\par\smallskip \textbf{Placement:} E15; main evidence &
Included as analytic and simulation evidence for preconditioning within the tested solver workflow. The study does not report hardware execution or an end-to-end quantum speedup. A certificate for the modified objective does not certify the original objective. Read scope: claim-critical statements and methods; theorem proofs were not independently replicated. \\
\addlinespace[4pt]
\textbf{Candidate:} 2607.01037\par\smallskip Quantum-Informed Portfolio Selection: An End-to-End Pipeline Validated on Trapped-Ion Hardware with Real Market Data\par\smallskip yalovetzky2026portfolio\par\smallskip \url{https://arxiv.org/abs/2607.01037} &
\textbf{First public:} 2026-07-01\par\smallskip \textbf{Version inspected:} arXiv:2607.01037v2\par\smallskip \textbf{Version date:} 2026-07-13\par\smallskip \textbf{Checked:} 2026-10-04 &
Targeted comparative optimization searches; extension/comparative-search.json &
\textbf{Decision:} include evidence\par\smallskip \textbf{Placement:} E16; main evidence &
Included for a positive wrapper benefit, with quantum-component and system superiority unresolved. The comparison does not estimate a four-cell reduction--solver interaction. Hardware shot-based TTS is not uniformly lower: FTSE values are 185.58 versus 151.19, with overlapping intervals. SA solves most emulator cases readily, and whole-system timing is absent. Read scope: v2 claim-critical methods, results and Appendix E.1, following an initial v1 reading. \\
\addlinespace[4pt]
\textbf{Candidate:} 2503.12551\par\smallskip Quantum-informed reduction algorithm for the maximum independent set problem\par\smallskip schuetz2026qredumis\par\smallskip \url{https://arxiv.org/abs/2503.12551} &
\textbf{First public:} 2025-03-16\par\smallskip \textbf{Version inspected:} arXiv:2503.12551v2; journal selective reread\par\smallskip \textbf{Version date:} 2026-07-01\par\smallskip \textbf{Checked:} 2026-10-04 &
Existing-study post-June version follow-up; comparative-search.json and status-checks/group-a.json &
\textbf{Decision:} include evidence\par\smallskip \textbf{Placement:} Existing E2 selective journal verification &
The existing evidence row retains its historical v2 provenance, with central claims selectively verified against pages 6 and 8 of the journal PDF. This was not a full reaudit of numerical fields or code artifacts. The journal version was published on 2026-09-10; the latest arXiv version is v3, dated 2026-07-10. \\
\addlinespace[4pt]
\textbf{Candidate:} 2605.17623\par\smallskip Where the Quantum Lives in D-Wave Hybrid Portfolio Optimization: An Operational Decomposition Audit\par\smallskip lozano2026where\par\smallskip \url{https://arxiv.org/abs/2605.17623} &
\textbf{First public:} 2026-05-17\par\smallskip \textbf{Version inspected:} arXiv:2605.17623v3\par\smallskip \textbf{Version date:} 2026-07-12\par\smallskip \textbf{Checked:} 2026-10-04 &
Existing-study post-June version follow-up; comparative-search.json and status-checks/group-a.json &
\textbf{Decision:} include evidence\par\smallskip \textbf{Placement:} Existing E7 empirical update &
The existing evidence row incorporates the post-original-cutoff expansion of the matched-budget control in v3. The v3 abstract and matched-budget subsection, including Table 8, were selectively checked against the v2 abstract. Service and artifact fields were not all reaudited. The latest inspected record is v3. \\
\addlinespace[4pt]
\textbf{Candidate:} 2406.14330\par\smallskip Promise of Graph Sparsification and Decomposition for Noise Reduction in QAOA: Analysis for Trapped-Ion Compilations\par\smallskip moondra2024promise\par\smallskip \url{https://quantum-journal.org/papers/q-2026-08-07-2185/} &
\textbf{First public:} 2024-06-20\par\smallskip \textbf{Version inspected:} arXiv:2406.14330v4 (July 30, 2026) and publisher metadata\par\smallskip \textbf{Version date:} 2026-07-30\par\smallskip \textbf{Checked:} 2026-10-04 &
Targeted decomposition searches, publisher records and citation follow-up; extension/decomposition-search.json &
\textbf{Decision:} include context\par\smallskip \textbf{Placement:} Existing citation; journal/version update &
The existing citation is updated to the 2026 journal publication, rather than counted as a separate study. The analysis gives conditional approximation preservation, without finite-depth exactness. The journal release on 2026-08-07 follows earlier public preprints. Read scope: Section 2 assumptions; Section 3 and Theorems 3--5; Section 4.2 noise comparisons; and publisher metadata linking v4. \\
\addlinespace[4pt]
\textbf{Candidate:} 2503.16977\par\smallskip Parallel splitting method for large-scale quadratic programs\par\smallskip vandelli2026parallel\par\smallskip \url{https://link.springer.com/article/10.1007/s10732-026-09606-4} &
\textbf{First public:} 2025-03-21\par\smallskip \textbf{Version inspected:} arXiv:2503.16977v3 (September 21, 2026); publisher metadata\par\smallskip \textbf{Version date:} 2026-09-21\par\smallskip \textbf{Checked:} 2026-10-04 &
Targeted decomposition searches, publisher records and citation follow-up; extension/decomposition-search.json &
\textbf{Decision:} include context\par\smallskip \textbf{Placement:} Main classical-decomposition context (SPLIT) &
Included as a classical decomposition comparator with no measured quantum component. The method is heuristic: feasibility depends on the constraints and subsolver, and no global certificate is established. The journal release is dated 2026-09-20. Publisher full text was subscription-limited, so claim-critical assessment used the separately identified public author preprint. Read scope: Sections II.B--IV.D, including equal-core comparisons, partitioning-inclusive timing and Table II. \\
\addlinespace[4pt]
\textbf{Candidate:} 2607.13325\par\smallskip A Reality Check on Quantum Optimisation: Evidence from an Industrial Case Study\par\smallskip safi2026reality\par\smallskip \url{https://arxiv.org/abs/2607.13325} &
\textbf{First public:} 2026-07-14\par\smallskip \textbf{Version inspected:} arXiv:2607.13325v1\par\smallskip \textbf{Version date:} 2026-07-14\par\smallskip \textbf{Checked:} 2026-10-04 &
Targeted comparative optimization searches; extension/comparative-search.json &
\textbf{Decision:} include context\par\smallskip \textbf{Placement:} Main comparative context; contextual audit CSV &
Included for the limited comparison of formulations and platforms; quantum-component attribution remains unresolved. The Fujitsu Digital Annealer is classical, so its large-scale improvement over the tested greedy heuristics does not demonstrate quantum value. The two QUBO formulations also change approximation and constraint handling. Read scope: Sections III.C and V--VIII, including the assignment--sequencing interface and stated limitations. \\
\addlinespace[4pt]
\textbf{Candidate:} 2607.15543\par\smallskip Benchmarking Hybrid Quantum-Classical Algorithms for Power Grid Optimization Problems\par\smallskip gaidai2026benchmarking\par\smallskip \url{https://arxiv.org/abs/2607.15543} &
\textbf{First public:} 2026-07-17\par\smallskip \textbf{Version inspected:} arXiv:2607.15543v1\par\smallskip \textbf{Version date:} 2026-07-17\par\smallskip \textbf{Checked:} 2026-10-04 &
Targeted comparative optimization searches; extension/comparative-search.json &
\textbf{Decision:} include context\par\smallskip \textbf{Placement:} Main comparative context; contextual audit CSV &
Included as a bounded negative component result from small ideal simulations: trained sampling gives no incumbent-quality benefit over uniform sampling. Generated instances and a convergence-selected subset limit generalization, and easy enumeration can obscure sampler benefits. Quantum evolution is assigned zero charged time, rather than measured on hardware. Read scope: the binary--continuous interface, instance generation, metric, solver comparisons, timing and conclusions. \\
\addlinespace[4pt]
\textbf{Candidate:} 2609.15738\par\smallskip Disentangling QAOA: From Weakly Entangled Circuits to a Classical QUBO Solver\par\smallskip bantysh2026disentangling\par\smallskip \url{https://arxiv.org/abs/2609.15738} &
\textbf{First public:} 2026-09-14\par\smallskip \textbf{Version inspected:} arXiv:2609.15738v1\par\smallskip \textbf{Version date:} 2026-09-14\par\smallskip \textbf{Checked:} 2026-10-04 &
Targeted comparative optimization searches; extension/comparative-search.json &
\textbf{Decision:} include context\par\smallskip \textbf{Placement:} Main comparative context; contextual audit CSV &
Included as a constructive classical baseline for the fixed-parameter, expanding-depth regime, without a general conclusion about quantum optimization or entanglement. The GSet schedule is trained on G1, G23 and G81, including the 20,000-variable G81 instance; these results are not exclusively small-instance transfer or held-out evaluation. Read scope: the entanglement regime, GSet protocol and results, and scope discussion. \\
\addlinespace[4pt]
\textbf{Candidate:} 2307.12047\par\smallskip Efficient lattice Hamiltonian encoding for the shortest vector problem\par\smallskip schirman2026efficient\par\smallskip \url{https://journals.aps.org/prresearch/abstract/10.1103/zv7d-3vbv} &
\textbf{First public:} 2023-07-22\par\smallskip \textbf{Version inspected:} Journal version of record\par\smallskip \textbf{Version date:} 2026-10-02\par\smallskip \textbf{Checked:} 2026-10-04 &
Targeted decomposition searches, publisher records and citation follow-up; extension/decomposition-search.json &
\textbf{Decision:} include context\par\smallskip \textbf{Placement:} Main structural-reduction context &
Included for a positive structural-wrapper benefit without quantum-component attribution. The subspace representation is exact, but selecting the subspace need not preserve the global optimum. The journal release on 2026-10-02 follows earlier public versions. Read scope: Section III kernel selection; Section VI noiseless VQE experiments on 200 lattices and resource comparisons; and Section VII classical applicability. \\
\addlinespace[4pt]
\textbf{Candidate:} 2607.28120\par\smallskip Approximate sampling from decoded quantum interferometry via Markov chain Monte Carlo methods\par\smallskip gilfuster2026mcmc\par\smallskip \url{https://arxiv.org/abs/2607.28120} &
\textbf{First public:} 2026-07-30\par\smallskip \textbf{Version inspected:} arXiv:2607.28120 v1\par\smallskip \textbf{Version date:} 2026-07-30\par\smallskip \textbf{Checked:} 2026-10-04 &
Targeted constraint/theory searches and citation chaining; extension/constraints-theory-search.json &
\textbf{Decision:} include context\par\smallskip \textbf{Placement:} Main theoretical/method context &
Included as empirical evidence for classical approximation of the target performance. The finite-size fit establishes neither a polynomial-time simulation theorem nor a measured practical crossover. Eligibility follows the arXiv v1 submission history dated 2026-07-30, rather than the August 24 date displayed in the HTML body. Read scope: claim-critical statements and methods; theorem proofs were not independently replicated. \\
\addlinespace[4pt]
\textbf{Candidate:} 2609.18209\par\smallskip The Expressive Power of Constrained QAOA: What You Might Have MISsed\par\smallskip tsvelikhovskiy2026expressive\par\smallskip \url{https://arxiv.org/abs/2609.18209} &
\textbf{First public:} 2026-09-16\par\smallskip \textbf{Version inspected:} arXiv:2609.18209 v1\par\smallskip \textbf{Version date:} 2026-09-16\par\smallskip \textbf{Checked:} 2026-10-04 &
Targeted constraint/theory searches and citation chaining; extension/constraints-theory-search.json &
\textbf{Decision:} include context\par\smallskip \textbf{Placement:} Main theoretical/method context &
Included for the controllability and expressivity theorem. Existence of an exact optimal-state preparation does not establish an efficient optimization algorithm. Read scope: claim-critical statements and methods; theorem proofs were not independently replicated. \\
\addlinespace[4pt]
\textbf{Candidate:} 2610.00736\par\smallskip Beyond Feasibility: Finite-Depth Accessibility in Constrained QAOA\par\smallskip ubale2026accessibility\par\smallskip \url{https://arxiv.org/abs/2610.00736} &
\textbf{First public:} 2026-09-30\par\smallskip \textbf{Version inspected:} arXiv:2610.00736 v1\par\smallskip \textbf{Version date:} 2026-09-30\par\smallskip \textbf{Checked:} 2026-10-04 &
Targeted constraint/theory searches and citation chaining; extension/constraints-theory-search.json &
\textbf{Decision:} include context\par\smallskip \textbf{Placement:} Main theoretical/method context &
Included for the ideal support restriction and controlled small-instance simulations. State accessibility bounds which states can be reached, rather than guaranteeing solution quality. The public listing is dated 2026-10-02; the submission date is recorded separately. Read scope: claim-critical statements and methods; theorem proofs were not independently replicated. \\
\addlinespace[4pt]
\textbf{Candidate:} 2610.02145\par\smallskip A provable quantum advantage for approximate optimization via decoded quantum interferometry\par\smallskip kramer2026provable\par\smallskip \url{https://arxiv.org/abs/2610.02145} &
\textbf{First public:} 2026-10-01\par\smallskip \textbf{Version inspected:} arXiv:2610.02145 v1 full text; v2 abstract and submission metadata\par\smallskip \textbf{Version date:} 2026-10-01\par\smallskip \textbf{Checked:} 2026-10-04 &
Targeted constraint/theory searches and citation chaining; extension/constraints-theory-search.json &
\textbf{Decision:} include context\par\smallskip \textbf{Placement:} Main theoretical/method context &
Included for an oracle-model quantum separation, without measured hardware or practical hybrid-system advantage. Claim-critical full text was inspected in v1. The latest v2 revision is dated 2026-10-02, but only its abstract and submission metadata were accessible: versioned HTML and PDF fetches failed, and unversioned HTML identified itself as v1. The citation preserves the inspected version; theorem proofs were not independently replicated. \\
\addlinespace[4pt]
\textbf{Candidate:} 2607.14650\par\smallskip Worst-Case Quantum Algorithm for Optimal Polynomial Intersection Beyond Decoded Quantum Interferometry\par\smallskip \url{https://arxiv.org/abs/2607.14650} &
\textbf{First public:} Not captured by this audit.\par\smallskip \textbf{Version inspected:} Primary record/abstract; version not extracted\par\smallskip \textbf{Version date:} Not captured by this audit.\par\smallskip \textbf{Checked:} 2026-10-04 &
Constraint/theory searches and citation chaining; extension/constraints-theory-search.json &
\textbf{Decision:} screen only\par\smallskip \textbf{Placement:} Discovery register only &
Retained as an OPI algorithm lead identified through the Kramer bibliography. The primary record was inspected, but the full proof was not audited. The bounded contextual selection emphasizes access-model distinctions and classical controls. Metadata comes from retained primary-record notes; calendar dates not retained in those notes remain blank. \\
\addlinespace[4pt]
\textbf{Candidate:} 2607.16541\par\smallskip Efficient Exact Quantum Sampling from the Sun-Wootters Distribution for Optimal Polynomial Intersection\par\smallskip \url{https://arxiv.org/abs/2607.16541} &
\textbf{First public:} Not captured by this audit.\par\smallskip \textbf{Version inspected:} Primary record/abstract; version not extracted\par\smallskip \textbf{Version date:} Not captured by this audit.\par\smallskip \textbf{Checked:} 2026-10-04 &
Constraint/theory searches and citation chaining; extension/constraints-theory-search.json &
\textbf{Decision:} screen only\par\smallskip \textbf{Placement:} Discovery register only &
Retained as an exact quantum OPI sampling lead identified through the Kramer bibliography. The primary record was inspected, but the full proof was not audited. Its omission from the main text is a scope decision, not a correctness assessment. Metadata comes from retained primary-record notes; calendar dates not retained in those notes remain blank. \\
\addlinespace[4pt]
\textbf{Candidate:} 2608.01121\par\smallskip Polynomial Time Quantum Approximation Schemes for Constrained Optimisation\par\smallskip \url{https://arxiv.org/abs/2608.01121} &
\textbf{First public:} 2026-08-02\par\smallskip \textbf{Version inspected:} arXiv:2608.01121v3 (claim-critical full-text screen)\par\smallskip \textbf{Version date:} 2026-08-06\par\smallskip \textbf{Checked:} 2026-10-04 &
Constraint/theory searches and citation chaining; extension/constraints-theory-search.json &
\textbf{Decision:} screen only\par\smallskip \textbf{Placement:} Discovery register only &
Retained after a claim-critical v3 screen. The approximation and optimal-hit claims depend on instance-specific inverse-polynomial optimal mass and a kernel-admissible promise family. These premises were not verified sufficiently to support broader efficiency language in this survey. The decision does not assert that the paper is incorrect. Bibliographic metadata comes from retained primary-record notes. \\
\addlinespace[4pt]
\textbf{Candidate:} 2608.07717\par\smallskip Implicit Differentiation for Measurement-Efficient Bilevel Quantum-Classical Optimization\par\smallskip \url{https://arxiv.org/abs/2608.07717} &
\textbf{First public:} 2026-08-07\par\smallskip \textbf{Version inspected:} Primary record/abstract; version not extracted\par\smallskip \textbf{Version date:} Not captured by this audit.\par\smallskip \textbf{Checked:} 2026-10-04 &
Constraint/theory searches and citation chaining; extension/constraints-theory-search.json &
\textbf{Decision:} screen only\par\smallskip \textbf{Placement:} Discovery register only &
Retained as a potentially relevant method for reusing measurements in bilevel optimization. Only the abstract was screened; a claim-critical full-text audit was not completed, so the study was not prioritized for the bounded contextual selection. Bibliographic metadata comes from retained primary-record notes; unretained dates remain blank. \\
\addlinespace[4pt]
\textbf{Candidate:} 2609.36999\par\smallskip Towards Scalable Quantum Unit Commitment in Power Systems: Encoding and Complexity Reduction Techniques\par\smallskip \url{https://arxiv.org/abs/2609.36999} &
\textbf{First public:} Not captured by this audit.\par\smallskip \textbf{Version inspected:} Not captured by this audit.\par\smallskip \textbf{Version date:} Not captured by this audit.\par\smallskip \textbf{Checked:} 2026-10-04 &
Targeted comparative searches; extension/comparative-search.json &
\textbf{Decision:} screen only\par\smallskip \textbf{Placement:} Discovery register only &
Retained as a secondary discovery for future follow-up. No primary metadata or empirical audit was completed in this search, so the record does not support size or performance claims in the survey. \\
\addlinespace[4pt]
\textbf{Candidate:} VOCAL2026\_41\par\smallskip A quantum-classical branch-and-bound algorithm\par\smallskip \url{https://www.euro-online.org/conf/vocal2026/treat_abstract?paperid=41} &
\textbf{First public:} Not captured by this audit.\par\smallskip \textbf{Version inspected:} Not captured by this audit.\par\smallskip \textbf{Version date:} Not captured by this audit.\par\smallskip \textbf{Checked:} 2026-10-04 &
Targeted decomposition screening; extension/decomposition-search.json &
\textbf{Decision:} screen only\par\smallskip \textbf{Placement:} Discovery register only &
Retained as a conference-abstract lead. No full paper or version date within the extension window was verified, so the abstract was not used as new guarantee evidence. \\
\addlinespace[4pt]
\textbf{Candidate:} doi:10.1103/568y-4264\par\smallskip Qubit-efficient variational quantum optimization via Pauli correlation encoding: Application to large-scale power demand portfolio optimization\par\smallskip \url{https://journals.aps.org/prapplied/abstract/10.1103/568y-4264} &
\textbf{First public:} Not captured by this audit.\par\smallskip \textbf{Version inspected:} Publisher record and abstract\par\smallskip \textbf{Version date:} 2026-10-01\par\smallskip \textbf{Checked:} 2026-10-04 &
Publisher article-date screening; targeted/publisher-search.json &
\textbf{Decision:} screen only\par\smallskip \textbf{Placement:} Discovery register only &
The October 1 journal release falls within the expanded July--October window. The abstract reports compressed portfolio encoding and a hardware demonstration, but controls and timing were not audited in full. The study remains a screened lead. Read scope: publisher abstract only; the first public preprint date was not extracted. \\
\addlinespace[4pt]
\textbf{Candidate:} 2404.17703\par\smallskip Quantum annealing approach to the minimum distance problem of quantum codes\par\smallskip ismail2026minimumdistance\par\smallskip \url{https://journals.aps.org/pra/accepted/10.1103/573l-6cc2} &
\textbf{First public:} 2024-04-26\par\smallskip \textbf{Version inspected:} arXiv:2404.17703v3; 2026 accepted abstract\par\smallskip \textbf{Version date:} 2024-07-03\par\smallskip \textbf{Checked:} 2026-10-04 &
Targeted comparative searches; extension/comparative-search.json &
\textbf{Decision:} screen only\par\smallskip \textbf{Placement:} Discovery register only &
Qbsolv with QA and SA subsolvers gives broadly comparable results under equal end-user time limits, with unknown remote latency affecting QA. The newly accepted abstract emphasizes solver-independent reformulation and embedding limits more cautiously than the 2024 preprint abstract. Heuristic upper bounds do not certify exact code distance. Read scope: the accepted abstract and v3 (2024), Sections V.1--V.2. DOI 10.1103/573l-6cc2 was verified, but volume, article number and publication date were not. The publisher spells the first author Refaat Ismail; the preprint uses Refat Ismail. \\
\addlinespace[4pt]
\textbf{Candidate:} 2509.15138\par\smallskip Sampled-Based Guided Quantum Walk: Non-variational quantum algorithm for combinatorial optimization\par\smallskip \url{https://arxiv.org/abs/2509.15138} &
\textbf{First public:} 2025-09-18\par\smallskip \textbf{Version inspected:} arXiv:2509.15138v2 (abstract and record)\par\smallskip \textbf{Version date:} 2026-10-01\par\smallskip \textbf{Checked:} 2026-10-04 &
October boundary category listings and primary-record follow-up; targeted/arxiv-search.json &
\textbf{Decision:} screen only\par\smallskip \textbf{Placement:} Discovery register only &
Retained as an updated method using offline classical spectral sampling to guide quantum walks without an online variational optimizer. The abstract reports binary quadratic and higher-order benchmarks up to 30 qubits. Empirical changes between v1 and v2 were not audited. Read scope: primary abstract and submission history. The 2026-10-02 announcement date is distinct from the version date. \\
\addlinespace[4pt]
\textbf{Candidate:} 2511.22251\par\smallskip A Framework for Handling and Exploiting Symmetry in Benders Decomposition\par\smallskip \url{https://arxiv.org/abs/2511.22251} &
\textbf{First public:} 2025-11-27\par\smallskip \textbf{Version inspected:} arXiv:2511.22251v4 (abstract and record)\par\smallskip \textbf{Version date:} 2026-10-01\par\smallskip \textbf{Checked:} 2026-10-04 &
October boundary category listings and primary-record follow-up; targeted/arxiv-search.json &
\textbf{Decision:} screen only\par\smallskip \textbf{Placement:} Discovery register only &
Retained as classical background on symmetry handling, aggregation and Benders-cut strengthening to reduce oracle calls. It does not study a quantum Benders component. Empirical changes between versions were not audited. Read scope: primary abstract and submission history. The 2026-10-02 announcement date is distinct from the version date. \\
\addlinespace[4pt]
\textbf{Candidate:} 2601.08085\par\smallskip Learning parameter curves in feedback-based quantum algorithms\par\smallskip \url{https://arxiv.org/abs/2601.08085} &
\textbf{First public:} 2026-01-13\par\smallskip \textbf{Version inspected:} arXiv:2601.08085v2 (abstract and record)\par\smallskip \textbf{Version date:} 2026-10-01\par\smallskip \textbf{Checked:} 2026-10-04 &
October boundary category listings and primary-record follow-up; targeted/arxiv-search.json &
\textbf{Decision:} screen only\par\smallskip \textbf{Placement:} Discovery register only &
Retained as an updated method using classical learning to predict circuit parameter curves and reduce measurement overhead, with weighted MaxCut and transverse-field Ising simulations. The primary record identifies v2 with the published article, but empirical changes were not audited individually. Read scope: primary abstract and submission history. The 2026-10-02 announcement date is distinct from the version date. \\
\addlinespace[4pt]
\textbf{Candidate:} 2604.07218\par\smallskip Improving Feasibility in Quantum Approximate Optimization Algorithm for Vehicle Routing via Constraint-Aware Initialization and Hybrid XY-X Mixing\par\smallskip \url{https://arxiv.org/abs/2604.07218} &
\textbf{First public:} 2026-04-08\par\smallskip \textbf{Version inspected:} arXiv:2604.07218v2 (abstract and record)\par\smallskip \textbf{Version date:} 2026-09-30\par\smallskip \textbf{Checked:} 2026-10-04 &
October boundary category listings and primary-record follow-up; targeted/arxiv-search.json &
\textbf{Decision:} screen only\par\smallskip \textbf{Placement:} Discovery register only &
Retained for its combination of structured initialization and hybrid XY-X mixing, compared with standard QAOA in ideal, shot-based and noisy simulations. This is a mechanism comparison rather than a classical-system benchmark. Read scope: primary abstract and submission history; no full empirical reaudit. The revision is dated September 30 and the listing October 2. ArXiv reports journal DOI 10.1016/j.trc.2026.106047, but the publisher release date was not verified. \\
\addlinespace[4pt]
\textbf{Candidate:} 2607.07336\par\smallskip Resource-Efficient Hybrid Quantum Neighborhood Selection for Large-Scale Molecular Diversity Optimization\par\smallskip \url{https://arxiv.org/abs/2607.07336} &
\textbf{First public:} 2026-07-08\par\smallskip \textbf{Version inspected:} arXiv:2607.07336v1\par\smallskip \textbf{Version date:} 2026-07-08\par\smallskip \textbf{Checked:} 2026-10-04 &
Followed the July companion recorded in the earlier publication-status audit; reread methods V-C, V-F, V-H and tables II-III. &
\textbf{Decision:} screen only\par\smallskip \textbf{Placement:} Discovery register only &
Retained as a screened companion to the cited HQNS mechanism. Claimed resource savings compare different workflows and stopping budgets without a matched classical frontier-solver substitution. Section V-H identifies no specific physical device or simulation backend, leaving the reported QPU timing boundary unverified. Quality is measured against observed annealing outputs rather than certified optima. The paper is eligible under the expanded window, but these limits prevent its use as additional hardware-advantage evidence. \\
\addlinespace[4pt]
\textbf{Candidate:} 2607.08708\par\smallskip Absence of quantum advantage for approximate spin glass optimization\par\smallskip \url{https://arxiv.org/abs/2607.08708} &
\textbf{First public:} 2026-07-09\par\smallskip \textbf{Version inspected:} Primary abstract; version not extracted\par\smallskip \textbf{Version date:} Not captured by this audit.\par\smallskip \textbf{Checked:} 2026-10-04 &
Targeted comparative searches; extension/comparative-search.json &
\textbf{Decision:} screen only\par\smallskip \textbf{Placement:} Discovery register only &
Retained as a narrow theoretical and numerical lead: semiclassical SK optimization can reproduce or improve the tested quantum schedules. The model and schedule restrictions preclude a general negative conclusion about quantum optimization. Read scope: primary abstract only. \\
\addlinespace[4pt]
\textbf{Candidate:} 2607.11637\par\smallskip From Circuits to Hardware: Benchmarking Standard and Qubit-Efficient Quantum Optimization on Real Hardware\par\smallskip sharma2026circuits\par\smallskip \url{https://arxiv.org/abs/2607.11637} &
\textbf{First public:} 2026-07-13\par\smallskip \textbf{Version inspected:} arXiv:2607.11637v2 (abstract and record)\par\smallskip \textbf{Version date:} 2026-08-07\par\smallskip \textbf{Checked:} 2026-10-04 &
Targeted comparative searches; extension/comparative-search.json &
\textbf{Decision:} screen only\par\smallskip \textbf{Placement:} Discovery register only &
Retained as a distinct paper from Adaptive Graph Shrinking (arXiv:2506.14250). Its abstract reports 247 method--instance combinations across four optimization families with matched random controls; one finite-sample warm-start MIS exception limits any blanket random-equivalence conclusion. Width reduction alone does not establish a performance benefit. Read scope: primary arXiv abstract and metadata; the official IOP page was blocked by robots restrictions. ArXiv verifies DOI 10.1088/2058-9565/ae94a4; volume 11, article 045029 and the September 1 publication date were only secondary-indexed in this audit. \\
\addlinespace[4pt]
\textbf{Candidate:} 2608.15655\par\smallskip Rand-SEMI-QAOA: Finite-Budget Depth-One MaxCut Ensembles on Compressed Quantum Registers\par\smallskip \url{https://arxiv.org/abs/2608.15655} &
\textbf{First public:} 2026-08-16\par\smallskip \textbf{Version inspected:} Primary abstract; version not extracted\par\smallskip \textbf{Version date:} Not captured by this audit.\par\smallskip \textbf{Checked:} 2026-10-04 &
Targeted comparative searches; extension/comparative-search.json &
\textbf{Decision:} screen only\par\smallskip \textbf{Placement:} Discovery register only &
Retained as a compressed depth-one MaxCut ensemble method on small random regular instances. Frozen selector-pool comparisons do not by themselves establish quantum-component or system superiority. Read scope: primary abstract only. \\
\addlinespace[4pt]
\textbf{Candidate:} 2608.30282\par\smallskip Best Annealing Path of Quantum Annealing via Efficient Adiabatic Phase Transition\par\smallskip \url{https://arxiv.org/abs/2608.30282} &
\textbf{First public:} 2026-08-31\par\smallskip \textbf{Version inspected:} Primary abstract; version not extracted\par\smallskip \textbf{Version date:} Not captured by this audit.\par\smallskip \textbf{Checked:} 2026-10-04 &
Targeted comparative searches; extension/comparative-search.json &
\textbf{Decision:} screen only\par\smallskip \textbf{Placement:} Discovery register only &
Retained as a simulation and annealing-mechanism lead, without an audited quantum-component attribution comparison. Read scope: primary abstract only. \\
\addlinespace[4pt]
\textbf{Candidate:} 2609.19392\par\smallskip Transformers as Intrinsic Optimizers for Quantum Approximate Optimization Algorithm\par\smallskip \url{https://arxiv.org/abs/2609.19392} &
\textbf{First public:} 2026-09-16\par\smallskip \textbf{Version inspected:} Primary abstract; version not extracted\par\smallskip \textbf{Version date:} Not captured by this audit.\par\smallskip \textbf{Checked:} 2026-10-04 &
Targeted comparative searches; extension/comparative-search.json &
\textbf{Decision:} screen only\par\smallskip \textbf{Placement:} Discovery register only &
Retained as a classical outer-optimizer enhancement for QAOA. It was not prioritized for the bounded contextual selection. Read scope: primary abstract only. \\
\addlinespace[4pt]
\textbf{Candidate:} 2609.39922\par\smallskip Fault-tolerant cost of shallow QAOA on near-symmetric optimization problems\par\smallskip \url{https://arxiv.org/abs/2609.39922} &
\textbf{First public:} 2026-09-30\par\smallskip \textbf{Version inspected:} Primary abstract; version not extracted\par\smallskip \textbf{Version date:} Not captured by this audit.\par\smallskip \textbf{Checked:} 2026-10-04 &
Targeted comparative searches; extension/comparative-search.json &
\textbf{Decision:} screen only\par\smallskip \textbf{Placement:} Discovery register only &
Retained as a theoretical contextual lead. Claim-critical methods and results were not audited for inclusion in the main text. Read scope: primary abstract only. \\
\addlinespace[4pt]
\textbf{Candidate:} 2609.40036\par\smallskip Cycle Codes and Decoded Quantum Interferometry\par\smallskip apte2026cycle\par\smallskip \url{https://arxiv.org/abs/2609.40036} &
\textbf{First public:} 2026-09-30\par\smallskip \textbf{Version inspected:} arXiv:2609.40036 v1\par\smallskip \textbf{Version date:} 2026-09-30\par\smallskip \textbf{Checked:} 2026-10-04 &
Targeted constraint/theory searches and citation chaining; extension/constraints-theory-search.json &
\textbf{Decision:} screen only\par\smallskip \textbf{Placement:} Discovery register only &
The full-text screen identified a restricted-family theoretical limitation and empirical classical controls. These results do not extend automatically to all DQI algorithms, Max-Cut families or the folded-OPI oracle model. The study was not selected for the bounded main-text expansion; this is not a correctness rejection. Read scope: claim-critical statements and methods; theorem proofs were not independently replicated. \\
\addlinespace[4pt]
\textbf{Candidate:} 2610.00502\par\smallskip Quantum Algorithms for OPI Variants Beyond Locality and Classical Decodability\par\smallskip \url{https://arxiv.org/abs/2610.00502} &
\textbf{First public:} 2026-09-30\par\smallskip \textbf{Version inspected:} arXiv:2610.00502v1 (abstract and record)\par\smallskip \textbf{Version date:} 2026-09-30\par\smallskip \textbf{Checked:} 2026-10-04 &
October boundary category listings and primary-record follow-up; targeted/arxiv-search.json &
\textbf{Decision:} screen only\par\smallskip \textbf{Placement:} Discovery register only &
Retained as adjacent theory extending coding and decoding approaches to OPI variants with global constraints. The abstract contains no empirical decomposition or system benchmark. Read scope: primary abstract and submission history; no full empirical audit. The 2026-10-02 announcement date is distinct from the submission date. \\
\addlinespace[4pt]
\textbf{Candidate:} 2610.00887\par\smallskip A Gate-Based Quantum Computing Framework for Codon Optimization\par\smallskip \url{https://arxiv.org/abs/2610.00887} &
\textbf{First public:} 2026-10-01\par\smallskip \textbf{Version inspected:} arXiv:2610.00887v1 (abstract and record)\par\smallskip \textbf{Version date:} 2026-10-01\par\smallskip \textbf{Checked:} 2026-10-04 &
October boundary category listings and primary-record follow-up; targeted/arxiv-search.json &
\textbf{Decision:} screen only\par\smallskip \textbf{Placement:} Discovery register only &
Retained as an application proof of concept comparing gate-based codon optimization with classical and annealing references, including IBM hardware. No decomposition contribution was apparent from the abstract. Read scope: primary abstract and submission history; no full empirical audit. The 2026-10-02 announcement date is distinct from the submission date. \\
\addlinespace[4pt]
\textbf{Candidate:} 2610.01253\par\smallskip Context-Aware Error Mitigation Orchestration for Hybrid Quantum Reinforcement Learning on NISQ Systems\par\smallskip \url{https://arxiv.org/abs/2610.01253} &
\textbf{First public:} 2026-10-01\par\smallskip \textbf{Version inspected:} arXiv:2610.01253v1 (claim-critical full text)\par\smallskip \textbf{Version date:} 2026-10-01\par\smallskip \textbf{Checked:} 2026-10-04 &
October boundary category listings and primary-record follow-up; targeted/arxiv-search.json; targeted/method-candidates.json full-text follow-up &
\textbf{Decision:} screen only\par\smallskip \textbf{Placement:} Discovery register only &
The inspected method uses Qiskit Aer simulations with 4--7 qubits for 15--100 customers; hardware times are projections. Results support classical mitigation-controller choices under the simulated noise model, while quantum-specific value remains unresolved. Table 2 fidelity uses ideal output probabilities, which would require a practical hardware proxy. No matched classical-policy substitution isolates the quantum contribution, and projected timing does not establish measured system advantage. \\
\addlinespace[4pt]
\textbf{Candidate:} 2610.01678\par\smallskip Safe Hypergraph Contraction via Capacity-Aware Repair Certificates\par\smallskip \url{https://arxiv.org/abs/2610.01678} &
\textbf{First public:} 2026-10-01\par\smallskip \textbf{Version inspected:} arXiv:2610.01678v1 (abstract and record)\par\smallskip \textbf{Version date:} 2026-10-01\par\smallskip \textbf{Checked:} 2026-10-04 &
October boundary category listings and primary-record follow-up; targeted/arxiv-search.json &
\textbf{Decision:} screen only\par\smallskip \textbf{Placement:} Discovery register only &
Retained as classical background on certified hypergraph coarsening under specified capacity and weight conditions. It informs transformation guarantees but does not test a quantum routine. Read scope: primary abstract and submission history; no full empirical audit. The 2026-10-02 announcement date is distinct from the submission date. \\
\addlinespace[4pt]
\textbf{Candidate:} 2610.01885\par\smallskip A two-stage approach to satellite constellation optimization: classical and QUBO formulations\par\smallskip \url{https://arxiv.org/abs/2610.01885} &
\textbf{First public:} 2026-10-01\par\smallskip \textbf{Version inspected:} arXiv:2610.01885v1 (claim-critical full text)\par\smallskip \textbf{Version date:} 2026-10-01\par\smallskip \textbf{Checked:} 2026-10-04 &
October boundary category listings and primary-record follow-up; targeted/arxiv-search.json; targeted/method-candidates.json full-text follow-up &
\textbf{Decision:} screen only\par\smallskip \textbf{Placement:} Discovery register only &
Section 7 reports continuous classical particle-swarm experiments in MATLAB. Discretized and quantum-inspired tests remained underway, with no quantum hardware or QUBO execution reported. The study therefore supplies preliminary classical-staging evidence and a proposed quantum-compatible formulation, without an isolated quantum-versus-classical component comparison. \\
\addlinespace[4pt]
\textbf{Candidate:} 2610.01902\par\smallskip Exponential quantum advantages for decoded quantum interferometry in the streaming setting\par\smallskip \url{https://arxiv.org/abs/2610.01902} &
\textbf{First public:} 2026-10-01\par\smallskip \textbf{Version inspected:} arXiv:2610.01902v1 (abstract and record)\par\smallskip \textbf{Version date:} 2026-10-01\par\smallskip \textbf{Checked:} 2026-10-04 &
October boundary category listings and primary-record follow-up; targeted/arxiv-search.json &
\textbf{Decision:} screen only\par\smallskip \textbf{Placement:} Discovery register only &
Retained as contextual theory on a DQI memory separation under streaming access and resource restrictions. This model differs from ordinary elapsed-time comparisons of hybrid optimizers. Read scope: primary abstract and submission history; no full empirical audit. The 2026-10-02 announcement date is distinct from the submission date. \\
\addlinespace[4pt]
\textbf{Candidate:} 2607.19581\par\smallskip Quantum-Resilient Distributed Optimization for Multi-Region Unit Commitment\par\smallskip \url{https://arxiv.org/abs/2607.19581} &
\textbf{First public:} Not captured by this audit.\par\smallskip \textbf{Version inspected:} Primary publisher/discovery record\par\smallskip \textbf{Version date:} 2026-07-21\par\smallskip \textbf{Checked:} 2026-10-04 &
Targeted decomposition screening; extension/decomposition-search.json &
\textbf{Decision:} exclude scope\par\smallskip \textbf{Placement:} Discovery register only &
Excluded from the extension because the method concerns post-quantum-secure communication in classical Benders decomposition, rather than a quantum optimization component. \\
\addlinespace[4pt]
\textbf{Candidate:} 2608.07047\par\smallskip Stationarity is not enough: tightness of the quantum mechanical bootstrap and the copositive cone\par\smallskip \url{https://arxiv.org/abs/2608.07047} &
\textbf{First public:} Not captured by this audit.\par\smallskip \textbf{Version inspected:} Not captured by this audit.\par\smallskip \textbf{Version date:} Not captured by this audit.\par\smallskip \textbf{Checked:} 2026-10-04 &
Targeted decomposition screening; extension/decomposition-search.json &
\textbf{Decision:} exclude scope\par\smallskip \textbf{Placement:} Discovery register only &
Excluded from the extension because the discovery record concerns quantum-state bootstrap and cone analysis without an identified direct hybrid combinatorial optimization method. \\
\addlinespace[4pt]
\textbf{Candidate:} GAS\_BD\par\smallskip Grover Adaptive Search Based Hybrid Benders Decomposition for Mixed-Integer Linear Programs\par\smallskip \url{https://tqe.ieee.org/2026/04/06/grover-adaptive-search-based-hybrid-benders-decomposition-for-mixed-integer-linear-programs/} &
\textbf{First public:} Not captured by this audit.\par\smallskip \textbf{Version inspected:} Primary publisher/discovery record\par\smallskip \textbf{Version date:} 2026-04-06\par\smallskip \textbf{Checked:} 2026-10-04 &
Targeted decomposition screening; extension/decomposition-search.json &
\textbf{Decision:} exclude scope\par\smallskip \textbf{Placement:} Discovery register only &
Excluded only from the new or revised extension: the verified publisher record is dated April, with no verified July--October release. \\
\addlinespace[4pt]
\textbf{Candidate:} doi:10.1007/s10791-026-10646-4\par\smallskip Energy-efficient clustering and mobile sink routing in wireless sensor networks using quantum-inspired optimization and deep learning\par\smallskip \url{https://link.springer.com/article/10.1007/s10791-026-10646-4} &
\textbf{First public:} Not captured by this audit.\par\smallskip \textbf{Version inspected:} Publisher record and abstract\par\smallskip \textbf{Version date:} 2026-10-03\par\smallskip \textbf{Checked:} 2026-10-04 &
Publisher article-date screening; targeted/publisher-search.json &
\textbf{Decision:} exclude scope\par\smallskip \textbf{Placement:} Discovery register only &
Excluded from the selected scope after a publisher-record and abstract screen. The wireless-sensor application combines quantum-inspired optimization and deep learning, but the abstract identifies no directly relevant decomposition or controlled native quantum-optimizer comparison. Hardware execution or quantum advantage was not inferred from terminology. No detailed empirical audit was performed, and the first public preprint date was not extracted. \\
\addlinespace[4pt]
\textbf{Candidate:} doi:10.1038/s40494-026-02863-3\par\smallskip A mixed-culture quantitative model of the restoration priority of traditional dwellings using Grover's search\par\smallskip \url{https://www.nature.com/articles/s40494-026-02863-3} &
\textbf{First public:} Not captured by this audit.\par\smallskip \textbf{Version inspected:} Publisher record and abstract\par\smallskip \textbf{Version date:} 2026-10-02\par\smallskip \textbf{Checked:} 2026-10-04 &
Publisher article-date screening; targeted/publisher-search.json &
\textbf{Decision:} exclude scope\par\smallskip \textbf{Placement:} Discovery register only &
Excluded from the selected scope after screening the indexed primary abstract and publisher subject listing. The Grover-filtering application uses a Boolean oracle and a 143-dwelling dataset, but the available material establishes no decomposition or controlled quantum-component comparison relevant to the evidence tables. Direct article and PDF access failed; the first public preprint date was not extracted. \\
\addlinespace[4pt]
\textbf{Candidate:} 2506.14250\par\smallskip Adaptive Graph Shrinking for Quantum Optimization of Constrained Combinatorial Problems\par\smallskip sharma2025adaptive\par\smallskip \url{https://arxiv.org/abs/2506.14250} &
\textbf{First public:} 2025-06-17\par\smallskip \textbf{Version inspected:} Primary record v1 (status only)\par\smallskip \textbf{Version date:} 2025-06-17\par\smallskip \textbf{Checked:} 2026-10-04 &
Targeted decomposition screening; extension/decomposition-search.json &
\textbf{Decision:} exclude scope\par\smallskip \textbf{Placement:} Existing survey record retained; no new extension entry &
Retained in the original corpus; excluded only from the new or revised extension because no revision within the extension window was identified on the primary record. This was a publication-status and version check, without a full numerical reaudit. \\
\addlinespace[4pt]
\textbf{Candidate:} 2601.14024\par\smallskip Performance enhancing of hybrid quantum-classical Benders approach for MILP optimization\par\smallskip lopezbanos2026performance\par\smallskip \url{https://arxiv.org/abs/2601.14024} &
\textbf{First public:} 2026-01-20\par\smallskip \textbf{Version inspected:} Primary record v3 (status only)\par\smallskip \textbf{Version date:} 2026-06-30\par\smallskip \textbf{Checked:} 2026-10-04 &
Targeted decomposition screening; extension/decomposition-search.json &
\textbf{Decision:} exclude scope\par\smallskip \textbf{Placement:} Existing survey record retained; no new extension entry &
Retained in the original corpus; excluded only from the new or revised extension because no revision within the extension window was identified. The latest primary version is dated 2026-06-30. This was a publication-status and version check, without a full numerical reaudit. \\
\addlinespace[4pt]
\textbf{Candidate:} 2603.28933\par\smallskip Iterative Optimization with Partial Convergence Guarantees on Neutral Atom Quantum Computers\par\smallskip perron2026iterative\par\smallskip \url{https://arxiv.org/abs/2603.28933} &
\textbf{First public:} 2026-03-30\par\smallskip \textbf{Version inspected:} Primary record v1 (status only)\par\smallskip \textbf{Version date:} 2026-03-30\par\smallskip \textbf{Checked:} 2026-10-04 &
Targeted decomposition screening; extension/decomposition-search.json &
\textbf{Decision:} exclude scope\par\smallskip \textbf{Placement:} Existing survey record retained; no new extension entry &
Retained in the original corpus; excluded only from the new or revised extension because no revision within the extension window was identified on the primary record. This was a publication-status and version check, without a full numerical reaudit. \\
\addlinespace[4pt]
\textbf{Candidate:} 2610.03161\par\smallskip Landscape-Dependent Performance of Photonic Quantum Solvers in QUBO Feature Selection for Financial Risk Detection\par\smallskip \url{https://arxiv.org/abs/2610.03161} &
\textbf{First public:} Not captured by this audit.\par\smallskip \textbf{Version inspected:} arXiv:2610.03161v1\par\smallskip \textbf{Version date:} 2026-10-02\par\smallskip \textbf{Checked:} 2026-10-12 &
Public arXiv October/recent listing plus date-filtered topic searches; bibliography follow-up for older ADMM leads; discovery.json &
\textbf{Decision:} screen only\par\smallskip \textbf{Placement:} Discovery log only &
Photonic feature-selection comparison is adjacent context; title/abstract checked, but no decomposition mechanism is identified and direct preparation/sampling attribution is not extracted. Announcement date not recovered from the recent listing. Date basis: boundary or prior submission announcement not verified. Not independently verified; submission 2026-10-02. \\
\addlinespace[4pt]
\textbf{Candidate:} 2610.03360\par\smallskip Rate-Optimal Quantum Discrete Simulation Optimization\par\smallskip \url{https://arxiv.org/abs/2610.03360} &
\textbf{First public:} Not captured by this audit.\par\smallskip \textbf{Version inspected:} arXiv:2610.03360v1\par\smallskip \textbf{Version date:} 2026-10-02\par\smallskip \textbf{Checked:} 2026-10-12 &
Public arXiv October/recent listing plus date-filtered topic searches; bibliography follow-up for older ADMM leads; discovery.json &
\textbf{Decision:} screen only\par\smallskip \textbf{Placement:} Discovery log only &
Oracle-based discrete simulation optimization; relevant optimization theory, but not a decomposition component or implementation benchmark; abstract screen only. Date basis: boundary or prior submission announcement not verified. Not independently verified; submission 2026-10-02. \\
\addlinespace[4pt]
\textbf{Candidate:} 2610.03379\par\smallskip Degree-conditioned anneal offsets reshape certified-optimum sampling for minimum vertex cover\par\smallskip \url{https://arxiv.org/abs/2610.03379} &
\textbf{First public:} Not captured by this audit.\par\smallskip \textbf{Version inspected:} arXiv:2610.03379v1\par\smallskip \textbf{Version date:} 2026-10-02\par\smallskip \textbf{Checked:} 2026-10-12 &
Public arXiv October/recent listing plus date-filtered topic searches; bibliography follow-up for older ADMM leads; discovery.json &
\textbf{Decision:} screen only\par\smallskip \textbf{Placement:} Discovery log only &
Controlled anneal-offset intervention changes diversity of certified optima; relevant hardware-control context, but not decomposition attribution. Abstract screen only; no new main evidence row. Date basis: boundary or prior submission announcement not verified. Not independently verified; submission 2026-10-02. \\
\addlinespace[4pt]
\textbf{Candidate:} 2610.03788\par\smallskip The Quantum Sphere: A Physically Realizable Optimization Benchmark with Provable Linear Convergence in White- and Black-Box Settings\par\smallskip \url{https://arxiv.org/abs/2610.03788} &
\textbf{First public:} 2026-10-06\par\smallskip \textbf{Version inspected:} arXiv:2610.03788v1\par\smallskip \textbf{Version date:} 2026-09-30\par\smallskip \textbf{Checked:} 2026-10-12 &
Public arXiv October/recent listing plus date-filtered topic searches; bibliography follow-up for older ADMM leads; discovery.json &
\textbf{Decision:} exclude scope\par\smallskip \textbf{Placement:} Discovery log only &
Benchmark optimizes a quantum-control landscape rather than using a quantum component to solve a classical mathematical program. Date basis: in window announcement. \\
\addlinespace[4pt]
\textbf{Candidate:} 2610.03890\par\smallskip Quantum-Assisted Optimization Guided by Machine-Learned Risk Maps for Aerial Surveillance Routing\par\smallskip lauwerier2026quantum\par\smallskip \url{https://arxiv.org/abs/2610.03890} &
\textbf{First public:} 2026-10-06\par\smallskip \textbf{Version inspected:} arXiv:2610.03890v1\par\smallskip \textbf{Version date:} 2026-10-02\par\smallskip \textbf{Checked:} 2026-10-12 &
Public arXiv October/recent listing plus date-filtered topic searches; bibliography follow-up for older ADMM leads; discovery.json &
\textbf{Decision:} include evidence\par\smallskip \textbf{Placement:} E18; main evidence &
E18: full-text read establishes fixed-candidate-pool selection controls and a distinction between projected QPU and measured classical time; no systematic quantum component superiority. Date basis: in window announcement. \\
\addlinespace[4pt]
\textbf{Candidate:} 2610.03949\par\smallskip Benchmarking Ansatze for Pauli Correlation Encoding in the Maximum Independent Set Problem\par\smallskip \url{https://arxiv.org/abs/2610.03949} &
\textbf{First public:} 2026-10-06\par\smallskip \textbf{Version inspected:} arXiv:2610.03949v1\par\smallskip \textbf{Version date:} 2026-10-02\par\smallskip \textbf{Checked:} 2026-10-12 &
Public arXiv October/recent listing plus date-filtered topic searches; bibliography follow-up for older ADMM leads; discovery.json &
\textbf{Decision:} screen only\par\smallskip \textbf{Placement:} Discovery log only &
Compressed PCE ansatz/decoding comparison is relevant context, but adds no independent matched-resource solver comparison beyond the manuscript's existing compression caution; abstract screened. Date basis: in window announcement. \\
\addlinespace[4pt]
\textbf{Candidate:} 2610.04068\par\smallskip Interaction-Aware Embedding Optimization for Neutral-Atom Quantum Processors\par\smallskip \url{https://arxiv.org/abs/2610.04068} &
\textbf{First public:} 2026-10-06\par\smallskip \textbf{Version inspected:} arXiv:2610.04068v1\par\smallskip \textbf{Version date:} 2026-10-02\par\smallskip \textbf{Checked:} 2026-10-12 &
Public arXiv October/recent listing plus date-filtered topic searches; bibliography follow-up for older ADMM leads; discovery.json &
\textbf{Decision:} screen only\par\smallskip \textbf{Placement:} Discovery log only &
Classical force-directed neutral-atom embedding method; relevant hardware-interface context, not a quantum optimization-component comparison. Abstract screened. Date basis: in window announcement. \\
\addlinespace[4pt]
\textbf{Candidate:} 2610.05524\par\smallskip Optimization Geometry of QAOA and Variational Quantum Algorithms\par\smallskip \url{https://arxiv.org/abs/2610.05524} &
\textbf{First public:} 2026-10-06\par\smallskip \textbf{Version inspected:} arXiv:2610.05524v1\par\smallskip \textbf{Version date:} 2026-10-04\par\smallskip \textbf{Checked:} 2026-10-12 &
Public arXiv October/recent listing plus date-filtered topic searches; bibliography follow-up for older ADMM leads; discovery.json &
\textbf{Decision:} screen only\par\smallskip \textbf{Placement:} Discovery log only &
Classical optimizer comparison within VQA training, using function-evaluation-matched multistart BFGS; does not isolate a quantum contribution to solving the classical problem. Abstract screened. Date basis: in window announcement. \\
\addlinespace[4pt]
\textbf{Candidate:} 2610.06000\par\smallskip Metastability in Quantum Dynamics for Optimization\par\smallskip \url{https://arxiv.org/abs/2610.06000} &
\textbf{First public:} 2026-10-06\par\smallskip \textbf{Version inspected:} arXiv:2610.06000v1\par\smallskip \textbf{Version date:} 2026-10-05\par\smallskip \textbf{Checked:} 2026-10-12 &
Public arXiv October/recent listing plus date-filtered topic searches; bibliography follow-up for older ADMM leads; discovery.json &
\textbf{Decision:} screen only\par\smallskip \textbf{Placement:} Discovery log only &
Quantum Langevin dynamics and metastability in continuous potentials; theoretical mechanism outside the focused decomposition update. Abstract screened. Date basis: in window announcement. \\
\addlinespace[4pt]
\textbf{Candidate:} 2610.06448\par\smallskip A Practical Introduction to VQE: Methods, Challenges, and Applications in Physics and Chemistry\par\smallskip \url{https://arxiv.org/abs/2610.06448} &
\textbf{First public:} 2026-10-06\par\smallskip \textbf{Version inspected:} arXiv:2610.06448v2\par\smallskip \textbf{Version date:} 2026-10-07\par\smallskip \textbf{Checked:} 2026-10-12 &
Public arXiv October/recent listing plus date-filtered topic searches; bibliography follow-up for older ADMM leads; discovery.json &
\textbf{Decision:} exclude scope\par\smallskip \textbf{Placement:} Discovery log only &
VQE tutorial emphasizes electronic structure and materials simulation rather than classical constrained/combinatorial optimization. Date basis: in window announcement. \\
\addlinespace[4pt]
\textbf{Candidate:} 2610.06753\par\smallskip Discovering New Problems for Decoded Quantum Interferometry\par\smallskip mazumder2026discovering\par\smallskip \url{https://arxiv.org/abs/2610.06753} &
\textbf{First public:} 2026-10-06\par\smallskip \textbf{Version inspected:} arXiv:2610.06753v1\par\smallskip \textbf{Version date:} 2026-10-05\par\smallskip \textbf{Checked:} 2026-10-12 &
Public arXiv October/recent listing plus date-filtered topic searches; bibliography follow-up for older ADMM leads; discovery.json &
\textbf{Decision:} include context\par\smallskip \textbf{Placement:} Minimal contextual citation &
DQI screening rules and decoding conditions sharpen the boundary between quality guarantees and tested classical comparators. Mechanism/theory context, not a physical-hardware evidence row. Date basis: in window announcement. \\
\addlinespace[4pt]
\textbf{Candidate:} 2610.06797\par\smallskip From Quantum Tunneling to Classical Optimization through Instantons\par\smallskip \url{https://arxiv.org/abs/2610.06797} &
\textbf{First public:} 2026-10-06\par\smallskip \textbf{Version inspected:} arXiv:2610.06797v1\par\smallskip \textbf{Version date:} 2026-10-05\par\smallskip \textbf{Checked:} 2026-10-12 &
Public arXiv October/recent listing plus date-filtered topic searches; bibliography follow-up for older ADMM leads; discovery.json &
\textbf{Decision:} screen only\par\smallskip \textbf{Placement:} Discovery log only &
Classical path-integral stochastic gradient descent can inherit tunneling escape-rate scaling in continuous potentials; useful adjacent reassessment but no discrete-decomposition comparison. Abstract screened; does not rule out all tunneling-based quantum advantage. Date basis: in window announcement. \\
\addlinespace[4pt]
\textbf{Candidate:} 2610.06815\par\smallskip CV-QAOA: Efficient Low-Depth Quantum Optimization of Continuous Variables\par\smallskip \url{https://arxiv.org/abs/2610.06815} &
\textbf{First public:} 2026-10-06\par\smallskip \textbf{Version inspected:} arXiv:2610.06815v1\par\smallskip \textbf{Version date:} 2026-10-05\par\smallskip \textbf{Checked:} 2026-10-12 &
Public arXiv October/recent listing plus date-filtered topic searches; bibliography follow-up for older ADMM leads; discovery.json &
\textbf{Decision:} screen only\par\smallskip \textbf{Placement:} Discovery log only &
Continuous-variable QAOA query guarantees and dequantization on constructed continuous objectives; outside the discrete-decomposition evidence focus. Abstract screened. Date basis: in window announcement. \\
\addlinespace[4pt]
\textbf{Candidate:} 2610.07144\par\smallskip Decoupling Spectral Gaps from Entanglement Complexity in Quantum Annealing\par\smallskip \url{https://arxiv.org/abs/2610.07144} &
\textbf{First public:} 2026-10-07\par\smallskip \textbf{Version inspected:} arXiv:2610.07144v1\par\smallskip \textbf{Version date:} 2026-10-05\par\smallskip \textbf{Checked:} 2026-10-12 &
Public arXiv October/recent listing plus date-filtered topic searches; bibliography follow-up for older ADMM leads; discovery.json &
\textbf{Decision:} screen only\par\smallskip \textbf{Placement:} Discovery log only &
Spectral-gap versus entanglement complexity theory supports hardware-context cautions; no decomposition attribution benchmark. Abstract screened. Date basis: in window announcement. \\
\addlinespace[4pt]
\textbf{Candidate:} 2610.09584\par\smallskip Quantum Approximation Complexity of Classical Optimization Problems\par\smallskip hadfield2026approximation\par\smallskip \url{https://arxiv.org/abs/2610.09584} &
\textbf{First public:} 2026-10-08\par\smallskip \textbf{Version inspected:} arXiv:2610.09584v1\par\smallskip \textbf{Version date:} 2026-10-07\par\smallskip \textbf{Checked:} 2026-10-12 &
Public arXiv October/recent listing plus date-filtered topic searches; bibliography follow-up for older ADMM leads; discovery.json &
\textbf{Decision:} include context\par\smallskip \textbf{Placement:} Minimal contextual citation &
Formal classical-output approximation guarantees require feasibility, score, success and full algorithmic cost. Definitional/theoretical context; no empirical row. Date basis: in window announcement. \\
\addlinespace[4pt]
\textbf{Candidate:} 2610.09798\par\smallskip Quantum Optimization of the Color-of-Money Problem on IBM Quantum Hardware: A Proof-of-Concept Study\par\smallskip \url{https://arxiv.org/abs/2610.09798} &
\textbf{First public:} 2026-10-08\par\smallskip \textbf{Version inspected:} arXiv:2610.09798v2\par\smallskip \textbf{Version date:} 2026-10-08\par\smallskip \textbf{Checked:} 2026-10-12 &
Public arXiv October/recent listing plus date-filtered topic searches; bibliography follow-up for older ADMM leads; discovery.json &
\textbf{Decision:} screen only\par\smallskip \textbf{Placement:} Discovery log only &
Small MILP-to-PCE hardware demonstration reproduces an exact classical solution; no decomposition or comparable-resource advantage. Abstract screened. Date basis: in window announcement. \\
\addlinespace[4pt]
\textbf{Candidate:} 2610.10254\par\smallskip Benchmarking Modular Optimization Strategies for Parameterized Quantum Circuits\par\smallskip \url{https://arxiv.org/abs/2610.10254} &
\textbf{First public:} 2026-10-08\par\smallskip \textbf{Version inspected:} arXiv:2610.10254v1\par\smallskip \textbf{Version date:} 2026-10-07\par\smallskip \textbf{Checked:} 2026-10-12 &
Public arXiv October/recent listing plus date-filtered topic searches; bibliography follow-up for older ADMM leads; discovery.json &
\textbf{Decision:} screen only\par\smallskip \textbf{Placement:} Discovery log only &
Modular classical direction/update training controls are relevant methodology, but the paper's quantum circuit is the object being trained rather than a swapped component solving a fixed classical mathematical program. Abstract and selected full-text headings/setup checked. Date basis: in window announcement. \\
\addlinespace[4pt]
\textbf{Candidate:} 2610.10351\par\smallskip Measurement-Efficient Differentiable Quantum Architecture Search for Combinatorial Optimization\par\smallskip \url{https://arxiv.org/abs/2610.10351} &
\textbf{First public:} 2026-10-08\par\smallskip \textbf{Version inspected:} arXiv:2610.10351v1\par\smallskip \textbf{Version date:} 2026-10-07\par\smallskip \textbf{Checked:} 2026-10-12 &
Public arXiv October/recent listing plus date-filtered topic searches; bibliography follow-up for older ADMM leads; discovery.json &
\textbf{Decision:} screen only\par\smallskip \textbf{Placement:} Discovery log only &
DQAS gradient-measurement reduction changes resource demand within circuit search; no decomposition or system-level comparative evidence. Abstract and selected full-text setup checked. Date basis: in window announcement. \\
\addlinespace[4pt]
\textbf{Candidate:} 2610.10542\par\smallskip Omnisolver: An extensible interface to Ising spin-glass and QUBO solvers: adding a distributed GPU brute-force plugin\par\smallskip \url{https://arxiv.org/abs/2610.10542} &
\textbf{First public:} 2026-10-09\par\smallskip \textbf{Version inspected:} arXiv:2610.10542v1\par\smallskip \textbf{Version date:} 2026-06-01\par\smallskip \textbf{Checked:} 2026-10-12 &
Public arXiv October/recent listing plus date-filtered topic searches; bibliography follow-up for older ADMM leads; discovery.json &
\textbf{Decision:} screen only\par\smallskip \textbf{Placement:} Discovery log only &
Omnisolver GPU plugin is relevant software context; primary record says June 1 submission despite October arXiv identifier/listing. Public-release provenance is ambiguous and no new benchmark claim is extracted. Date basis: in window announcement. \\
\addlinespace[4pt]
\textbf{Candidate:} 2610.10772\par\smallskip Scaling Quantum Optimization to the Thousand-Qubit Scale with Distributed Quantum Sampling\par\smallskip liaqat2026scaling\par\smallskip \url{https://arxiv.org/abs/2610.10772} &
\textbf{First public:} 2026-10-09\par\smallskip \textbf{Version inspected:} arXiv:2610.10772v1\par\smallskip \textbf{Version date:} 2026-10-07\par\smallskip \textbf{Checked:} 2026-10-12 &
Public arXiv October/recent listing plus date-filtered topic searches; bibliography follow-up for older ADMM leads; discovery.json &
\textbf{Decision:} include evidence\par\smallskip \textbf{Placement:} E17; main evidence &
E17: partitioned QAOA/SA/random fixed-wrapper controls with physical IBM execution, repair and packing ablations; no matched-cost/system advantage. Date basis: in window announcement. \\
\addlinespace[4pt]
\textbf{Candidate:} 2610.11262\par\smallskip Gap Sensitivity and Bottleneck Migration in Sparse Constraint-Preserving Quantum Annealing\par\smallskip \url{https://arxiv.org/abs/2610.11262} &
\textbf{First public:} 2026-10-09\par\smallskip \textbf{Version inspected:} arXiv:2610.11262v1\par\smallskip \textbf{Version date:} 2026-10-08\par\smallskip \textbf{Checked:} 2026-10-12 &
Public arXiv October/recent listing plus date-filtered topic searches; bibliography follow-up for older ADMM leads; discovery.json &
\textbf{Decision:} screen only\par\smallskip \textbf{Placement:} Discovery log only &
Constraint-preserving driver adaptation uses gap sensitivity and finite-size oracle information; no quantum-versus-classical solver comparison. Abstract and methods/limitations inspected. Date basis: in window announcement. \\
\addlinespace[4pt]
\textbf{Candidate:} 2610.11710\par\smallskip Parameter-Dependent Noise Resilience in the Quantum Approximate Optimization Algorithm\par\smallskip \url{https://arxiv.org/abs/2610.11710} &
\textbf{First public:} 2026-10-09\par\smallskip \textbf{Version inspected:} arXiv:2610.11710v1\par\smallskip \textbf{Version date:} 2026-10-08\par\smallskip \textbf{Checked:} 2026-10-12 &
Public arXiv October/recent listing plus date-filtered topic searches; bibliography follow-up for older ADMM leads; discovery.json &
\textbf{Decision:} screen only\par\smallskip \textbf{Placement:} Discovery log only &
Parameter-dependent QAOA noise bounds and training modifications; useful quantum-method context but not decomposition or attribution evidence. Abstract screened. Date basis: in window announcement. \\
\addlinespace[4pt]
\textbf{Candidate:} 2610.11759\par\smallskip Beyond QAOA: A Review of AI and Quantum Computing for Adaptive Combinatorial Optimization\par\smallskip lau2026beyond\par\smallskip \url{https://arxiv.org/abs/2610.11759} &
\textbf{First public:} 2026-10-09\par\smallskip \textbf{Version inspected:} arXiv:2610.11759v2\par\smallskip \textbf{Version date:} 2026-10-09\par\smallskip \textbf{Checked:} 2026-10-12 &
Public arXiv October/recent listing plus date-filtered topic searches; bibliography follow-up for older ADMM leads; discovery.json &
\textbf{Decision:} include context\par\smallskip \textbf{Placement:} Minimal contextual citation &
New review of AI-controlled optimization workflows provides related-survey context and independently distinguishes system and quantum-component claims. Its corpus counts are not pooled with this survey. Date basis: in window announcement. \\
\addlinespace[4pt]
\textbf{Candidate:} 2610.13283\par\smallskip Circuit Design Problems in Binary Polynomial Quantum Optimization\par\smallskip \url{https://arxiv.org/abs/2610.13283} &
\textbf{First public:} 2026-10-12\par\smallskip \textbf{Version inspected:} arXiv:2610.13283v1\par\smallskip \textbf{Version date:} 2026-10-09\par\smallskip \textbf{Checked:} 2026-10-12 &
Public arXiv October/recent listing plus date-filtered topic searches; bibliography follow-up for older ADMM leads; discovery.json &
\textbf{Decision:} screen only\par\smallskip \textbf{Placement:} Discovery log only &
Classical MIP/heuristic circuit synthesis for higher-order binary objectives; quantum compilation rather than a hybrid optimization evidence comparison. Abstract screened. Date basis: in window announcement. \\
\addlinespace[4pt]
\textbf{Candidate:} 2609.29856\par\smallskip A Block Decomposed QUBO Workflow for Chromosome-Y Phylogeny Reconstruction\par\smallskip \url{https://arxiv.org/abs/2609.29856} &
\textbf{First public:} Not captured by this audit.\par\smallskip \textbf{Version inspected:} arXiv:2609.29856v1\par\smallskip \textbf{Version date:} 2026-09-24\par\smallskip \textbf{Checked:} 2026-10-12 &
Public arXiv October/recent listing plus date-filtered topic searches; bibliography follow-up for older ADMM leads; discovery.json &
\textbf{Decision:} screen only\par\smallskip \textbf{Placement:} Discovery log only &
September 24 v1 only; no in-window revision or publication found on the primary record. Relevant QUBO decomposition lead remains outside this bounded extension. Date basis: no in window event verified. Not independently verified; submission 2026-09-24. \\
\addlinespace[4pt]
\textbf{Candidate:} 2608.04159\par\smallskip A Distributed Quantum Approximate Optimization Algorithm For Unit Commitment\par\smallskip \url{https://arxiv.org/abs/2608.04159} &
\textbf{First public:} Not captured by this audit.\par\smallskip \textbf{Version inspected:} arXiv:2608.04159v1\par\smallskip \textbf{Version date:} 2026-08-04\par\smallskip \textbf{Checked:} 2026-10-12 &
Public arXiv October/recent listing plus date-filtered topic searches; bibliography follow-up for older ADMM leads; discovery.json &
\textbf{Decision:} screen only\par\smallskip \textbf{Placement:} Discovery log only &
August 4 v1 only; no in-window revision/publication on the inspected record. Relevant distributed QAOA/ADMM mechanism is an older lead, not a new October release. Date basis: no in window event verified. Not independently verified; submission 2026-08-04. \\
\addlinespace[4pt]
\textbf{Candidate:} 2603.14744\par\smallskip Towards Exponential Quantum Improvements in Solving Cardinality-Constrained Binary Optimization\par\smallskip \url{https://arxiv.org/abs/2603.14744} &
\textbf{First public:} Not captured by this audit.\par\smallskip \textbf{Version inspected:} arXiv:2603.14744v1\par\smallskip \textbf{Version date:} 2026-03-16\par\smallskip \textbf{Checked:} 2026-10-12 &
Public arXiv October/recent listing plus date-filtered topic searches; bibliography follow-up for older ADMM leads; discovery.json &
\textbf{Decision:} screen only\par\smallskip \textbf{Placement:} Discovery log only &
March 16 v1 only; no in-window revision/publication on inspected record. Older quantum ADMM/cardinality theory lead. Date basis: no in window event verified. Not independently verified; submission 2026-03-16. \\
\addlinespace[4pt]
\textbf{Candidate:} 2610.13621\par\smallskip Encoding Matters: A Qubit-Efficient Reformulation for Quantum Resilient Routing\par\smallskip \url{https://arxiv.org/abs/2610.13621} &
\textbf{First public:} 2026-10-12\par\smallskip \textbf{Version inspected:} arXiv:2610.13621v1\par\smallskip \textbf{Version date:} 2026-10-09\par\smallskip \textbf{Checked:} 2026-10-12 &
Public arXiv October/recent listing plus date-filtered topic searches; bibliography follow-up for older ADMM leads; discovery.json &
\textbf{Decision:} screen only\par\smallskip \textbf{Placement:} Discovery log only &
Main text inspected. Redundant classical-presolve example under the page budget; optimum-aware parameter training, different circuit depths and a published baseline limit attribution. No quantum component swap. The full author surname is Dalla Torre. Date basis: in window announcement. \\
\addlinespace[4pt]
\textbf{Candidate:} 2610.07252\par\smallskip Resource-Aware Grover Search for Minimum Vertex Cover\par\smallskip \url{https://arxiv.org/abs/2610.07252} &
\textbf{First public:} 2026-10-07\par\smallskip \textbf{Version inspected:} arXiv:2610.07252v1\par\smallskip \textbf{Version date:} 2026-10-05\par\smallskip \textbf{Checked:} 2026-10-12 &
Public arXiv October/recent listing plus date-filtered topic searches; bibliography follow-up for older ADMM leads; discovery.json &
\textbf{Decision:} screen only\par\smallskip \textbf{Placement:} Discovery log only &
Grover oracle-resource accounting for minimum vertex cover; not a hybrid decomposition benchmark. Abstract screened. Date basis: in window announcement. \\
\addlinespace[4pt]
\textbf{Candidate:} 2610.05128\par\smallskip Training Variational Quantum Algorithms Is NP-Hard, Even Locally\par\smallskip \url{https://arxiv.org/abs/2610.05128} &
\textbf{First public:} 2026-10-06\par\smallskip \textbf{Version inspected:} arXiv:2610.05128v1\par\smallskip \textbf{Version date:} 2026-10-04\par\smallskip \textbf{Checked:} 2026-10-12 &
Public arXiv October/recent listing plus date-filtered topic searches; bibliography follow-up for older ADMM leads; discovery.json &
\textbf{Decision:} screen only\par\smallskip \textbf{Placement:} Discovery log only &
VQA-training NP-hardness addresses the classical parameter-selection layer already qualified in the manuscript; context-only abstract screen, not an empirical optimization result. Date basis: in window announcement. \\
\addlinespace[4pt]
\textbf{Candidate:} 2610.06638\par\smallskip Hybrid Classical-Quantum Solutions to Accelerate the Adoption of Quantum Computing\par\smallskip \url{https://arxiv.org/abs/2610.06638} &
\textbf{First public:} Not captured by this audit.\par\smallskip \textbf{Version inspected:} arXiv:2610.06638v1\par\smallskip \textbf{Version date:} 2026-10-05\par\smallskip \textbf{Checked:} 2026-10-12 &
Public arXiv October/recent listing plus date-filtered topic searches; bibliography follow-up for older ADMM leads; discovery.json &
\textbf{Decision:} screen only\par\smallskip \textbf{Placement:} Discovery log only &
Broad hybrid software-engineering review; no specific new decomposition implementation or controlled result extracted. Abstract screened. Date basis: in window version submission and observed public record. Not independently verified; submission 2026-10-05. \\
\addlinespace[4pt]
\textbf{Candidate:} doi:10.1007/s10791-026-10589-w\par\smallskip A convex-Feynman hybrid optimization paradigm for energy-efficient scheduling in large-scale multi-tier cloud data centers\par\smallskip \url{https://link.springer.com/article/10.1007/s10791-026-10589-w} &
\textbf{First public:} 2026-10-09\par\smallskip \textbf{Version inspected:} Publisher version of record\par\smallskip \textbf{Version date:} 2026-10-09\par\smallskip \textbf{Checked:} 2026-10-12 &
Publisher-targeted October quantum optimization query; discovery.json &
\textbf{Decision:} exclude scope\par\smallskip \textbf{Placement:} Discovery log only &
Despite quantum/Feynman terminology, Sections 7--9 specify classical Metropolis sampling, convex optimization/ADMM and conventional server experiments; no quantum processor or emulated quantum optimization component identified. Date basis: in window journal publication. \\
\addlinespace[4pt]
\end{longtable}
\endgroup


\clearpage
%
\begingroup
\fontsize{7.5}{9.0}\selectfont
\setlength{\tabcolsep}{3.5pt}
\setlength{\AuditTableWidth}{\dimexpr\linewidth-21pt\relax}
\begin{longtable}{@{}>{\raggedright\arraybackslash}p{0.12\AuditTableWidth}>{\raggedright\arraybackslash}p{0.23\AuditTableWidth}>{\raggedright\arraybackslash}p{0.37\AuditTableWidth}>{\raggedright\arraybackslash}p{0.28\AuditTableWidth}@{}}
\caption{Source-location index for the 30 audited studies, including all 25 main evidence-table entries. Pointers identify the inspected version and relevant passage; retained prior inspection and selective new checks are distinguished. Access limits do not imply a fresh verification of every numerical field or reproduction of experiments. For satisfiability instances, \(n\) counts variables and \(m\) counts clauses.}\label{tab:s6}\\
\toprule
\textbf{Study} & \textbf{Claim or comparison to locate} & \textbf{Inspected source and passage} & \textbf{Verification scope} \\
\midrule\endfirsthead
\multicolumn{4}{l}{\normalsize Table~\thetable\ continued.}\\
\toprule
\textbf{Study} & \textbf{Claim or comparison to locate} & \textbf{Inspected source and passage} & \textbf{Verification scope} \\
\midrule\endhead
\multicolumn{4}{r}{Continues on the next page.}\\\endfoot
\bottomrule\endlastfoot
[E1] Ponce et al., 2025 &
Reduced-arm mean approximation ratios near 0.96 versus monolithic QAOA near 0.75; ten 100-vertex hardware demonstrations. &
\textbf{Version:} Journal version, published 2025-02-18\par\smallskip \textbf{Location:} Abstract; Appendix B, Table 1 (approximation ratios by solver).\par\smallskip \url{https://link.springer.com/article/10.1007/s11128-025-04675-z} &
\textbf{Checked:} 2026-10-06\par\smallskip \textbf{Basis:} selective primary passage check\par\smallskip \textbf{Limits:} Publisher abstract and Table 1 caption accessible; main article is subscription content. Detailed 25-instance and 12-16-qubit setup retained from prior audit, not freshly reverified. \\
\addlinespace[4pt]
[E2] Schuetz et al. (2026) &
On 39 instances, SaReduMIS matches or exceeds the quantum arm using 500 replicas and two SA steps versus 1,000 measurements; quantum samples beat random within the shared wrapper. &
\textbf{Version:} arXiv v2 (2026-07-01), with central SaReduMIS claims selectively rechecked against journal version of record (2026-09-10) on 2026-10-04\par\smallskip \textbf{Location:} Journal PDF pp. 6 and 8 (SaReduMIS comparator, success, and sampling counts).\par\smallskip \url{https://journals.aps.org/prresearch/pdf/10.1103/nh3b-1wv5} &
\textbf{Checked:} 2026-10-04\par\smallskip \textbf{Basis:} retained prior inspection\par\smallskip \textbf{Limits:} Retained selective journal check in status-checks/group-a.json; the historical arXiv v2 audit is preserved. Nominal counts are not matched time. Remaining numerical fields and artifacts were not fully re-audited. \\
\addlinespace[4pt]
[E3] Dupont et al., 2025 &
Mean approximation ratio near 99\%; projected 100-qubit runtime trails stronger classical heuristics. &
\textbf{Version:} Physical Review Applied 23, 014045 (2025)\par\smallskip \textbf{Location:} Abstract: approximation-ratio and projected-runtime comparison sentences.\par\smallskip \url{https://journals.aps.org/prapplied/abstract/10.1103/PhysRevApplied.23.014045} &
\textbf{Checked:} 2026-10-06\par\smallskip \textbf{Basis:} selective primary passage check\par\smallskip \textbf{Limits:} Selective check of the cited passages, not a reproduction of the experiments. \\
\addlinespace[4pt]
[E4] Kim et al., 2025 &
Over 50 large dense instances; accuracy normalized by best observed objective; QPU programming/sampling timing excludes reading, queueing and local communication. &
\textbf{Version:} npj Quantum Information 11, 77 (2025), journal PDF\par\smallskip \textbf{Location:} Abstract (p. 1); Performance metrics: relative accuracy and computational time (pp. 5-6), Eq. (5).\par\smallskip \url{https://www.nature.com/articles/s41534-025-01020-1.pdf} &
\textbf{Checked:} 2026-10-06\par\smallskip \textbf{Basis:} selective primary passage check\par\smallskip \textbf{Limits:} Selective check of the cited passages, not a reproduction of the experiments. \\
\addlinespace[4pt]
[E5] Tasseff et al., 2024 &
Restricted runtime benefit on constructed sparse instance classes; no universal performance claim. &
\textbf{Version:} Journal of Heuristics (2024), version of record\par\smallskip \textbf{Location:} Sections 3 and 4; Appendix D, Figures 7-9.\par\smallskip \url{https://link.springer.com/article/10.1007/s10732-024-09530-5} &
\textbf{Checked:} 2026-10-06\par\smallskip \textbf{Basis:} selective primary passage check\par\smallskip \textbf{Limits:} Qualitative comparison, with no numerical speedup multiplier asserted in the manuscript. \\
\addlinespace[4pt]
[E6] Vodeb et al., 2026 &
139 maximum-cut instances with 100-10,000 vertices; hybrid and simulated annealing attain certified optima where available, through 251 graph vertices; larger-instance quality and heterogeneous timing require separate interpretations. &
\textbf{Version:} journal version (July 3, 2026)\par\smallskip \textbf{Location:} Sections 4.1, 4.2 and 5; Section 4.1 explains the extra graph vertex and Section 4.2 the bqp250 certification.\par\smallskip \url{https://link.springer.com/article/10.1007/s11128-026-05263-5} &
\textbf{Checked:} 2026-10-04\par\smallskip \textbf{Basis:} retained prior inspection\par\smallskip \textbf{Limits:} Retained selective check in status-checks/group-a.json; not a fresh reread or full re-audit of every instance. Hardware and some literature-sourced timings differ. \\
\addlinespace[4pt]
[E7] Lozano, 2026 &
Service matches 54 certified optima; mean QPU access is 0.034 seconds within a five-second nominal budget. Version 3 compares Tabu on 36 instances and 108 seeded five-second runs, with mean absolute gap 0.00078. &
\textbf{Version:} arXiv v3 (2026-07-12); post-cutoff expanded Tabu control selectively rechecked on 2026-10-04\par\smallskip \textbf{Location:} Version 3 abstract; matched-budget Tabu subsection and Table 8, checked against the version 2 abstract.\par\smallskip \url{https://arxiv.org/html/2605.17623v3} &
\textbf{Checked:} 2026-10-04\par\smallskip \textbf{Basis:} retained prior inspection\par\smallskip \textbf{Limits:} Retained selective check in status-checks/group-a.json. Expanded control is from the July 12 revision; service internals and all artifact fields were not re-audited. It is an external control, not an internal QPU ablation. \\
\addlinespace[4pt]
[E8] Zhu et al., 2026 &
20 hardest of 4000 instances; 104 qubits; 2.23-fold sampling gain and projected 2.34-fold TTS at 0.4 MHz. &
\textbf{Version:} National Science Review 13(9), nwag124 (2026)\par\smallskip \textbf{Location:} Benchmark-selection paragraph before Fig. 1; Fig. 3 and accompanying normalization discussion; Fig. 4 and TTS discussion; Supplementary Section S4D.\par\smallskip \url{https://academic.oup.com/nsr/article/13/9/nwag124/8503488} &
\textbf{Checked:} 2026-10-06\par\smallskip \textbf{Basis:} selective primary passage check\par\smallskip \textbf{Limits:} Measured samples and projected hardware throughput are distinct; GPU and planar-specific baseline comparisons are absent. \\
\addlinespace[4pt]
[E9] Pennington et al. (2025) &
Hardware QAOA outperforms annealing on 50- and 70-node graphs but not 100-node graphs; nine heavy-hex graphs require 52-111 qubits. &
\textbf{Version:} arXiv:2509.10657v2\par\smallskip \textbf{Location:} Section V.3.3, Fig. 10 and Results/Discussion paragraphs; Sections V.3.1-V.3.2 for small complete/bipartite graphs.\par\smallskip \url{https://arxiv.org/html/2509.10657v2} &
\textbf{Checked:} 2026-10-06\par\smallskip \textbf{Basis:} selective primary passage check\par\smallskip \textbf{Limits:} Selective check of the cited passages, not a reproduction of the experiments. \\
\addlinespace[4pt]
[E10] Dupont et al. (2025) &
200-instance averages; original versus preconditioned convergence; preconditioning accounting and limited Ankaa-3 demonstration. &
\textbf{Version:} Author preprint arXiv:2502.18570v2, supplementing the recorded October 3 journal version\par\smallskip \textbf{Location:} Sections III.1.2-III.1.4 and IV; Figs. 1-3, 5 and their captions.\par\smallskip \url{https://arxiv.org/html/2502.18570v2} &
\textbf{Checked:} 2026-10-06\par\smallskip \textbf{Basis:} selective primary passage check\par\smallskip \textbf{Limits:} Publisher abstract checked for identity; publisher PDF fetch failed. These locator labels refer explicitly to author v2, not an assumed journal pagination match. \\
\addlinespace[4pt]
[E11] Wagner and Liers (2026) &
Master-relaxation pricing/separation experiments are motivated by a complete branch-price-and-cut framework that is not implemented. Specialized exact pricing and separation run faster than the tested quantum alternatives. &
\textbf{Version:} arXiv:2412.15665v2 (April 24, 2026)\par\smallskip \textbf{Location:} Section 3 opening (implementation scope); Section 5 Pricing, Figures 4-7 (Figure 6 method times, Figure 7 relaxation time); Section 5 Separation, Figures 8-10; Table 1.\par\smallskip \url{https://arxiv.org/html/2412.15665v2} &
\textbf{Checked:} 2026-10-06\par\smallskip \textbf{Basis:} selective primary passage check\par\smallskip \textbf{Limits:} Selective scope/outcome reread. Faster specialized classical methods does not mean every classical comparator wins every metric: QA often beats SA pricing time and improves on CVRPSEP gaps for three cases. No complete integer-solver comparison. \\
\addlinespace[4pt]
[E12] Perron et al. (2026) &
Up to 300 vertices; gaps within 5-10\%; equal nominal sampling budgets and simulated annealing strongest. &
\textbf{Version:} arXiv:2603.28933v1\par\smallskip \textbf{Location:} Abstract; Section III Results; Fig. 3 caption and accompanying discussion.\par\smallskip \url{https://arxiv.org/html/2603.28933v1} &
\textbf{Checked:} 2026-10-06\par\smallskip \textbf{Basis:} selective primary passage check\par\smallskip \textbf{Limits:} Selective check of the cited passages, not a reproduction of the experiments. \\
\addlinespace[4pt]
[E13] de Valk et al., 2026 &
Ten permutations of one 20-customer instance and five 5-customer toys; MPS-QAOA cut selection consumes 99.3\% of toy-pipeline wall time, with IBM component execution demonstrated separately. &
\textbf{Version:} arXiv:2607.26771v1\par\smallskip \textbf{Location:} Sections III.B-III.C; Section IV, Tables I-II; Sections V.A-V.C, Figures 2-5; Section VI.\par\smallskip \url{https://arxiv.org/html/2607.26771v1} &
\textbf{Checked:} 2026-10-04\par\smallskip \textbf{Basis:} retained prior inspection\par\smallskip \textbf{Limits:} Retained claim-critical inspection in discovery-search/extension/decomposition-search.json. Ten permutations are not ten independent routing instances; MPS wall time is classical emulation time. This location record excludes the separately corrected author surname. \\
\addlinespace[4pt]
[E14] Naghmouchi et al., 2026 &
Two routing topologies up to 30 nodes; initialization/refinement ablations bring the quantum-emulated arm mean gap below 1\% to CPLEX, with 95\% intervals. Warm-start alone narrows backend differences. &
\textbf{Version:} arXiv:2609.37900v1\par\smallskip \textbf{Location:} Sections V.1-V.4 and V.5.2 (pricing, reductions, recovery); Section VI.1 (benchmarks and comparator); Sections VI.2-VI.3 and Figures 4-7 (mean gaps and ablations).\par\smallskip \url{https://arxiv.org/html/2609.37900v1} &
\textbf{Checked:} 2026-10-04\par\smallskip \textbf{Basis:} retained prior inspection\par\smallskip \textbf{Limits:} Retained inspection in discovery-search/extension/decomposition-search.json. Below 1\% refers to plotted means, not every instance. No hardware execution, wall-clock speedup, or exact-pricing certificate is established. \\
\addlinespace[4pt]
[E15] Ramesh et al., 2026 &
Fifty balanced-bipartition instances per size, n=12-40, use analytic or emulated QAOA correlations. Earlier near-optimal incumbent discovery excludes correlation cost; main penalties use oracle information and alternative penalty choices are also tested. &
\textbf{Version:} arXiv:2608.28842v1\par\smallskip \textbf{Location:} Sections 3.2-3.4 (objective and hard balance), 3.5 and 3.7 (scale and cost exclusions), 4.3/Figure 5 and Appendix B (penalty choices), 4.4-5 (optimum limitations).\par\smallskip \url{https://arxiv.org/html/2608.28842v1} &
\textbf{Checked:} 2026-10-04\par\smallskip \textbf{Basis:} retained prior inspection\par\smallskip \textbf{Limits:} Retained inspection in discovery-search/extension/constraints-theory-search.json. No hardware execution or whole-system quantum speedup; optimizing the surrogate does not certify original optimality. \\
\addlinespace[4pt]
[E16] Yalovetzky et al., 2026 &
Four hardware market graphs up to 225 assets and 78-qubit first kernels, plus 73 emulator graphs, compare iterative qReduMIS with standalone QAOA on shared first kernels. Shot-TTS excludes classical and communication time. &
\textbf{Version:} arXiv:2607.01037v2\par\smallskip \textbf{Location:} Sections IV-V; Appendix E.1-E.2 (source anchors \#S4, \#S5, \#A5.SS1 and \#A5.SS2).\par\smallskip \url{https://arxiv.org/html/2607.01037v2} &
\textbf{Checked:} 2026-10-04\par\smallskip \textbf{Basis:} retained prior inspection\par\smallskip \textbf{Limits:} Retained selective version 2 check in discovery-search/extension/comparative-search.json. Raw-success budgets differ; FTSE shot-TTS is not lower (185.58 versus 151.19, overlapping intervals). No hardware classical-sampler substitution or complete system timing. \\
\addlinespace[4pt]
[E17] Liaqat et al., 2026 &
Repaired QAOA candidates lower residual error relative to the tested SA/random samplers in a shared partitioned E-FCFW wrapper; packing reduces reported hardware execution time, without a matched-cost or system advantage claim. &
\textbf{Version:} arXiv:2610.10772v1 (2026-10-07)\par\smallskip \textbf{Location:} Sections III.2--III.3, IV.1--IV.4; Figures 2--3; Table 1; Appendix A.1--A.4 and Algorithm 1.\par\smallskip \url{https://arxiv.org/html/2610.10772v1} &
\textbf{Checked:} 2026-10-12\par\smallskip \textbf{Basis:} new study inspection\par\smallskip \textbf{Limits:} Specified unequal sampling budgets; packed17 vs unpacked20 rounds; no comprehensive runtime or repeated-instance uncertainty; no benchmark-code reproduction. \\
\addlinespace[4pt]
[E18] Lauwerier et al., 2026 &
On one fixed tour-selection instance at 10, 15, 20 vertices, the emulated quantum selector matches the exact restricted optimum only at 15 and is below it at 10 and 20; all arms return feasible selections. &
\textbf{Version:} arXiv:2610.03890v1 (submitted 2026-10-02; announced 2026-10-06)\par\smallskip \textbf{Location:} Sections 3.3--3.4; 4.1; 5.5.1--5.5.3; Table 3 and its runtime footnote.\par\smallskip \url{https://arxiv.org/html/2610.03890v1} &
\textbf{Checked:} 2026-10-12\par\smallskip \textbf{Basis:} new study inspection\par\smallskip \textbf{Limits:} One fixed-date controlled instance; noiseless quantum emulation; QPU time projected from 100 shots at 1 Hz; LP bound applies to restricted pool and may be loose; no end-to-end hardware claim. \\
\addlinespace[4pt]
[A1] Munoz-Bauza and Lidar, 2025 &
Tailored high-precision 2D spin-glass family; over 1,300 error-suppressed logical qubits; time-to-epsilon scaling comparison against parallel tempering with isoenergetic cluster moves; favorable reported scaling at optimality gaps of at least 1.0\%. &
\textbf{Version:} Published journal record, Physical Review Letters 134, 160601 (published April 23, 2025), DOI 10.1103/PhysRevLett.134.160601\par\smallskip \textbf{Location:} Abstract: benchmark/comparator passage, QAC logical-qubit count passage, and final optimality-gap sentence.\par\smallskip \url{https://journals.aps.org/prl/abstract/10.1103/PhysRevLett.134.160601} &
\textbf{Checked:} 2026-10-06\par\smallskip \textbf{Basis:} selective primary passage check\par\smallskip \textbf{Limits:} The numerical thresholds and comparator are explicit in the published abstract. This selective source check does not independently reproduce the scaling fits or establish end-to-end advantage. \\
\addlinespace[4pt]
[A2] Pawlowski et al., 2026 &
The same Sidon-28 instance family and time-to-epsilon target are reassessed with GPU simulated bifurcation. Total solver runtime and pure GPU runtime are distinguished; comparison closes the reported quantum-classical scaling gap. Small available quantum sizes and overhead sensitivity limit asymptotic inference. &
\textbf{Version:} Author arXiv:2505.22514v2 (July 10, 2026), accompanying the journal article Physical Review Applied 26, 014024 (July 8, 2026), DOI 10.1103/qcws-8kgm\par\smallskip \textbf{Location:} Figures 1-2 and their captions; Thermal vs Chaotic vs Quantum Annealing section, Eq. (1) and following runtime-definition paragraphs; Supplemental Material S2 and S5.\par\smallskip \url{https://arxiv.org/html/2505.22514v2} &
\textbf{Checked:} 2026-10-06\par\smallskip \textbf{Basis:} selective primary passage check\par\smallskip \textbf{Limits:} Journal abstract independently checked at \url{https://journals.aps.org/prapplied/abstract/10.1103/qcws-8kgm}; journal body requires authorization, so the detailed figure/section locators refer specifically to author v2. No performance reproduction was undertaken. \\
\addlinespace[4pt]
[A3] Chandarana et al., 2026 &
Twenty higher-order instances, with target energy matching the independently reported ground-state energy for 14; solver-side accounting includes preprocessing, QPU execution interval and postprocessing, while deployment/launch/communication components are treated separately. Figure 1 and its discussion compare matched-time quality against CPU baselines, PT+ and GPU ABS3; 128 vCPUs and eight A100 GPUs are reported. &
\textbf{Version:} arXiv:2603.13607v1, March 13, 2026\par\smallskip \textbf{Location:} Abstract; Benchmark design and metrics section (runtime-accounting and target-energy paragraphs); Figure 1 caption and Sub-second latency regime; Time-to-solution distributions, Table 1.\par\smallskip \url{https://arxiv.org/html/2603.13607v1} &
\textbf{Checked:} 2026-10-06\par\smallskip \textbf{Basis:} selective primary passage check\par\smallskip \textbf{Limits:} The time-to-solution metric in Table 1 is distinct from Figure 1 matched-time quality; PT+ is strongest on the reported TTS comparison. The audit preserves this distinction and does not equate the authors' solver-side budget with full deployment latency. \\
\addlinespace[4pt]
[A4] Ebadi et al., 2022 &
Native unit-disk MIS experiments range from 39 to 289 atoms. The positive comparison against optimized simulated annealing is restricted to hard instances in the deep-circuit regime; the text explicitly says the speedup is absent when the depth condition is not met. &
\textbf{Version:} Author arXiv:2202.09372v1 (February 18, 2022), corresponding to Science 376, 1209-1215 (2022), DOI 10.1126/science.abo6587\par\smallskip \textbf{Location:} Introduction paragraph beginning Here we address this challenge (39-289 qubits); Figure 4E caption; Quantum speedup on the hardest graphs discussion before Discussion and Outlook.\par\smallskip \url{https://arxiv.org/html/2202.09372v1} &
\textbf{Checked:} 2026-10-06\par\smallskip \textbf{Basis:} selective primary passage check\par\smallskip \textbf{Limits:} Locators refer to the inspected author v1, not the final journal typesetting. The result is a restricted hardware-evolution/success comparison, not a complete application-pipeline comparison. \\
\addlinespace[4pt]
[A5] Andrist et al., 2023 &
Structured quasiplanar instances with up to thousands of nodes are solved to optimality within minutes on commodity hardware; relaxing restrictions on the classical simulated-annealing method produces competitive scaling on the studied family. &
\textbf{Version:} Published journal record, Physical Review Research 5, 043277 (December 21, 2023), DOI 10.1103/PhysRevResearch.5.043277\par\smallskip \textbf{Location:} Abstract, sentences describing Union-Jack-like quasiplanar instances and the subsequent simulated-annealing scaling comparison.\par\smallskip \url{https://journals.aps.org/prresearch/abstract/10.1103/PhysRevResearch.5.043277} &
\textbf{Checked:} 2026-10-06\par\smallskip \textbf{Basis:} selective primary passage check\par\smallskip \textbf{Limits:} The size/time description and comparator qualification are explicit in the published abstract. They do not extend to arbitrary unit-disk instances; the same abstract reports much harder families with less structure or higher connectivity. \\
\addlinespace[4pt]
[A6] Cazals et al., 2026 &
Hardness is varied through density and treewidth; quantum exact-solution times are slower in the small-instance comparison. Extended classical benchmarks motivate approximately 1,000 atoms at a 1-kHz repetition rate as a necessary step toward possible advantage. &
\textbf{Version:} Published journal record, Physical Review Applied 25, 034085 (March 26, 2026), DOI 10.1103/8gjv-ll1f\par\smallskip \textbf{Location:} Abstract: density/treewidth construction passage, small-instance comparison sentence, and final 1,000-atom/1-kHz estimate.\par\smallskip \url{https://journals.aps.org/prapplied/abstract/10.1103/8gjv-ll1f} &
\textbf{Checked:} 2026-10-06\par\smallskip \textbf{Basis:} selective primary passage check\par\smallskip \textbf{Limits:} The 1,000-atom/1-kHz value is a projection from classical benchmarks, not an observed quantum advantage or a sufficient guarantee. Journal body requires authorization; the exact numerical projection and qualitative conclusions checked here are explicit in the public abstract. \\
\addlinespace[4pt]
[A7] Lu et al., 2026 &
At m/n=0.626, fits through n=70 give bases 1.0077\textasciicircum{}n and 1.0070\textasciicircum{}n versus 1.0128\textasciicircum{}n, using inverse quantum success probability versus classical conflict/propagation counts. A 13-qubit experiment reaches 22-variable instances. &
\textbf{Version:} Nature Computational Science journal version (June 19, 2026)\par\smallskip \textbf{Location:} Main item (4); Results: Scalings of time complexities, Table 1; Figure 1 (co-design changes), Figure 2c-e (simulation scale); Results: Experiment on superconducting quantum processors, second series, Figure 5.\par\smallskip \url{https://www.nature.com/articles/s43588-026-01007-8} &
\textbf{Checked:} 2026-10-06\par\smallskip \textbf{Basis:} selective primary passage check\par\smallskip \textbf{Limits:} Selective claim check. The experiment consistently improves on its reported baseline only through n=20; n=22 is the reached size. Cost units and simultaneous construction changes prevent a common-resource interaction estimate. \\
\addlinespace[4pt]
Context He et al., 2026 &
Five 96-node graphs on Helios compare warm-started QAOA with implemented GW and HLZ; local search improves quality and success. The unrestricted Burer-Monteiro runtime crossover is projected. &
\textbf{Version:} arXiv v1 (2026-03-10); hardware-comparator characterization selectively rechecked on 2026-10-04\par\smallskip \textbf{Location:} Figure 2; Sections II.3, II.5 and IV.5 (hardware quality, projected crossover and classical implementations).\par\smallskip \url{https://arxiv.org/html/2603.10191v1} &
\textbf{Checked:} 2026-10-04\par\smallskip \textbf{Basis:} retained prior inspection\par\smallskip \textbf{Limits:} Retained selective check in status-checks/group-a.json. Not every simulation/resource assumption was re-audited; current hardware quality comparisons do not establish the projected fault-tolerant speedup. \\
\addlinespace[4pt]
Context Tao et al., 2026 &
GPU quantum-annealing-inspired samplers provide a classical counterpoint on the related multi-objective benchmarks. Sampling is compared with published quantum results; full model-construction, sampling and filtering time is compared with classical baselines. &
\textbf{Version:} arXiv:2604.26477v2 (May 20, 2026)\par\smallskip \textbf{Location:} Results and discussion: Sampling time advantage, Figure 1/Table 1; Sample quality and end-to-end solution time, Figure 2 and accompanying text; Comparison to evolutionary algorithm-based heuristics, Table 2.\par\smallskip \url{https://arxiv.org/html/2604.26477v2} &
\textbf{Checked:} 2026-10-06\par\smallskip \textbf{Basis:} selective primary passage check\par\smallskip \textbf{Limits:} Bounded qualitative reread, not a full empirical audit. Quantum sampling comparisons reuse reported results; Figure 2 complete-pipeline comparisons are against classical DCM/DPA-a. The paper allows that QA may show better end-to-end performance on moderate hardware-native graphs. \\
\addlinespace[4pt]
Context Safi et al., 2026 &
Industrial scheduling separates assignment QUBOs and classical sequencing. Larger-scale results come from the classical Digital Annealer; quantum devices address smaller subsets and reported time excludes cloud integration costs. &
\textbf{Version:} arXiv:2607.13325v1\par\smallskip \textbf{Location:} Section III.C and Sections V-VIII; primary HTML anchors \#S5, \#S6 and \#S7 for formulations, platforms and results.\par\smallskip \url{https://arxiv.org/html/2607.13325v1} &
\textbf{Checked:} 2026-10-04\par\smallskip \textbf{Basis:} retained prior inspection\par\smallskip \textbf{Limits:} Retained inspection in discovery-search/extension/comparative-search.json. Dataset scale (2-32 machines, 2-85 rigs, 2-256 jobs) is not full quantum-execution scale; model and platform effects overlap; communication, API and queueing are omitted. \\
\addlinespace[4pt]
Context Bantysh et al., 2026 &
BOND-1 supplies a classical product-state comparator in a fixed-parameter expanding-depth regime. GSet cuts use best-of-ten reporting and a trained schedule; preprocessing is excluded from timing. &
\textbf{Version:} arXiv:2609.15738v1\par\smallskip \textbf{Location:} Sections III, VI and VII (primary HTML anchors \#S3, \#S6 and \#S7): entanglement regime, GSet protocol/results and scope discussion.\par\smallskip \url{https://arxiv.org/html/2609.15738v1} &
\textbf{Checked:} 2026-10-04\par\smallskip \textbf{Basis:} retained prior inspection\par\smallskip \textbf{Limits:} Retained inspection in discovery-search/extension/comparative-search.json. Schedule training includes G1, G23 and G81, so those cases are not held out. Best cuts and per-run mean times are different summaries; this does not show entanglement is generally unnecessary. \\
\addlinespace[4pt]
Context Gaidai and Mukherjee, 2026 &
Across 5-13-generator instances, ideal variational sampling fails to improve incumbent quality over uniform sampling with the same continuous solver, even with quantum simulation time excluded from the charged budget. &
\textbf{Version:} arXiv:2607.15543v1\par\smallskip \textbf{Location:} Sections III.4, IV.2 and IV.4 (primary HTML anchors \#S3.SS4, \#S4.SS2 and \#S4.SS4): interface, generated-instance selection, comparison and timing.\par\smallskip \url{https://arxiv.org/html/2607.15543v1} &
\textbf{Checked:} 2026-10-04\par\smallskip \textbf{Basis:} retained prior inspection\par\smallskip \textbf{Limits:} Retained inspection in discovery-search/extension/comparative-search.json. All arms use the fastest 100 variational-converging cases of 120 generated per size; small selected instances and free quantum evolution do not support a universal conclusion. \\
\addlinespace[4pt]
\end{longtable}
\endgroup


\end{landscape}
%

\section*{Contextual records and reconciliation notes}

These 5 cited empirical counterpoints have no row in the main evidence tables. Their records remain in the CSV and Table~S4; their bounded interpretations are given here. Theory and method citations without an empirical audit appear in Table~S5 where discovered by the update.

\paragraph{He et al., 2026}

\textbf{Source:} arXiv v1 (2026-03-10); hardware-comparator characterization selectively rechecked on 2026-10-04\par\smallskip \textbf{Design:} Five 96-node 3-regular graphs: hardware RWS-QAOA versus implemented GW and HLZ. Without local search it beats GW and approaches HLZ in approximation ratio; with local-search postprocessing also applied to classical outputs it surpasses HLZ, especially in success probability. Unrestricted strengthened Burer-Monteiro crossover is projected.\par\smallskip \textbf{Scope:} Five 96-node 3-regular Max-Cut graphs for the hardware comparison.\par\smallskip \textbf{Supported inference:} Measured solution-quality comparison with implemented GW and HLZ; unrestricted-system timing advantage remains projected.\par\smallskip \textbf{Timing:} Measured hardware quality comparison; future fault-tolerant runtime crossover against unrestricted strengthened Burer-Monteiro is a projection, not a measured current-hardware speedup.\par\smallskip \textbf{Certificate:} Optimal reference cuts define comparison metrics; the quantum routine does not itself certify optimality.

\paragraph{Tao et al., 2026}

\textbf{Source:} arXiv:2604.26477v2 (May 20, 2026); selective timing-scope reread October 6, 2026\par\smallskip \textbf{Design:} contextual counterpoint: sampling-time comparison with reported quantum hardware results (Figure 1); complete-pipeline comparison with classical DCM and DPA-a (Figure 2); complete quantum-system comparison unresolved\par\smallskip \textbf{Scope:} closely related multi-objective benchmarks (cited counterpoint)\par\smallskip \textbf{Supported inference:} Classical sampling and pipeline counterpoint; complete quantum-system comparison unresolved.\par\smallskip \textbf{Timing:} Figure 1 compares sampling time; Figure 2 includes Ising model construction, sampling and Pareto filtering for the classical quantum-inspired pipeline versus classical baselines; complete quantum-pipeline times are not established\par\smallskip \textbf{Certificate:} Not applicable

\paragraph{Safi et al., 2026}

\textbf{Source:} arXiv:2607.13325v1\par\smallskip \textbf{Design:} Industrial scheduling assignment QUBOs with classical sequencing, varied formulations/platforms; small exact baseline, time-limited references and classical heuristics. No isolated quantum-component replacement.\par\smallskip \textbf{Scope:} Industrial scheduling with alternative compact and constraint-rich assignment formulations.\par\smallskip \textbf{Supported inference:} Limited platform/formulation comparison; quantum component attribution unresolved.\par\smallskip \textbf{Timing:} Paper labels execution end-to-end but explicitly omits communication latency, API overhead, queueing. Shared scaling trend is not matched elapsed time.\par\smallskip \textbf{Certificate:} Quantum assignment/sequence split heuristic; exact branch-and-bound reference only on tractable cases; larger references time-limited.

\paragraph{Bantysh et al., 2026}

\textbf{Source:} arXiv:2609.15738v1\par\smallskip \textbf{Design:} Classical product-state projection after two-qubit updates; approximate MPS analysis and GSet quality versus literature best-known cuts. Not a matched whole-system solver race.\par\smallskip \textbf{Scope:} GSet Max-Cut and small random-graph entanglement studies.\par\smallskip \textbf{Supported inference:} Constructive classical baseline in fixed-parameter expanding-depth regime; no general conclusion about quantum optimization or entanglement.\par\smallskip \textbf{Timing:} Fixed-schedule timing covers projected evolution and excludes graph coloring. Quality is the best cut among 10 runs; reported time is the mean per run. Tuning times identify the first best hit during an approximately one-hour search.\par\smallskip \textbf{Certificate:} Heuristic product-state solver; best-known reference cuts do not certify optimum.

\paragraph{Gaidai and Mukherjee, 2026}

\textbf{Source:} arXiv:2607.15543v1\par\smallskip \textbf{Design:} Ideal variational binary sampler versus Uniform and SMAC, sharing the continuous CasADi/IPOPT subproblem; separate monolithic SCIP comparison at 5 generators.\par\smallskip \textbf{Scope:} Generated AC optimal-power-flow/unit-commitment instances.\par\smallskip \textbf{Supported inference:} Negative component result on small ideal simulations: trained sampling gives no incumbent-quality benefit over uniform. Not a universal no-advantage claim.\par\smallskip \textbf{Timing:} Equal charged budget per size with quantum simulation time removed entirely. This favorable zero-quantum-time assumption is not hardware wall-clock evidence.\par\smallskip \textbf{Certificate:} Heuristic nonconvex inner optimization; no global original-problem certificate established by hybrid sampling.

\paragraph{Reconciliation of the recovered records.}

\textbf{[E3]}\par Dupont et al., 2025: Corrected the historical CSV phrase "pipeline fixed" to distinguish the reported whole-system comparison from a fixed-wrapper component swap, in line with the revised main article.\par

\textbf{[E6]}\par Vodeb et al., 2026: Corrected the below-251 graph-size threshold: QUBO-to-Max-Cut adds one vertex, so bqp250 has 251 graph vertices. Main-table wording uses the certified subset without conflating variables and vertices.\par

\textbf{[E7]}\par Lozano, 2026: Updated empirical source from v2 (June 2) to v3 (July 12), an explicitly identified post-cutoff expansion of an eligible study. The new Tabu control uses 36 instances/108 seeded runs; service optimality anchors and mean QPU access remain unchanged.\par

\textbf{[E9]}\par Pennington et al. (2025): Reconciled supplement availability wording to the source-code repository already recorded in the original CSV.\par

\textbf{[E13]}\par de Valk et al., 2026: New empirical inclusion in the October literature extension; no prior-version record applies because the study was absent from the historical July audit. Bibliographic correction on October 6, 2026: the primary arXiv v1 author metadata gives Sebastian de Bon (accented given name in the bibliography); the bibliography follows the surname de Bon. The earlier October 4 status note about a title-page spelling is preserved as a historical observation, not the current name decision. Source: https://arxiv.org/abs/2607.26771v1. Round 2 reconciliation on October 6, 2026: direct inspection of the same version's PDF title page confirms Sebastian de Bone (accented given name in the source), while the abstract metadata still gives de Bon. The bibliography now follows the PDF spelling de Bone, superseding the earlier metadata-based decision without treating compliance with that earlier review as an unaddressed request. Both surname particles remain protected for journal rendering. Sources: https://arxiv.org/pdf/2607.26771v1 (page 1 author block) and https://arxiv.org/abs/2607.26771v1 (author metadata).\par

\textbf{[E14]}\par Naghmouchi et al., 2026: New empirical inclusion in the October literature extension; no prior-version record applies because the study was absent from the historical July audit.\par

\textbf{[E15]}\par Ramesh et al., 2026: New October inclusion. Survey author David E. Bernal Neira is a coauthor and corresponding author; overlap disclosed in the main article.\par

\textbf{[E16]}\par Yalovetzky et al., 2026: New empirical inclusion in the October literature extension; no prior-version record applies because the study was absent from the historical July audit.\par

\textbf{[E17]}\par Liaqat et al., 2026: New E17 is a separate paper extending Pennington E9, not a new version of E9; E9 historical empirical version remains intact. 1710 is the assembled binary-variable count, not a coherent device width. Preliminary 4096-variable conclusions lack the same detailed accounting and are not used to enlarge the audited scale.\par

\textbf{[E18]}\par Lauwerier et al., 2026: New E18. In Table 3, quantum reward 20.74/28.37/21.25 versus exact 21.18/28.37/27.14 at 10/15/20 vertices. Positive LP-bound gaps for exact solutions are relaxation gaps, not solver suboptimality. Selection uses an unweighted MIS quantum objective with classical weighted recovery; more candidate tours changes the restricted search space as well as qubit count.\par

\textbf{Context}\par He et al., 2026: Corrected the historical characterization as only exceeding worst-case guarantees: the inspected v1 includes head-to-head implemented GW/HLZ comparisons. No new version was needed for this correction.\par

\textbf{Context}\par Safi et al., 2026: New empirical inclusion in the October literature extension; no prior-version record applies because the study was absent from the historical July audit.\par

\textbf{Context}\par Bantysh et al., 2026: New empirical inclusion in the October literature extension; no prior-version record applies because the study was absent from the historical July audit.\par

\textbf{Context}\par Gaidai and Mukherjee, 2026: New empirical inclusion in the October literature extension; no prior-version record applies because the study was absent from the historical July audit.\par


\clearpage
\small
%
\section*{Data dictionary}
\subsection*{Dates and scope}
Expanded evidence eligibility cutoff: October 12, 2026. Publication-status check date: October 12, 2026. Original cutoff: June 30, 2026. Historical audit snapshot: July 27, 2026. The targeted literature extension covers July 1, 2026 through October 12, 2026, with searches assessed in America/New\_York. These dates govern different operations. Original-cutoff statuses and historical source versions are preserved; new inclusions and selectively reread claims record their inspected version. The update includes targeted discovery, publication checks, and empirical inspection of selected sources, not a claim of exhaustive coverage. Source-location verification dates are recorded per claim, distinguishing retained earlier inspections from later selective passage checks; they do not extend the eligibility cutoff or imply a new corpus-wide publication-status check.

\subsection*{Reconciliation and generation}
The 30 study records comprise 25 main evidence-table records and 5 contextual records with row\_id "-". citation\_key is unique even for contextual records. Tables S1 and S2 render the evidence-table audits; Table S4 renders all study version/status records. Table S3 renders 4 historical exclusions. Table S5 renders 82 current discovery decisions. Table S6 indexes the inspected source versions and passages for all 30 studies. The CSVs are the source of truth. Tables S1-S5 and this dictionary are rendered deterministically after checking identifiers, schema coverage, inference-label agreement with the main tables, and publication-status fields. Table S6 is rendered from supplementary/claim-locations.csv after checking complete study coverage. The source download supplies the standalone LaTeX source and ancillary CSV records; the ancillary README.txt explains compilation and packaged filenames, and SOURCE\_MANIFEST.json records original paths and source checksums.

\subsection*{Supported effects}
A component effect compares solvers inside a fixed wrapper; a system effect compares complete pipelines at a stated resource boundary and target quality; a reduction-solver interaction asks whether the reduction changes relative performance. Wrapper benefit and negative or unresolved comparisons are informative findings. Each claim is assessed at its own scope; a single blanket component-resolved end-to-end standard is not applied to every study.

\subsection*{Label conventions}
"Not reported by the authors" asserts an absence only when the source was inspected for that field. "Not applicable" marks a field without meaning for that comparison. "Not found by this survey's documented search" and "not found by this survey's search of the study text" describe search limitations. "Not extracted from the study by this survey's audit" describes an audit limitation. Empty artifact identifiers mean no identifier was captured. Historical artifact labels and source-version descriptions are retained unless a documented correction says otherwise; a current publication check does not reverify them.

\subsection*{Inclusion policy and provenance}
A study enters the evidence synthesis when its comparative design can inform a component, system, or interaction effect, regardless of publication status. Mechanism studies may be discussed without a table row. The exclusion log records the historical purposive selection decisions and their present disposition, not an exhaustive census. The original source package remains in the recovery archive. No archival DOI for this manuscript package has been assigned by this update; a future deposit must be cited only after it exists.

\subsection*{Columns of \texttt{supplementary\_studies.csv}}
\noindent\texttt{row\_id}: Evidence-table identifier E followed by a positive integer or A followed by a positive integer. A hyphen marks contextual empirical records without a main-table row; citation\_key remains unique.\par\smallskip
\noindent\texttt{citation\_key}: BibTeX key in references.bib; the unique record key, including for contextual rows.\par\smallskip
\noindent\texttt{authors\_year}: Short author-year label for the inspected empirical source. A publication-status change does not silently replace that source version.\par\smallskip
\noindent\texttt{venue}: Venue attached to the empirical audit; current\_venue separately records a later verified bibliographic venue.\par\smallskip
\noindent\texttt{version\_inspected}: Empirical source version underlying the current retained claims. A publication-status refresh alone does not replace this version; an explicitly documented selective reread may do so. Unpinned historical labels remain identified as such.\par\smallskip
\noindent\texttt{prior\_version\_inspected}: Exact version\_inspected text preserved from the recovered July 2026 CSV. Empty for studies newly added by the October literature extension.\par\smallskip
\noindent\texttt{status\_at\_original\_cutoff}: Historical publication status at the original June 30, 2026 eligibility cutoff, preserved when the current cutoff expands. Newer studies record that they were not yet public then. Retrospective acceptance dates alone do not establish public availability.\par\smallskip
\noindent\texttt{status\_at\_prior\_verification}: Status recorded in the historical July 27 audit snapshot, formerly status\_at\_submission. New additions are marked as absent from that historical audit rather than reconstructing a check that was not performed.\par\smallskip
\noindent\texttt{identifier}: DOI or arXiv identifier recorded with the historical audit/citation. For post-cutoff journal citations the empirical version is still distinguished in version\_inspected.\par\smallskip
\noindent\texttt{search\_route}: How the study entered this purposive survey, including historical review routes or the targeted literature extension. The latter discovery decisions are recorded in supplementary/literature-search.csv.\par\smallskip
\noindent\texttt{inclusion\_reason}: Reason the study informs the decomposition lens or a component, system, or reduction-solver interaction effect.\par\smallskip
\noindent\texttt{quantum\_role}: Task assigned to the quantum routine, or an explicit absence of a quantum routine in classical reassessments.\par\smallskip
\noindent\texttt{classical\_role}: Tasks assigned to classical routines and controls.\par\smallskip
\noindent\texttt{comparison\_design}: What the study holds fixed and varies. Whole-system comparisons do not become component comparisons merely because a pipeline is described.\par\smallskip
\noindent\texttt{instance\_family}: Problem family, sizes, and counts retained in the audit.\par\smallskip
\noindent\texttt{target\_quality\_definition}: Objective-quality or success target used in the comparison; its definition is needed before comparing resource use.\par\smallskip
\noindent\texttt{timing\_boundary}: Included and excluded resources. Standard terms are full end-to-end wall clock, solver-side pipeline time, device execution time, projected time, and nominal operation or sample count. No timing claim is distinct from any timed comparison.\par\smallskip
\noindent\texttt{uncertainty\_reporting}: Variability or uncertainty reporting captured by the audit. An audit limitation is not evidence that the study lacks uncertainty analysis.\par\smallskip
\noindent\texttt{hardware\_software\_versions}: Hardware, solver, software-version, and parameter information at the resolution actually captured by the audit.\par\smallskip
\noindent\texttt{certificate\_status}: Independently checkable witnesses or certificates, their provider, and limits. The main evidence tables do not contain a certificate column; this field supplies that detail.\par\smallskip
\noindent\texttt{artifact\_status}: Controlled vocabulary: archival dataset; source-code repository; input-data repository; available on request; publisher supplement; not reported; not found by documented search. This identifies the reported artifact class, not an independent audit of its contents; artifact\_note retains scope and access limits.\par\smallskip
\noindent\texttt{artifact\_identifier}: Known DOI or repository reference. Empty means no identifier was captured; never substitute an invented archival DOI or a promise of a future record.\par\smallskip
\noindent\texttt{supported\_inference}: Bounded interpretation aligned with the main evidence tables: wrapper benefit, component comparison, system comparison, interaction evidence, or unresolved attribution, with direction and qualifications where supported. These effects are assessed separately.\par\smallskip
\noindent\texttt{scale\_and\_repetitions}: Scale and repetition details supporting the audited comparison, recovered from the original supplement or inspected for a newly included study.\par\smallskip
\noindent\texttt{baseline\_hardware\_and\_stopping\_rule}: Baseline implementation and stopping-rule details captured by the empirical audit; uncaptured details remain explicitly unspecified.\par\smallskip
\noindent\texttt{artifact\_note}: Descriptive availability statement or explicit reconciliation note. artifact\_status and artifact\_identifier govern the displayed availability classification.\par\smallskip
\noindent\texttt{audit\_snapshot\_date}: Date of the audit snapshot when this record was created, not a claim that every empirical source was independently reread on that exact date.\par\smallskip
\noindent\texttt{status\_current}: Publication status supported by the current check, or an explicit unresolved/not-reverified description. It may differ from the empirical source version.\par\smallskip
\noindent\texttt{status\_check\_date}: ISO date of the current publication-status check or attempted check. A date alone does not mean verification succeeded.\par\smallskip
\noindent\texttt{status\_check\_result}: Controlled result: verified\_primary\_record; limited\_primary\_record; unresolved. The result concerns publication metadata, not an independent verification of all empirical claims.\par\smallskip
\noindent\texttt{status\_sources}: Primary-record URLs supporting the check, separated by " | ". Inaccessible or inconclusive sources are identified in status\_notes, not treated as successful verification.\par\smallskip
\noindent\texttt{status\_notes}: Evidence, limitations, version availability, or explanation of the current status. No-journal-version-found is search-limited, not proof that no such version exists.\par\smallskip
\noindent\texttt{current\_identifier}: Identifier of the publication record checked in the current refresh, when supported. Historical identifier and empirical version fields are retained separately.\par\smallskip
\noindent\texttt{current\_venue}: Venue of the publication record checked in the current refresh, when supported.\par\smallskip
\noindent\texttt{latest\_version\_at\_status\_check}: Latest source version observed during the metadata check, when captured. Availability of a newer version does not make it the empirical source inspected.\par\smallskip
\noindent\texttt{empirical\_refresh\_date}: ISO date of an explicitly documented primary-text inspection for a new empirical inclusion or a selective reread; empty for metadata-only refreshes.\par\smallskip
\noindent\texttt{empirical\_refresh\_note}: Claims inspected, source URLs, corrections or new findings, and limitations. Scope is local to the stated claims, not a comprehensive empirical re-audit.\par\smallskip
\noindent\texttt{reconciliation\_notes}: Corrections of drift between the historical CSV, supplement, and revised main article. The recovered original package remains available as provenance.\par\smallskip
\subsection*{Columns of \texttt{supplementary\_exclusion\_log.csv}}
\noindent\texttt{study}: Authors and short title of the close but omitted study.\par\smallskip
\noindent\texttt{identifier}: Persistent identifier retained in the exclusion log.\par\smallskip
\noindent\texttt{first\_version\_date}: Historical first-version date used for eligibility; retain the recorded precision instead of inventing an exact day.\par\smallskip
\noindent\texttt{why\_considered}: Recorded reason the study was close to the survey scope.\par\smallskip
\noindent\texttt{why\_excluded}: Historical reason for omission from the evidence tables. This describes the purposive selection decision; absence of hardware, matched timing, or a component-isolating control is not a universal bar to inclusion of system or interaction evidence.\par\smallskip
\noindent\texttt{exclusion\_id}: Stable historical exclusion-log identifier X followed by a positive integer, independent of evidence-table identifiers.\par\smallskip
\noindent\texttt{citation\_key}: Main-bibliography key when the excluded study is cited there; empty when no corresponding citation is used.\par\smallskip
\noindent\texttt{status\_current}: Publication status supported by the current check, or an explicit unresolved statement.\par\smallskip
\noindent\texttt{status\_check\_date}: ISO date of the publication-status check or attempted check.\par\smallskip
\noindent\texttt{status\_check\_result}: verified\_primary\_record, limited\_primary\_record, or unresolved; metadata scope only.\par\smallskip
\noindent\texttt{status\_sources}: Supporting primary-record URLs separated by " | ".\par\smallskip
\noindent\texttt{status\_notes}: Evidence and limitations of the check, without treating unavailable material as verified.\par\smallskip
\noindent\texttt{review\_disposition}: Whether the historical exclusion was retained or changed in this revision, and the scope of that decision.\par\smallskip
\subsection*{Columns of \texttt{supplementary/literature-search.csv}}
\noindent\texttt{candidate\_id}: Unique stable identifier for a candidate considered in the targeted literature extension.\par\smallskip
\noindent\texttt{citation\_key}: Main-bibliography key for included evidence or context; may be blank for screened-only or out-of-scope candidates.\par\smallskip
\noindent\texttt{title}: Title of the candidate source used in the discovery decision.\par\smallskip
\noindent\texttt{source\_url}: Primary record or full-text URL supporting candidate identity and the stated decision.\par\smallskip
\noindent\texttt{first\_public\_date}: First public date verified for the candidate, ISO YYYY-MM-DD; blank if not verified. Submission and announcement dates can differ and are explained in reason when material.\par\smallskip
\noindent\texttt{version\_inspected}: Exact source version inspected, with any access limit stated explicitly. A current record alone does not imply full-text inspection.\par\smallskip
\noindent\texttt{version\_date}: Date of the inspected version, ISO YYYY-MM-DD; blank if not verified.\par\smallskip
\noindent\texttt{discovery\_route}: Query family, citation-chaining route, existing-record refresh, or other documented route by which the candidate was considered.\par\smallskip
\noindent\texttt{decision}: Controlled vocabulary: include\_evidence; include\_context; screen\_only; exclude\_scope. This is a purposive editorial decision, not a verdict on correctness.\par\smallskip
\noindent\texttt{placement}: Main evidence row, manuscript discussion, or screen-only location; blank when no placement applies.\par\smallskip
\noindent\texttt{reason}: Relevance, supported inference, limitations, and reason for the decision. Distinguish source findings from audit judgments and date ambiguity.\par\smallskip
\noindent\texttt{source\_check\_date}: ISO date of the discovery check, not a claim of universal literature completeness or a full proof/empirical verification.\par\smallskip
\subsection*{Columns of \texttt{supplementary/claim-locations.csv}}
\noindent\texttt{row\_id}: Evidence-table identifier E1-E18 or A1-A7, matching supplementary\_studies.csv. A hyphen marks a contextual record; citation\_key remains unique.\par\smallskip
\noindent\texttt{citation\_key}: Unique BibTeX/study key linking this source-location record to supplementary\_studies.csv and the main manuscript.\par\smallskip
\noindent\texttt{claim}: Bounded numerical claim or comparison to locate, including its relevant comparator, resource boundary or scope qualification.\par\smallskip
\noindent\texttt{source\_version}: Exact version used for the stated passage check. An author preprint used for detailed locators is distinguished from the cited journal record; this field does not silently replace the historical audit version.\par\smallskip
\noindent\texttt{source\_url}: HTTPS primary-source URL for the inspected record or text. Multiple sources, when needed, are separated by " | ".\par\smallskip
\noindent\texttt{locator}: Specific abstract passage, section, table, figure, appendix or stable source anchor supporting the claim in source\_version.\par\smallskip
\noindent\texttt{verification\_date}: ISO date of the supporting inspection or retained documented check. This claim-level date is distinct from the corpus eligibility cutoff, publication-status date and the date the index was assembled.\par\smallskip
\noindent\texttt{verification\_basis}: selective\_primary\_passage\_check identifies a targeted primary-text inspection on verification\_date; retained\_prior\_inspection identifies an earlier documented inspection carried forward without a new reread. new\_study\_inspection identifies a full-text extraction for a newly included study. None of these labels implies experimental reproduction.\par\smallskip
\noindent\texttt{limitations}: Access restrictions, version differences, retained-audit boundaries, resource-accounting qualifications and other limits on what the cited passage establishes.\par\smallskip


\clearpage
\begin{landscape}
%
\section{Software implementations of decomposition and reduction}
Table~\ref{tab:s7} records 23 public implementations inspected on October 9, 2026.
This purposive software inventory is separate from the empirical audit and does not extend its eligibility cutoff.
Inclusion requires a concrete implementation of subproblem optimization, reconstructable reduction, multilevel optimization, structural preprocessing, or a constrained outer/inner interface.
Packages, research prototypes, and an educational implementation are identified separately; multiple algorithms within one package are not counted as separate tools.
Source and documentation were inspected, but no package was installed, executed, independently validated, or benchmarked for this inventory.
The short commit links identify inspected implementation files, not necessarily released versions; full hashes, additional source paths, and evidence notes are in \nolinkurl{supplementary/software-catalog.json}.
\texttt{QUBODecomposition.jl} version 0.1.0 was separately verified as registered and released on October 12, 2026; this does not change the implementation-inspection dates.
The accompanying \nolinkurl{supplementary/software-screening.json} records selected leads outside the table; no verified public source in this search does not imply that none exists.
For example, \texttt{Quafu-Qcover} decomposes QAOA expectation evaluation into local neighborhoods; this is distinct from optimizing independent sub-QUBOs and merging their assignments \citep{xu2024qcover}.
License labels transcribe the inspected project declarations rather than assess legal compatibility; a missing project license is reported as unidentified.
\small
\setlength{\tabcolsep}{5pt}
\setlength{\AuditTableWidth}{\dimexpr\linewidth-8\tabcolsep\relax}
\begin{longtable}{@{}>{\RaggedRight\arraybackslash}p{.19\AuditTableWidth}>{\RaggedRight\arraybackslash}p{.33\AuditTableWidth}>{\RaggedRight\arraybackslash}p{.30\AuditTableWidth}>{\RaggedRight\arraybackslash}p{.18\AuditTableWidth}@{}}
\caption{Public software implementing decomposition and reduction. Source inspection establishes scope, not comparative performance or certified software correctness.}\label{tab:s7}\\
\toprule Software and artifact type & Implemented mechanism & Scope and guarantee limits & Inspected source / declared license \\
\midrule\endfirsthead
\multicolumn{4}{l}{\textit{Table S7 continued}}\\
\toprule Software and artifact type & Implemented mechanism & Scope and guarantee limits & Inspected source / declared license \\
\midrule\endhead
\midrule\multicolumn{4}{r}{\textit{Continued on next page}}\\\endfoot
\bottomrule\endlastfoot
\texttt{qbsolv}\par Archived solver & Alternates bounded, incumbent-conditioned sub-QUBO optimization with tabu search. The default subsolver is classical; an annealer can be substituted. \citep{booth2017partitioning} & Heuristic for coupled problems; exact subsolves alone do not certify the original optimum. Archived dependencies may require reconstruction. & \href{https://github.com/dwavesystems/qbsolv/blob/7b71f110ad9837d1aa4bd80efaaad692feda810f/src/solver.cc}{\texttt{7b71f11}}\par Apache-2.0 \\
\addlinespace[5pt]
\texttt{dwave-hybrid}\par Python framework & Composable decomposers, samplers and composers implement energy-impact, graph-traversal and random neighborhoods, conditioned subproblems, and parallel solver portfolios. \citep{dwave_hybrid_framework_docs} & A workflow framework, not a global certificate. Kerberos and LatticeLNLS workflows for geometrically structured interaction graphs are templates within this package, not independent software entries. & \href{https://github.com/dwavesystems/dwave-hybrid/blob/ec17a700b0250123da9909ec82db4ecb2516993d/hybrid/decomposers.py}{\texttt{ec17a70}}\par Apache-2.0 \\
\addlinespace[5pt]
\texttt{QSplit}\par Workflow framework & Matrix, interaction-neighborhood, graph-partition and macrovariable splitting feed distributed solves and voting, with optional conditioned, mean-field or soft-consensus refinement. \citep{bifulco2027qsplit,qsplit_software} & Aggregation and compression are heuristic. Matrix-tile dimensions need not equal logical variable counts; reconstruction must be evaluated on the original QUBO. & \href{https://github.com/alpha-unito/QSplit/blob/279aff6e1b1a499c36120cd0dbd3e85e57f83080/cpp/splitting.cpp}{\texttt{279aff6}}\par MIT \\
\addlinespace[5pt]
\texttt{Pasqal QUBO Solver}\par Python solver package & DecomposeQuboSolver builds device-layout-compatible subproblems and delegates remaining variables to a classical solver. \citep{pasqal_qubo_decomposition_docs} & The inspected decomposition rejects negative off-diagonal coefficients. Layout-compatible local solutions do not supply a general global certificate. & \href{https://github.com/pasqal-io/qubo-solver/blob/5211de10a7f2db2e3f8b3d2349d62385b9d94566/qubosolver/solver.py}{\texttt{5211de1}}\par Pasqal MIT-derived agreement (research-limited patent grant) \\
\addlinespace[5pt]
\texttt{QUBODecomposition.jl}\par Julia package & Solves fitting connected components and otherwise performs bounded conditioned neighborhood sweeps through optimizer interfaces, with full-model objective checks and budgets. \citep{qubodecomposition_software} & Component-wise exactness requires exact child solves; coupled sweeps remain heuristic. The manuscript author contributes to this package. & \href{https://github.com/JuliaQUBO/QUBODecomposition.jl/blob/c760d09594adf2dcfc6884523e908d8076af6b4c/src/serial.jl}{\texttt{c760d09}}\par MPL-2.0 \\
\addlinespace[5pt]
\texttt{HADOF}\par Public research prototype & Iteratively forms sub-QUBOs using sampled marginal means for external variables, updates those means with sampling subsolvers, and combines final block samples. Sequential and parallel implementations support annealing and variational/circuit backends. \citep{Sankar2025HADOF,Sankar2026ParallelHADOF} & Mean-field boundary substitution and final sample aggregation are heuristic; they do not certify the original QUBO optimum or preserve its full joint distribution. & \href{https://github.com/Namasi-PhD/HADOF/blob/ced7ee409ee0571362cf7628928765ab1fe26b94/HADOF/sequentialHADOF.py}{\texttt{ced7ee4}}\par No project license identified \\
\addlinespace[5pt]
\texttt{DQAOA (ORNL research implementation)}\par Public research prototype & Solves overlapping then randomly selected principal sub-QUBOs with QAOA, distributes subsolves with MPI, and incorporates proposed bits only when the full QUBO objective improves; a single-core variant is also provided. \citep{Kim2026DistributedQAOA} & Principal submatrices omit fixed-exterior linear contributions; global objective checks safeguard accepted bits, but neither heuristic proposals nor MPI parallelism establish optimality or quantum advantage. & \href{https://github.com/Seongmin-Kim-ornl/DQAOA/blob/9ed1842760fee1f20e2e25c8eabd97c841a51233/DQAOA/DQAOA.py}{\texttt{9ed1842}}\par No project license identified \\
\addlinespace[5pt]
\texttt{QAOA-in-QAOA (QAOA squared)}\par Public research implementation & Partitions a weighted MaxCut graph, solves each small subgraph, then constructs a coarser graph whose binary variables choose whether to flip entire subgraph assignments. Repeats this contraction until the remaining graph fits the subsolver. \citep{Zhou2023QAOAinQAOA} & Merging exploits MaxCut complement symmetry and optimizes over retained subgraph assignments; it is not unrestricted exact QUBO decomposition. Paper approximation guarantees require their stated assumptions. & \href{https://github.com/ZeddTheGoat/QAOA_in_QAQA/blob/7704cfd2c2cbfac58a11fbf4f8beebb1efb9c04c/QAOA_square.py}{\texttt{7704cfd}}\par Apache-2.0 \\
\addlinespace[5pt]
\texttt{Neural QAOA squared}\par Research code and model weights & Extends QAOA-in-QAOA with a learned, capacity-constrained graph partition and joint QAOA parameter initialization, followed by subgraph solution and hierarchical orientation merging. Includes classical partition policies, training code, evaluation drivers, and pretrained weights. \citep{Zheng2026NeuralQAOA} & A learned heuristic: training, inference, and subsolver costs belong in comparisons. Its evaluated instance results do not establish general superiority or quantum advantage. & \href{https://github.com/0SliverBullet/Neural-QAOA-Squared/blob/fcbc0ab3adca5afc35d83d6f363dadc1512cfbea/src/models/generator.py}{\texttt{fcbc0ab}}\par MIT \\
\addlinespace[5pt]
\texttt{CUDA-Q academic divide-and-conquer QAOA}\par Executable educational example & Recursively partitions MaxCut graphs, solves subgraphs with CUDA-Q QAOA, and solves a merger Hamiltonian that chooses whole-subgraph color flips; MPI distributes initial subgraph workloads across workers. \citep{NVIDIA2026DCQAOA} & Teaching implementation with heuristic subsolves and aggregation; it is not a generic packaged QUBO decomposer, nor does distribution certify solution quality. & \href{https://github.com/NVIDIA/cuda-q-academic/blob/847ec2de0aa81bbe2a0d9d8be5a2d11928fabc4f/qaoa-for-max-cut/qaoa-divide-and-conquer.py}{\texttt{847ec2d}}\par Apache-2.0 AND CC-BY-NC-4.0 (source header) \\
\addlinespace[5pt]
\texttt{ML-QLS (ml\_qls)}\par research implementation & Coarsens graph-partitioning and community-detection instances, solves a coarse model, projects the solution, and refines bounded local Ising subproblems through classical, annealing, or QAOA interfaces. \citep{ushijima2021multilevel} & Multilevel heuristic; local solves do not certify the original optimum. Historical dependencies require environment reconstruction. & \href{https://github.com/rsln-s/ml_qls/blob/bc376276ba684460aeccaa371b4fc38003139e34/multilevel/GP/ml_gp.py}{\texttt{bc37627}}\par No project license identified \\
\addlinespace[5pt]
\texttt{MLQAOA}\par Research code and notebooks & Coarsens Max-Cut graphs, solves the coarse problem, projects assignments, and refines small subproblems; graph representations select transferred QAOA parameters from stored reference instances. \citep{bach2024mlqaoa} & Heuristic Max-Cut workflow; spectral similarity and approximate local solves provide no global certificate. Code uses legacy simulator calls. & \href{https://github.com/bachbao/MLQAOA/blob/41ba8b8b31185e620254293db9ec9fb3da79b280/maxcut.py}{\texttt{41ba8b8}}\par No project license identified \\
\addlinespace[5pt]
\texttt{OpenQAOA}\par Released Python SDK & RQAOA ranks one- and two-spin expectation values, fixes spins or pair relations, rewrites the Ising problem, solves a small residual classically, and reconstructs complete solutions. Custom schedules and adaptive multi-elimination are implemented. \citep{sharma2022openqaoa,bravyi2020obstacles} & Correlation fixing is heuristic; exact solution of the residual does not certify the original optimum. Sampling accuracy and reduction choices matter. & \href{https://github.com/entropicalabs/openqaoa/blob/c6e7dd759392fcca180bab9bd29c5f5e7dab3662/src/openqaoa-core/openqaoa/algorithms/rqaoa/rqaoa_workflow.py}{\texttt{c6e7dd7}}\par MIT \\
\addlinespace[5pt]
\texttt{QIRO reference implementation}\par Public research prototype & Quantum correlations guide problem-specific MIS and MAX-2-SAT simplification. The MIS implementation solves small disconnected components by enumeration; the MAX-2-SAT implementation offers backtracking. A separate RQAOA implementation handles quadratic Hamiltonians. \citep{finzgar2024qiro} & Correlation-based choices can remove an optimum. Feasibility-preserving rules and local enumeration do not provide a global optimality certificate; optional backtracking is bounded heuristic search. & \href{https://github.com/jernejrudifinzgar/qiro/blob/8fc8e9564c8694a1e21583617ed424ec4ff7bf4d/QIRO_MIS.py}{\texttt{8fc8e95}}\par No project license identified \\
\addlinespace[5pt]
\texttt{Qrisp QIRO}\par Python framework module & QIROProblem alternates QAOA sampling with user-supplied replacement routines and rebuilds cost operators, mixers and initialization. Supplied routines cover MIS, clique and MaxSAT; set packing is transformed to an MIS graph. \citep{qrisp2026qiro,finzgar2024qiro} & Heuristic, problem-specific reductions require separate validation. The inspected workflow keeps the original register length, so logical simplification does not itself establish reduced hardware width. & \href{https://github.com/eclipse-qrisp/Qrisp/blob/30a266c461a195375bbf8cb1e12aa7ea53256261/src/qrisp/algorithms/qiro/qiro_problem.py}{\texttt{30a266c}}\par EPL-2.0; conditional secondary GPL license \\
\addlinespace[5pt]
\texttt{qReduMIS}\par Python research package & Alternates exposed-clique MIS reductions with decisions inferred from sampled independent sets. Inset/outset strategies select likely included/excluded vertices, maintain an incumbent, and restart classical reduction on the remaining graph. \citep{schuetz2026qredumis} & Sample-guided fixing is heuristic despite exact classical rules. The inspected solver accepts atom-position MIS instances and implements Rydberg backends, not generic QUBO or weighted MIS. & \href{https://github.com/jpmorganchase/qReduMIS/blob/43feb7502ca4c90e10bde24e077c1d54427f8fef/qReduMIS_package/qReduMIS/mis_solver.py}{\texttt{43feb75}}\par Apache-2.0 \\
\addlinespace[5pt]
\texttt{DASQA (DBK/DBR)}\par research implementation & Branches maximum-clique and minimum-vertex-cover graphs, applies problem-specific reductions and bounds, and delegates size-limited leaves to a user-supplied solver before reconstructing a solution. \citep{pelofske2019maxclique,pelofske2019vertexcover,pelofske2021decomposition} & Exactness is conditional on exact leaf solves and valid pruning; substituting a heuristic annealer does not supply an optimality certificate. & \href{https://github.com/lanl/Decomposition-Algorithms-for-Scalable-Quantum-Annealing/blob/78f010906a6234a2e23a4da0015429162680a0b6/DBK.py}{\texttt{78f0109}}\par No project license identified \\
\addlinespace[5pt]
\texttt{dwave-preprocessing}\par Python/C++ utility package & Roof duality computes a lower bound and persistent assignments; fixing composites reduce and restore variables, while a component composite dispatches disconnected BQM components to child samplers. \citep{dwave_preprocessing_docs} & Strong persistencies hold in every optimum; weak choices preserve at least one. Global component optimality still requires exact child solves. & \href{https://github.com/dwavesystems/dwave-preprocessing/blob/5e653e1624ce2f1b5ad9d2508abd25a4c3ca8b08/dwave/preprocessing/lower_bounds.py}{\texttt{5e653e1}}\par Apache-2.0 \\
\addlinespace[5pt]
\texttt{Qoolchain}\par Python/C++ research toolchain & Combines persistency, probing and implication-graph reductions with strongly connected components and Shannon branching; exposes residual models, variable maps, persistent assignments and a lower bound. \citep{orlandi2025qoolchain} & Inspected C++ interface accepts integer coefficients. Shannon branching can grow exponentially; preprocessing and a lower bound alone do not solve all residuals. & \href{https://github.com/ilRenato/Qoolchain/blob/537490c948bd3c952e9ff1a85b924ba8de9860b0/Toolchain_implementation/PreprocessingToolchain.pyx}{\texttt{537490c}}\par MIT \\
\addlinespace[5pt]
\texttt{QUBOLite}\par Python utility package & Implements QPRO+ implication/persistency rules and partial assignments that apply fixings to an upper-triangular QUBO, retain an objective offset, and expand residual solutions. \citep{mucke2025qubolite} & Reduction utility, not a hardware decomposition engine. Residual optimization remains necessary; coefficient-compression utilities are a separate transformation. & \href{https://github.com/smuecke/qubolite/blob/cbe57e0bea68e9f751f8591e3ba0254030ac5965/qubolite/preprocessing.py}{\texttt{cbe57e0}}\par No project license identified \\
\addlinespace[5pt]
\texttt{Qiskit Optimization}\par Python optimization library & ADMMOptimizer alternates binary QUBO and continuous solves. RecursiveMinimumEigenOptimizer separately eliminates correlation-selected variable relations and reconstructs solutions. \citep{qiskit_optimization_docs,gambella2020multiblock} & ADMM rejects binary-continuous objective coupling and requires compatible subsolvers; convergence is conditional. Correlation elimination is heuristic, not certified persistency. & \href{https://github.com/qiskit-community/qiskit-optimization/blob/65ce191a4c3ed48f203d8dd6eea476d0deea9e71/qiskit_optimization/algorithms/admm_optimizer.py}{\texttt{65ce191}}\par Apache-2.0 \\
\addlinespace[5pt]
\texttt{HQC-Bend (HQCMCBD)}\par Python research solver & Separates a binary master from continuous linear subproblems, generates one or multiple Benders cuts, and offers Gurobi, D-Wave CQM/BQM, OpenJij, and simulated-bifurcation master backends. \citep{zhao2025hqcbend} & Heuristic master samples and finite encodings do not themselves certify MILP optimality; Gurobi and hardware access are external dependencies. & \href{https://github.com/djzts/HQCMCBD-API/blob/43e9d2951a16da67c769419d72b596c5c042a1a5/HQCMCBD.py}{\texttt{43e9d29}}\par MIT \\
\addlinespace[5pt]
\texttt{Logistiqs (hqc-mcms-benders)}\par paper reproduction package & Keeps Benders master and dual subproblems classical while selecting candidate cuts through minimum-set-cover or maximum-coverage models solved by classical or QAOA backends. \citep{devalk2026hybrid} & Cut selection is a heuristic acceleration interface; MILP certificates depend on classical bounds and valid cuts. Some archived hardware data require paid services to regenerate. & \href{https://github.com/camilledevalk/hqc-mcms-benders/blob/c4fb78901b1b5ffe97143ab10a046a8a3cff950c/milp_engine/solvers/bender_milp_solver.py}{\texttt{c4fb789}}\par No project license identified \\
\addlinespace[5pt]
\end{longtable}
\normalsize


\end{landscape}

%
